\documentclass[reqno,12pt,a4paper]{amsart}

\usepackage{mathrsfs}
\usepackage{amssymb}
\usepackage{amsfonts}
\usepackage{latexsym}
\usepackage{amsthm}
\usepackage{graphicx}
\def\lb{\label}

\newcommand{\er}[1]{\textrm{(\ref{#1})}}

\begin{document}


\renewcommand{\theequation}{\arabic{section}.\arabic{equation}}
\theoremstyle{plain}
\newtheorem{theorem}{\bf Theorem}[section]
\newtheorem{lemma}[theorem]{\bf Lemma}
\newtheorem{corollary}[theorem]{\bf Corollary}
\newtheorem{proposition}[theorem]{\bf Proposition}
\newtheorem{definition}[theorem]{\bf Definition}
\newtheorem{remark}[theorem]{\bf Remark}

\def\a{\alpha}  \def\cA{{\mathcal A}}     \def\bA{{\bf A}}  \def\mA{{\mathscr A}}
\def\b{\beta}   \def\cB{{\mathcal B}}     \def\bB{{\bf B}}  \def\mB{{\mathscr B}}
\def\g{\gamma}  \def\cC{{\mathcal C}}     \def\bC{{\bf C}}  \def\mC{{\mathscr C}}
\def\G{\Gamma}  \def\cD{{\mathcal D}}     \def\bD{{\bf D}}  \def\mD{{\mathscr D}}
\def\d{\delta}  \def\cE{{\mathcal E}}     \def\bE{{\bf E}}  \def\mE{{\mathscr E}}
\def\D{\Delta}  \def\cF{{\mathcal F}}     \def\bF{{\bf F}}  \def\mF{{\mathscr F}}
\def\c{\chi}    \def\cG{{\mathcal G}}     \def\bG{{\bf G}}  \def\mG{{\mathscr G}}
\def\z{\zeta}   \def\cH{{\mathcal H}}     \def\bH{{\bf H}}  \def\mH{{\mathscr H}}
\def\e{\eta}    \def\cI{{\mathcal I}}     \def\bI{{\bf I}}  \def\mI{{\mathscr I}}
\def\p{\psi}    \def\cJ{{\mathcal J}}     \def\bJ{{\bf J}}  \def\mJ{{\mathscr J}}
\def\vT{\Theta} \def\cK{{\mathcal K}}     \def\bK{{\bf K}}  \def\mK{{\mathscr K}}
\def\k{\kappa}  \def\cL{{\mathcal L}}     \def\bL{{\bf L}}  \def\mL{{\mathscr L}}
\def\l{\lambda} \def\cM{{\mathcal M}}     \def\bM{{\bf M}}  \def\mM{{\mathscr M}}
\def\L{\Lambda} \def\cN{{\mathcal N}}     \def\bN{{\bf N}}  \def\mN{{\mathscr N}}
\def\m{\mu}     \def\cO{{\mathcal O}}     \def\bO{{\bf O}}  \def\mO{{\mathscr O}}
\def\n{\nu}     \def\cP{{\mathcal P}}     \def\bP{{\bf P}}  \def\mP{{\mathscr P}}
\def\r{\rho}    \def\cQ{{\mathcal Q}}     \def\bQ{{\bf Q}}  \def\mQ{{\mathscr Q}}
\def\s{\sigma}  \def\cR{{\mathcal R}}     \def\bR{{\bf R}}  \def\mR{{\mathscr R}}
\def\S{\Sigma}  \def\cS{{\mathcal S}}     \def\bS{{\bf S}}  \def\mS{{\mathscr S}}
\def\t{\tau}    \def\cT{{\mathcal T}}     \def\bT{{\bf T}}  \def\mT{{\mathscr T}}
\def\f{\phi}    \def\cU{{\mathcal U}}     \def\bU{{\bf U}}  \def\mU{{\mathscr U}}
\def\F{\Phi}    \def\cV{{\mathcal V}}     \def\bV{{\bf V}}  \def\mV{{\mathscr V}}
\def\P{\Psi}    \def\cW{{\mathcal W}}     \def\bW{{\bf W}}  \def\mW{{\mathscr W}}
\def\o{\omega}  \def\cX{{\mathcal X}}     \def\bX{{\bf X}}  \def\mX{{\mathscr X}}
\def\x{\xi}     \def\cY{{\mathcal Y}}     \def\bY{{\bf Y}}  \def\mY{{\mathscr Y}}
\def\X{\Xi}     \def\cZ{{\mathcal Z}}     \def\bZ{{\bf Z}}  \def\mZ{{\mathscr Z}}
\def\be{{\bf e}}
\def\bv{{\bf v}} \def\bu{{\bf u}}
\def\Om{\Omega}

\newcommand{\mc}{\mathscr {c}}

\newcommand{\gA}{\mathfrak{A}}          \newcommand{\ga}{\mathfrak{a}}
\newcommand{\gB}{\mathfrak{B}}          \newcommand{\gb}{\mathfrak{b}}
\newcommand{\gC}{\mathfrak{C}}          \newcommand{\gc}{\mathfrak{c}}
\newcommand{\gD}{\mathfrak{D}}          \newcommand{\gd}{\mathfrak{d}}
\newcommand{\gE}{\mathfrak{E}}
\newcommand{\gF}{\mathfrak{F}}           \newcommand{\gf}{\mathfrak{f}}
\newcommand{\gG}{\mathfrak{G}}           
\newcommand{\gH}{\mathfrak{H}}           \newcommand{\gh}{\mathfrak{h}}
\newcommand{\gI}{\mathfrak{I}}           \newcommand{\gi}{\mathfrak{i}}
\newcommand{\gJ}{\mathfrak{J}}           \newcommand{\gj}{\mathfrak{j}}
\newcommand{\gK}{\mathfrak{K}}            \newcommand{\gk}{\mathfrak{k}}
\newcommand{\gL}{\mathfrak{L}}            \newcommand{\gl}{\mathfrak{l}}
\newcommand{\gM}{\mathfrak{M}}            \newcommand{\gm}{\mathfrak{m}}
\newcommand{\gN}{\mathfrak{N}}            \newcommand{\gn}{\mathfrak{n}}
\newcommand{\gO}{\mathfrak{O}}
\newcommand{\gP}{\mathfrak{P}}             \newcommand{\gp}{\mathfrak{p}}
\newcommand{\gQ}{\mathfrak{Q}}             \newcommand{\gq}{\mathfrak{q}}
\newcommand{\gR}{\mathfrak{R}}             \newcommand{\gr}{\mathfrak{r}}
\newcommand{\gS}{\mathfrak{S}}              \newcommand{\gs}{\mathfrak{s}}
\newcommand{\gT}{\mathfrak{T}}             \newcommand{\gt}{\mathfrak{t}}
\newcommand{\gU}{\mathfrak{U}}             \newcommand{\gu}{\mathfrak{u}}
\newcommand{\gV}{\mathfrak{V}}             \newcommand{\gv}{\mathfrak{v}}
\newcommand{\gW}{\mathfrak{W}}             \newcommand{\gw}{\mathfrak{w}}
\newcommand{\gX}{\mathfrak{X}}               \newcommand{\gx}{\mathfrak{x}}
\newcommand{\gY}{\mathfrak{Y}}              \newcommand{\gy}{\mathfrak{y}}
\newcommand{\gZ}{\mathfrak{Z}}             \newcommand{\gz}{\mathfrak{z}}

\def\ve{\varepsilon}   \def\vt{\vartheta}    \def\vp{\varphi}    \def\vk{\varkappa}

\def\A{{\mathbb A}} \def\B{{\mathbb B}} \def\C{{\mathbb C}}
\def\dD{{\mathbb D}} \def\E{{\mathbb E}} \def\dF{{\mathbb F}} \def\dG{{\mathbb G}} \def\H{{\mathbb H}}\def\I{{\mathbb I}} \def\J{{\mathbb J}} \def\K{{\mathbb K}} \def\dL{{\mathbb L}}\def\M{{\mathbb M}} \def\N{{\mathbb N}} \def\O{{\mathbb O}} \def\dP{{\mathbb P}} \def\R{{\mathbb R}}\def\S{{\mathbb S}} \def\T{{\mathbb T}} \def\U{{\mathbb U}} \def\V{{\mathbb V}}\def\W{{\mathbb W}} \def\X{{\mathbb X}} \def\Y{{\mathbb Y}} \def\Z{{\mathbb Z}}


\def\la{\leftarrow}              \def\ra{\rightarrow}            \def\Ra{\Rightarrow}
\def\ua{\uparrow}                \def\da{\downarrow}
\def\lra{\leftrightarrow}        \def\Lra{\Leftrightarrow}


\def\lt{\biggl}                  \def\rt{\biggr}
\def\ol{\overline}               \def\wt{\widetilde}
\def\no{\noindent}


\let\ge\geqslant                 \let\le\leqslant
\def\lan{\langle}                \def\ran{\rangle}
\def\/{\over}                    \def\iy{\infty}
\def\sm{\setminus}               \def\es{\emptyset}
\def\ss{\subset}                 \def\ts{\times}
\def\pa{\partial}                \def\os{\oplus}
\def\om{\ominus}                 \def\ev{\equiv}
\def\iint{\int\!\!\!\int}        \def\iintt{\mathop{\int\!\!\int\!\!\dots\!\!\int}\limits}
\def\el2{\ell^{\,2}}             \def\1{1\!\!1}
\def\sh{\sharp}
\def\wh{\widehat}
\def\bs{\backslash}
\def\intl{\int\limits}

\def\na{\mathop{\mathrm{\nabla}}\nolimits}
\def\sh{\mathop{\mathrm{sh}}\nolimits}
\def\ch{\mathop{\mathrm{ch}}\nolimits}
\def\where{\mathop{\mathrm{where}}\nolimits}
\def\all{\mathop{\mathrm{all}}\nolimits}
\def\as{\mathop{\mathrm{as}}\nolimits}
\def\Area{\mathop{\mathrm{Area}}\nolimits}
\def\arg{\mathop{\mathrm{arg}}\nolimits}
\def\const{\mathop{\mathrm{const}}\nolimits}
\def\det{\mathop{\mathrm{det}}\nolimits}
\def\diag{\mathop{\mathrm{diag}}\nolimits}
\def\diam{\mathop{\mathrm{diam}}\nolimits}
\def\dim{\mathop{\mathrm{dim}}\nolimits}
\def\dist{\mathop{\mathrm{dist}}\nolimits}
\def\Im{\mathop{\mathrm{Im}}\nolimits}
\def\Iso{\mathop{\mathrm{Iso}}\nolimits}
\def\Ker{\mathop{\mathrm{Ker}}\nolimits}
\def\Lip{\mathop{\mathrm{Lip}}\nolimits}
\def\rank{\mathop{\mathrm{rank}}\limits}
\def\Ran{\mathop{\mathrm{Ran}}\nolimits}
\def\Re{\mathop{\mathrm{Re}}\nolimits}
\def\Res{\mathop{\mathrm{Res}}\nolimits}
\def\res{\mathop{\mathrm{res}}\limits}
\def\sign{\mathop{\mathrm{sign}}\nolimits}
\def\span{\mathop{\mathrm{span}}\nolimits}
\def\supp{\mathop{\mathrm{supp}}\nolimits}
\def\Tr{\mathop{\mathrm{Tr}}\nolimits}
\def\BBox{\hspace{1mm}\vrule height6pt width5.5pt depth0pt \hspace{6pt}}


\newcommand\nh[2]{\widehat{#1}\vphantom{#1}^{(#2)}}
\def\dia{\diamond}

\def\Oplus{\bigoplus\nolimits}

%





\def\qqq{\qquad}
\def\qq{\quad}
\let\ge\geqslant
\let\le\leqslant
\let\geq\geqslant
\let\leq\leqslant
\newcommand{\ca}{\begin{cases}}
\newcommand{\ac}{\end{cases}}
\newcommand{\ma}{\begin{pmatrix}}
\newcommand{\am}{\end{pmatrix}}
\renewcommand{\[}{\begin{equation}}
\renewcommand{\]}{\end{equation}}
\def\eq{\begin{equation}}
\def\qe{\end{equation}}
\def\[{\begin{equation}}
\def\bu{\bullet}

\title[{Laplacians on periodic discrete graphs}]
{Laplacians on periodic discrete graphs}

\date{\today}

\author[Andrey Badanin]{Andrey Badanin}
\address{Department of Mathematical Analysis, Institute of Mathematics, Information and Space Technologies, Uritskogo St. 68, Northern (Arctic) Federal University,
Arkhangelsk, 163002,
 \ an.badanin@gmail.com}
\author[Evgeny Korotyaev]{Evgeny Korotyaev}
\address{Mathematical Physics Department, Faculty of Physics, Ulianovskaya 2,
St. Petersburg State University, St. Petersburg, 198904,
 and Pushkin Leningrad State University, Russia,
 \ korotyaev@gmail.com,}
\author[Natalia Saburova]{Natalia Saburova}
\address{Department of Algebra and Geometry, Institute of Mathematics, Information and Space Technologies, Uritskogo St. 68, Northern (Arctic) Federal University,
Arkhangelsk, 163002,
 \ n.saburova@gmail.com}


\begin{abstract}
We consider Laplacians on $\Z^2$-periodic discrete  graphs. The
following results are obtained: 1) The Floquet-Bloch decomposition
is constructed and basic properties are derived. 2) The estimates of the
Lebesgue measure of the spectrum in terms of geometric parameters of the graph
are obtained. 3) The spectrum of the
Laplacian is described, when the so-called fundamental graph
consists of one or two vertices and any number of edges.
4) We consider the hexagonal lattice perturbed  by adding one edge to the fundamental graph.
There exist two cases: a) if the perturbed  hexagonal lattice is  bipartite,
then the spectrum of the perturbed  Laplacian coincides with the spectrum
$[-1,1]$ for the  unperturbed   case, b) if the perturbed  hexagonal lattice is not  bipartite, then there is a gap in the spectrum of the perturbed Laplacian.
Moreover, some deeper results are obtained for the perturbation
of the square lattice.

\end{abstract}

\maketitle


\vskip 0.25cm

\section {Introduction and main results}

\subsection{Introduction.}
\setcounter{equation}{0}  We
discuss the spectrum of Laplacians on periodic discrete
graphs.  Laplacians on periodic graphs are of interest due to their applications to problems of physics and chemistry. They are used to describe and to study properties of different periodic media, including nanomedia, see \cite{NG04}, \cite{Ha85},
 \cite{SDD98}.

There are a lot of papers, and even books,  on the spectrum of the discrete Laplacian on an infinite graph. One of the main problems is to describe the spectral properties of the Laplacian in terms of geometric parameters of the graph. Useful
sources concerning operators acting on infinite graphs are the books \cite{Ba98}, \cite{Bo08}, \cite{BH12}, \cite{CDS95}, \cite{CDGT88}, \cite{P12}
and the papers \cite{HN09}, \cite{Me94},
\cite{MW89}, \cite{MRA07}. See also the references therein.

There are papers about the spectrum of discrete Laplacians on periodic graphs. Higuchi and Shirai \cite{HS04} (see also \cite{RR07}) obtain
the decomposition of the Laplacian into a constant fiber direct
integral. Higuchi and Nomura \cite{HN09} prove that the spectrum of
the Laplacian consists  of an absolutely continuous part and a
finite number of flat bands (i.e., eigenvalues with infinite
multiplicity). The absolutely continuous spectrum consists of a
finite number of intervals (spectral bands) separated by gaps.
Moreover, they also show that for each flat band there exists a
finitely supported eigenfunction.

There are results about spectral properties of the
discrete Laplace and Schr\"odinger operators on specific
periodic graphs. The hexagonal lattice can be viewed as a discrete
model of graphene, which is two-dimensional single-layered carbon
sheet with honeycomb structure. It was recently discovered by Geim and
Novoselov \cite{NG04}.
Graphene is a very hot subject in physics, where tight binding  models for the Schr\"odinger operator are standard and give
interesting band structures.  Ando \cite{A12} considers the spectral theory for the discrete Schr\"odinger operators with finitely supported potentials on the hexagonal lattice and their inverse scattering problem. Korotyaev and Kutsenko
\cite{KK10} -- \cite{KK10b} study the spectra of the discrete
Schr\"odinger operators on graphene nano-tubes and nano-ribbons in
external fields. See more about graphene in Section 5.
Schr\"odinger operators with decreasing potentials on the lattice $\Z^d$ are considered by Boutet de Monvel-Sahbani \cite{BS99}, Isozaki-Korotyaev \cite{IK12},
Rosenblum-Solomjak \cite{RoS09}.
Gieseker-Kn\"orrer-Trubowitz \cite{GKT93} consider
 Schr\"odinger operators with periodic potentials on the lattice $\Z^2$, the
simplest example of $\Z^2$-periodic graphs.  They study its
Bloch variety and its integrated density of states.

We describe our main goals of the paper:

1) to construct the Floquet theory  for Laplacians on periodic graphs.

2) to estimate the Lebesgue measure of the spectrum of the Laplacian
 in terms of geometric parameters of the graph.

3) to describe the spectrum of the Laplacian under the perturbation
of the square lattice and the hexagonal lattice.

4) It is known \cite{C97} that the investigation of the spectrum of the Laplacian on an equilateral quantum graph (i.e., a graph consisting of identical segments) can be reduced to the study of the spectrum of the discrete Laplacian. Thus, we have to do the needed spectral analysis of the discrete Laplacian given in this paper,
in order to describe spectral properties (including the Bethe-Sommerfeld conjecture) of the Laplacians on quantum graphs  \cite{KS}.

It should be noted that the results obtained in this work can be generalized to the case of $\Z^d$, $d\geq3$, periodic graphs \cite{KS1}.

\subsection{The definitions of periodic graphs and fundamental
graphs.} Let $\G=(V,\cE)$ be a connected graph, possibly  having loops and multiple edges, where $V$ is the set of its vertices and  $\cE$ is the
set of its unoriented edges. The graphs under consideration are embedded into $\R^2$. Considering each edge in $\cE$ to have two orientations, we can introduce the set $\cA$ of all oriented edges. The inverse edge of $\be\in\cA$ is denoted by $\bar\be$. The oriented edge starting at $u\in V$ and ending at $v\in V$ will be denoted as the ordered pair $(u,v)$. Vertices $u,v\in V$ will be called \emph{adjacent}
and denoted by $u\sim v$, if $(u,v)\in \cA$. We define the degree
${\vk}_v$ of the vertex $v\in V$ as the number of all oriented edges from $\cA$ starting at $v$. Below we consider $\Z^2$-periodic graphs $\G$,
satisfying the following conditions:

1) {\it the number of vertices from $V$ in any bounded domain $\ss\R^2$ is
finite;

2) the degree of each vertex is finite;

3) there exists a basis $a_1,a_2$ in $\R^2$ such that $\G$ is invariant under translations through the vectors $a_1$ and $a_2$:
$$
\G+a_1=\G, \qqq  \G+a_2=\G.
$$
The vectors $a_1,a_2$ are called the periods of $\G$.}

In the plane $\R^2$ we consider a coordinate system with the origin at some point $O$. The coordinate axes of this system are directed along the vectors $a_1$ and $a_2$. Below the coordinates of all vertices of $\G$ will be expressed  in this coordinate system. Then it follows from the definition of $\Z^2$-periodic graph  that $\G$ is invariant under translations through any integer vector:
$$
\G+p=\G,\qqq \forall\, p\in\Z^2.
$$
Examples of periodic graphs are shown in Figures
\ref{ff.0.3}\emph{a},  \ref{ff.0.1}\emph{a}.

We define \emph{the
fundamental graph} $\Gamma_0$ of the periodic graph $\Gamma$ as a graph
on the surface $\R^2/\Z^2$ by
\[
\lb{G0} \G_0=\G/{\Z}^2\ss T^2, \quad \textrm{where}\quad
T^2=\R^2/\Z^2.
\]
The vertex set $V_0$, the set $\cE_0$ of unoriented edges and the set $\cA_0$ of oriented edges of $\G_0$ are finite (see Proposition \ref{pro0}).

We introduce {\it an edge "index"}, which is important in the study of the Laplace operator. We identify the vertices of the fundamental graph $\G_0=(V_0,\cE_0)$ with the vertices of the periodic graph $\G=(V,\cE)$ in the set $[0,1)^2$. Then for any $v\in V$ the following unique representation holds true:
\[
\lb{Dv}
v=[v]+\tilde v, \qquad [v]\in\Z^2,\qquad \tilde v\in V_0\subset[0,1)^2.
\]
In other words, each vertex $v$ can be represented uniquely as the sum of an integer part $[v]\in \Z^2$ and a fractional part $\tilde v$ that is a vertex of the fundamental graph $\G_0$. For any oriented edge $\be=(u,v)\in\cA$
we define {\bf the edge "index"}  $\t({\bf e})$ as the integer vector
\[
\lb{in}
\t({\bf e})=[v]-[u]\in\Z^2,
\]
where
$$
u=[u]+\tilde{u},\qquad v=[v]+\tilde{v}, \qquad [u], [v]\in\Z^2,\qquad \tilde{u},\tilde{v}\in V_0.
$$
If $\be=(u,v)$ is an oriented edge of the graph $\G$, then
by the definition of the fundamental graph there is an oriented edge $\tilde\be=(\tilde u,\tilde v\,)$ on $\G_0$. For an edge $\tilde\be\in\cA_0$ we define the edge index $\t(\tilde{\bf e})$ by
\[
\lb{inf}
\t(\tilde{\bf e})=\t(\be).
\]
In other words, indices of periodic graph edges are inherited by edges of the fundamental graph. The edge indices, generally speaking, depend on the choice of the coordinate origin $O$. But in a fixed coordinate system
the index of the fundamental graph edge is uniquely determined by \er{inf}, since
due to Proposition \ref{pro1}.iii we have
$$
\t(\be+p)=\t(\be),\qqq \forall\, (\be, p)\in\cA \ts \Z^2.
$$
Edges with nonzero indices will be called {\bf bridges}. They  are
important, when we describe the spectrum of the Laplacian.
 The bridges provide the connectivity of the periodic graph and the removal of all bridges disconnects the graph into infinitely many connected components.

\subsection{Laplace operators on graphs.} Let $\ell^2(V)$ be the
Hilbert space of all square summable functions $f:V\to \C$, equipped
with the norm
$$
\|f\|^2_{\ell^2(V)}=\sum_{v\in V}|f(v)|^2<\infty.
$$
Recall that the degree ${\vk}_v$ of the vertex $v\in V$ is the number of oriented edges starting at the
vertex $v$. We define the Laplacian (or the Laplace operator) $\D$ acting on the Hilbert space
$\ell^2(V)$ by
\[
\lb{DOL}
 \big(\D f\big)(v)=\frac1{\sqrt{\vk_v}} \sum\limits_{(v,\,u)\in\cA}\frac1{\sqrt{\vk_u}}\,f(u), \qquad \forall f\,\in \ell^2(V).
\]
Note that sometimes $1-\D\ge 0$ is also called
the Laplace operator.

A graph is called \emph{bipartite} if its vertex set is divided into two disjoint sets (called \emph{parts} of the graph) such that each edge connects vertices from distinct sets (see p.105 in \cite{Or62}). Examples of bipartite graphs are the square lattice (Figure \ref{ff.0.1}\emph{a}) and the hexagonal lattice (Figure \ref{ff.0.3}\emph{a}). The  triangular lattice (Figure \ref{f.0.2}\emph{b}) is non-bipartite. Note that for a bipartite periodic graph there exists a bipartite
fundamental graph (see Lemma \ref{l1}.ii), but not every fundamental graph is bipartite. Indeed, the square lattice (Figure \ref{ff.0.1}\emph{a}) is bipartite, but its fundamental graph shown in Figure \ref{ff.0.1}\emph{b} is non-bipartite.

We recall well-known properties of the Laplacian $\Delta$, which
hold true  for  finite and infinite graphs (in particular, for
periodic graphs) (see \cite{Ch97}, \cite{Me94}, \cite{M92}, \cite{MW89}):

\textbf{Main properties of the Laplacian}:
\begin{itemize}
  \item[1)] \emph{The operator $\Delta$ is self-adjoint and bounded.}
  \item[2)] \emph{The spectrum $\s(\Delta)$ is contained in  $[-1,1]$.}
  \item[3)] \emph{The point 1 belongs to the spectrum $\sigma(\Delta)$.}
  \item[4)] \emph{The graph is bipartite $\Leftrightarrow $ the point $-1\in\sigma(\Delta)$ $\Leftrightarrow$ the  spectrum $\sigma(\Delta)$ is symmetric with
  respect to the point zero.}
\end{itemize}

Denote by $v_1,\ldots,v_\nu$ the vertices of $V_0$, where
$\nu<\iy$ is a number of vertices of $V_0$. Thus we have
$$
V_0=\{v_j: j\in \N_\n\},\qqq \N_\n=\{1,\ldots,\nu\}.
$$

We present our first results about the decomposition of the Laplacian $\D$ into a direct integral.

\medskip

\begin{theorem}\label{pro2}
i) The operator $\D$ acting on $\ell^2(V)$ has the decomposition
into a constant fiber direct integral
\[
\lb{raz}
\begin{aligned}
& \ell^2(V)={1\/(2\pi)^2}\int^\oplus_{\T^2}\ell^2(V_0)\,d\vt ,\qqq
\T^2=\R^2/(2\pi\Z)^2=[-\pi;\pi]^2,\\
& U\D U^{-1}={1\/(2\pi)^2}\int^\oplus_{\T^2}\D(\vt)d\vt,
\end{aligned}
\]
for some unitary operator $U$. Here the fiber space
$\ell^2(V_0)=\C^\nu$ and $\D(\vt )=\{\D_{jk}(\vt )\}$, $\vt\in \T^2$, is a
Floquet  $\nu\ts\nu$ matrix given by
\medskip
\[
\label{l2.15} \D_{jk}(\vt )=\left\{\begin{array}{cl}
{1\/\sqrt{\vk_j\vk_k}}\sum\limits_{{\bf e}=(v_j,\,v_k)\in{\cA}_0}e^{\,i\lan\t
({\bf e}),\,\vt\ran }, \qqq &  {\rm if}\  \ (v_j,v_k)\in \cA_0 \\
[18pt]
 0, &  {\rm if}\  \ (v_j,v_k)\notin \cA_0  \\
\end{array}\right.,
\]
where $\vk_j$ is the degree of $v_j$ and
$\lan\cdot\,,\cdot\ran$ denotes the standard inner product in
$\R^2$. In particular,
\[
\label{l2.15'} \D_{jj}(\vt )=\left\{\begin{array}{cl}
\frac1{\vk_j}\sum\limits_{\mathbf{e}=
(v_j,\,v_j)\in\cA_0}\cos \lan\t (\mathbf{e}),\vt\ran , \quad &
\ {\rm if \ a\  loop}\  \ (v_j,v_j)\in \cA_0\\[18pt]
 0, & \ {\rm if}\  \ (v_j,v_j)\notin \cA_0  \\
\end{array}\right..
\]

ii) Let $\D^{(1)}(\vt )$ be the Floquet matrix for $\D$ defined by (\ref{l2.15}) in another coordinate system with the origin $O_1$. Then the matrices $\D^{(1)}(\vt )$ and $\D(\vt)$ are unitarily equivalent for each $\vt \in\T^2$.

iii) The spectrum $\s\big(\D(\vt )\big)$ is contained in  $[-1,1]$ for
each $\vt \in\T^2$.

iv) The entry $\D_{jk}(\cdot)$ is constant iff there is no bridge $(v_j,v_k)\in\cA_0$.

v) The point 1 is never a flat band of $\D$.

vi) The Floquet matrix  $\D(\cdot)$ has at least one non-constant entry $\D_{jk}(\cdot)$ for some $j\le k$.

vii) A fundamental graph $\G_0$ is bipartite iff  the spectrum $\s(\D(\vt))$ is symmetric with respect to the point zero for each $\vt \in\T^2$.
\end{theorem}

\no \textbf{Remark.} 1) The representations \er{l2.15}, \er{l2.15'} are
new. Recall that the existence of $\D(\vt)$ was proved in \cite{HS04}.

\setlength{\unitlength}{1.0mm}
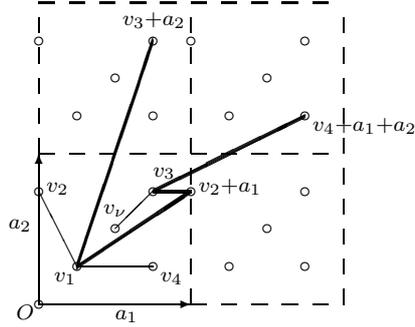
\begin{figure}[h]
\centering
\unitlength 1mm 
\linethickness{0.4pt}
\ifx\plotpoint\undefined\newsavebox{\plotpoint}\fi 

\begin{picture}(60,50)(0,0)

\multiput(10,10)(4,0){10}{\line(1,0){2}}
\multiput(10,30)(4,0){10}{\line(1,0){2}}
\multiput(10,50)(4,0){10}{\line(1,0){2}}

\multiput(10,10)(0,4){10}{\line(0,1){2}}
\multiput(30,10)(0,4){10}{\line(0,1){2}}
\multiput(50,10)(0,4){10}{\line(0,1){2}}
\put(10,10){\circle{1}}

\put(10,10){\vector(1,0){20.00}}
\put(10,10){\vector(0,1){20.00}}

\put(7.0,8.0){$\scriptstyle O$}

\put(20,8){$\scriptstyle a_1$}
\put(6,20){$\scriptstyle a_2$}

\put(15,15){\line(1,0){10.00}}
\put(25,25){\line(1,0){5.00}}
\put(25,25.1){\line(1,0){5.00}}
\put(25,24.9){\line(1,0){5.00}}
\put(25,24.8){\line(1,0){5.00}}

\put(25,25.1){\line(2,1){20.00}}
\put(25,25.2){\line(2,1){20.00}}
\put(25,24.9){\line(2,1){20.00}}
\put(25,24.8){\line(2,1){20.00}}
\put(25,25){\line(2,1){20.00}}

\put(15,15){\line(-1,2){5.00}}

\put(15,15.1){\line(3,2){15.00}}
\put(15,15.2){\line(3,2){15.00}}
\put(15,14.9){\line(3,2){15.00}}
\put(15,14.8){\line(3,2){15.00}}
\put(15,15){\line(3,2){15.00}}

\put(15,15.4){\line(1,3){10.00}}
\put(15,15.2){\line(1,3){10.00}}
\put(15,15){\line(1,3){10.00}}
\put(15,14.8){\line(1,3){10.00}}
\put(15,14.6){\line(1,3){10.00}}
\put(15,15){\line(1,3){10.00}}

\put(20,20){\line(1,1){5.00}}

\put(15,15){\circle{1}}
\put(10,25){\circle{1}}
\put(20,20){\circle{1}}
\put(25,25){\circle{1}}
\put(25,15){\circle{1}}
\put(12,13){$\scriptstyle v_1$}
\put(11,25){$\scriptstyle v_2$}\put(31,25){$\scriptstyle v_2+a_1$}
\put(26,13){$\scriptstyle v_4$}\put(46,33){$\scriptstyle v_4+a_1+a_2$}
\put(25,27){$\scriptstyle v_3$}\put(21,47){$\scriptstyle v_3+a_2$}
\put(18.5,22){$\scriptstyle v_\nu$}


\put(35,15){\circle{1}}
\put(30,25){\circle{1}}
\put(40,20){\circle{1}}
\put(45,25){\circle{1}}
\put(45,15){\circle{1}}

\put(15,35){\circle{1}}
\put(10,45){\circle{1}}
\put(20,40){\circle{1}}
\put(25,45){\circle{1}}
\put(25,35){\circle{1}}


\put(35,35){\circle{1}}
\put(30,45){\circle{1}}
\put(40,40){\circle{1}}
\put(45,45){\circle{1}}
\put(45,35){\circle{1}}
\end{picture}

\vspace{-0.5cm}
\caption{A periodic graph $\G$; only edges of the fundamental graph $\G_0$ are shown; the bridges of $\G$
are marked by bold lines.} \label{ff.0.11}
\end{figure}

 2) The graph $\G$ with $\n=5$ on Figure \ref{ff.0.11} has the following bridges
$$
(v_1,v_2+a_1),\quad (v_1,v_3+a_2), \quad (v_3,v_2+a_1), \quad (v_3,v_4+a_1+a_2)
$$
because they connect vertices with different integer parts in the sense of the
identity (\ref{Dv}).

3) There are always some bridges on a graph $\G$,
if $\G$ is connected.
Otherwise, if there are no bridges on some graph $\G$, then this
graph is not connected and is a union of infinitely many connected components.
In this case the spectrum of the Laplacian on $\G$ consists of only
a finite number of eigenvalues with infinite multiplicity.
Due to Theorem \ref{Tfb} i, it is impossible. The presence of bridges on the graph "enlarges" the spectrum of the Laplacian.

\

\subsection{Spectral properties of Laplacians.}
Theorem \ref{pro2} and standard arguments (see Theorem XIII.85 in
\cite{RS78})  describe the spectrum of the Laplacian. Each Floquet
${\nu\ts\nu}$ matrix  $\D(\vt)$, $\vt\in\T^2$, has $\n$ eigenvalues
$\l_n(\vt)$, $n=1,\ldots,\n$. They are real and locally analytic functions on the torus $\T^2$, since $\D(\vt)$ is  self-adjoint and analytic in $\vt$ on the torus $\T^2$. The parameter $\vt$ is called {\it quasimomentum}.
If some $\l_n(\cdot)=C_n=\const$ on some set $\cB\ss\T^2$ of positive Lebesgue measure, then the operator $\Delta$ on $\G$ has  the eigenvalue $C_n$
with infinite multiplicity. We will call $C_n$ a \emph{flat band}.
Due to \cite{HN09} we have that  $\l_*$ is an eigenvalue of $\D$ iff $\l_*$ is an eigenvalue of $\D(\vt)$ for any $\vt\in\T^2$.
 Thus, we can define the multiplicity of a flat band by:
an eigenvalue $\l_*$  of $\D$ has the multiplicity $m$ iff
$\l_*=\const$ is an eigenvalue of $\D(\vt)$ for  each $\vt\in\T^2$ with the multiplicity $m$ (except maybe for a finite number of $\vt\in\T^2$).
Thus if the operator $\D$ has $r\ge 0$ flat bands, then we denote them by
$$
\m_j=\l_{\n-j+1}(\vt), \qqq j\in\N_r,
$$
and they are labeled by
\[
\m_1\le \m_2\le\ldots\le  \m_r,
\]
counting multiplicities. Thus, all other eigenvalues
$\l_n(\vt)$, $n=1,\ldots,\n-r$ are not constant. They can be enumerated
in decreasing order (counting multiplicities)
\[
\label{eq.3}
\l_{\nu-r}(\vt )\leq\l_{\nu-r-1}(\vt )\leq\ldots\leq\l_1(\vt ), \qqq \forall\ \vt\in\T^2.
\]
Each $\l_n(\vt )$ is a piecewise analytic function on the torus $\T^2$
and defines \emph{a dispersion relation}. Its plot is a
\emph{dispersion curve}. Define the \emph{spectral bands} $\s_n$  by
\[
\lb{ban}
\s_n=[\l_n^-,\l_n^+]=\l_n(\T^2), \quad \textrm{where}\quad\l_n^-=\min_{\vt
\in\T^2}\l_n(\vt ),\quad \l_n^+=\max_{\vt \in\T^2}\l_n(\vt ),\qq \forall\,n\in\N_{\n}.
\]
For each $n=1,\ldots, \n-r$ we have $\l_n^-<\l_n^+$ and the spectral band $\s_n$ is open (non-degenerate). For each $n=\n-r+1,\ldots,\n$ we have $\l_n^\pm=\m_r$ and the spectral band $\s_n$ is close (degenerate).
If $\l_{n+1}^+<\l_{n}^-$ for some $n\in\N_{\n-r-1}$, then the interval $(\l_{n+1}^+,\l_{n}^-)$ is called a \emph{gap}.

Thus the spectrum of the operator $\D$ on the periodic graph $\G$ has the form
\[
\lb{r0}
\begin{aligned}
&\s(\D)=\s_{ac}(\D)\cup \s_{fb}(\D),\\
&\s_{ac}(\D)=\bigcup_{n=1}^{\nu-r}\s_n,\qqq \s_{fb}(\D)=\{\m_1,\ldots,\m_r\},
\end{aligned}
\]
where $\s_{ac}(\D)$ is the absolutely continuous spectrum and $\s_{fb}(\D)$ is a set of all flat bands (eigenvalues with infinite multiplicity).
We now describe precisely all bands for specific graphs.

\begin{theorem}\lb{T100}
i) Let all bridges of the fundamental graph $\G_0$ be loops, i.e., the
indices of all edges $(v_j,v_k)$, $1\leq j<k\leq \nu$, of $\G_0$ be zero. Then
\[
\lb{eq.5} \l_n^+=\l_n(0),\qqq \forall \;n\in \N_\n.
\]
Moreover, the eigenvalue $\l_n(0)$ is a flat band of $\D$ iff $\l_n^-=\l_n^+$. The multiplicity of the flat band $\m=\l_n^+$ of the Laplacian $\D$ is the multiplicity of  $\l_n^+$ as the eigenvalue of the operator $\D(0)$.

ii) Let, in addition, $\cos\lan\t
({\bf e}),\,\vt_0\ran=-1$ for all bridges $\be\in\cE_0$ and some $\vt_0\in\T^2$.
Then
\[
\lb{es}
\s_n=[\l_n^-,\l_n^+],\qqq \l_n^-=\l_n(\vt_0), \qqq \forall \;n\in \N_\n.
\]

iii) Let $\G$ be a bipartite periodic graph, satisfying  the condition i)     ($\G_0$ is not bipartite, since there is a loop on $\G_0$).
Then $\l_n^\pm$ are the eigenvalues of the matrix $\pm\D(0)$.
\end{theorem}

\no\textbf{Remark.} 1) $\l_n^+$, $n\in\N_\n$, are the eigenvalues  of the
operator  $\D(0)$ defined by (\ref{l2.15}) that is the Laplacian on the fundamental graph $\G_0$.
The case of item ii) is similar to the case of $N$-periodic Jacobi
 matrices on the lattice $\Z$ (and for Hill operators).
The spectrum of these operators is absolutely continuous and is a union
of spectral bands, separated by gaps. The endpoints of the bands are so-called $2N$-periodic eigenvalues.

2) For some classes of graphs one can easily determine $\vt_0$ such that $\cos\lan\t({\bf e}),\,\vt_0\ran=-1$ for all bridges $\be\in\cA_0$. Let $\t(\be)=(\t_1(\be),\t_2(\be))\in\Z^2$ be the index of a bridge $\be$ of $\G_0$. Then
\[
\lb{l-}
\vt_0=\left\{\begin{array}{ll}
  (\pi,0), & \textrm{if $\t_1(\be)$ is odd for all bridges $\be\in\cA_0$}\\
  (0,\pi), & \textrm{if $\t_2(\be)$ is odd for all bridges $\be\in\cA_0$}\\
  (\pi,\pi),\qqq & \textrm{if $\t_1(\be)+\t_2(\be)$ is odd for all bridges
$\be\in\cA_0$}\\
\end{array}\right..
\]
\setlength{\unitlength}{1.0mm}
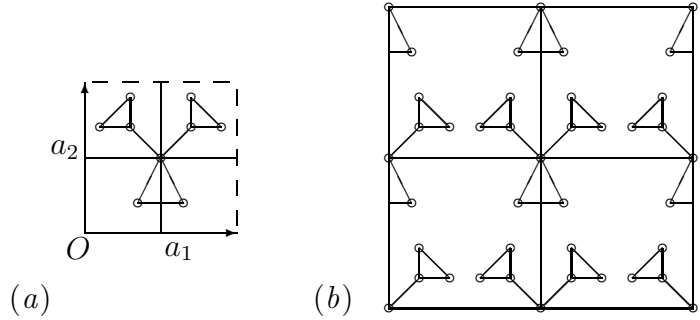
\begin{figure}[h]
\centering
\unitlength 1mm 
\linethickness{0.4pt}
\ifx\plotpoint\undefined\newsavebox{\plotpoint}\fi 
\begin{picture}(120,50)(0,0)

\put(60,10){\line(1,0){40.00}}
\put(60,30){\line(1,0){40.00}}
\put(60,50){\line(1,0){40.00}}
\put(60,10){\line(0,1){40.00}}
\put(80,10){\line(0,1){40.00}}
\put(100,10){\line(0,1){40.00}}

\put(60,10){\circle{1}}
\put(80,10){\circle{1}}
\put(100,10){\circle{1}}

\put(60,30){\circle{1}}
\put(80,30){\circle{1}}
\put(100,30){\circle{1}}

\put(60,50){\circle{1}}
\put(80,50){\circle{1}}
\put(100,50){\circle{1}}

\put(30,30){\circle{1}}
\put(20,20){\vector(1,0){20.00}}
\put(20,20){\vector(0,1){20.00}}
\multiput(21,40)(4,0){5}{\line(1,0){2}}
\multiput(40,21)(0,4){5}{\line(0,1){2}}
\put(17.5,16.5){$O$}
\put(30.5,17.0){$a_1$}
\put(15.5,30.5){$a_2$}
\put(20,30){\line(1,0){20.00}}
\put(30,20){\line(0,1){20.00}}

\put(30,30){\line(1,1){4.00}}
\put(34,34){\line(1,0){4.00}}
\put(34,34){\line(0,1){4.00}}
\put(34,38){\line(1,-1){4.00}}
\put(34,34){\circle{1}}
\put(34,38){\circle{1}}
\put(38,34){\circle{1}}

\put(30,30){\line(-1,1){4.00}}
\put(26,34){\line(-1,0){4.00}}
\put(26,34){\line(0,1){4.00}}
\put(26,38){\line(-1,-1){4.00}}
\put(26,34){\circle{1}}
\put(26,38){\circle{1}}
\put(22,34){\circle{1}}

\put(30,30){\line(-1,-2){3.00}}
\put(30,30){\line(1,-2){3.00}}
\put(27,24){\line(1,0){6.00}}
\put(27,24){\circle{1}}
\put(33,24){\circle{1}}

\put(10,10){(\emph{a})}
\put(50,10){(\emph{b})}

\put(80,30){\line(1,1){4.00}}
\put(84,34){\line(1,0){4.00}}
\put(84,34){\line(0,1){4.00}}
\put(84,38){\line(1,-1){4.00}}
\put(84,34){\circle{1}}
\put(84,38){\circle{1}}
\put(88,34){\circle{1}}

\put(80,30){\line(-1,1){4.00}}
\put(76,34){\line(-1,0){4.00}}
\put(76,34){\line(0,1){4.00}}
\put(76,38){\line(-1,-1){4.00}}
\put(76,34){\circle{1}}
\put(76,38){\circle{1}}
\put(72,34){\circle{1}}

\put(80,30){\line(-1,-2){3.00}}
\put(80,30){\line(1,-2){3.00}}
\put(77,24){\line(1,0){6.00}}
\put(77,24){\circle{1}}
\put(83,24){\circle{1}}

\put(100,30){\line(-1,1){4.00}}
\put(96,34){\line(-1,0){4.00}}
\put(96,34){\line(0,1){4.00}}
\put(96,38){\line(-1,-1){4.00}}
\put(96,34){\circle{1}}
\put(96,38){\circle{1}}
\put(92,34){\circle{1}}

\put(100,30){\line(-1,-2){3.00}}
\put(97,24){\line(1,0){3.00}}
\put(97,24){\circle{1}}
\put(60,30){\line(1,1){4.00}}
\put(64,34){\line(1,0){4.00}}
\put(64,34){\line(0,1){4.00}}
\put(64,38){\line(1,-1){4.00}}
\put(64,34){\circle{1}}
\put(64,38){\circle{1}}
\put(68,34){\circle{1}}

\put(60,30){\line(1,-2){3.00}}
\put(60,24){\line(1,0){3.00}}
\put(63,24){\circle{1}}
\put(80,10){\line(1,1){4.00}}
\put(84,14){\line(1,0){4.00}}
\put(84,14){\line(0,1){4.00}}
\put(84,18){\line(1,-1){4.00}}
\put(84,14){\circle{1}}
\put(84,18){\circle{1}}
\put(88,14){\circle{1}}

\put(80,10){\line(-1,1){4.00}}
\put(76,14){\line(-1,0){4.00}}
\put(76,14){\line(0,1){4.00}}
\put(76,18){\line(-1,-1){4.00}}
\put(76,14){\circle{1}}
\put(76,18){\circle{1}}
\put(72,14){\circle{1}}

\put(100,10){\line(-1,1){4.00}}
\put(96,14){\line(-1,0){4.00}}
\put(96,14){\line(0,1){4.00}}
\put(96,18){\line(-1,-1){4.00}}
\put(96,14){\circle{1}}
\put(96,18){\circle{1}}
\put(92,14){\circle{1}}

\put(60,10){\line(1,1){4.00}}
\put(64,14){\line(1,0){4.00}}
\put(64,14){\line(0,1){4.00}}
\put(64,18){\line(1,-1){4.00}}
\put(64,14){\circle{1}}
\put(64,18){\circle{1}}
\put(68,14){\circle{1}}
\put(80,50){\line(-1,-2){3.00}}
\put(80,50){\line(1,-2){3.00}}
\put(77,44){\line(1,0){6.00}}
\put(77,44){\circle{1}}
\put(83,44){\circle{1}}
\put(100,50){\line(-1,-2){3.00}}
\put(97,44){\line(1,0){3.00}}
\put(97,44){\circle{1}}
\put(60,50){\line(1,-2){3.00}}
\put(60,44){\line(1,0){3.00}}
\put(63,44){\circle{1}}
\end{picture}
\vspace{-0.5cm}
\caption{\emph{a}) The fundamental graph $\G_0$;\quad \emph{b}) the periodic graph $\G$.} \label{ffS'}
\end{figure}

3) The fundamental graph $\G_0$ on Figure \ref{ffS'}\emph{a} has only 4 oriented bridges, which are loops with indices $(0,\pm1)$ and $(\pm1,0)$. All other edges of $\G_0$ have zero indices.
Thus, $\t_1(\be)+\t_2(\be)$ is odd for all bridges $\be\in\cA_0$ and according to Theorem \ref{T100}.i--ii and the identity \er{l-} the spectrum of the Laplacian on $\G$, see Figure \ref{ffS'}\emph{b}, has the form
$
\s(\D)=\bigcup_{n=1}^{\nu}\s_n$, \ $\s_n=[\l_n(\pi,\pi),\l_n(0)].
$

\

We estimate the spectrum $\s(\D)$ of the Laplacian.

\begin{theorem}
\lb{T1} i) The Lebesgue measure
$|\s(\D)|$ of the spectrum of the Laplacian $\D$ satisfies
\[
\lb{eq.7}
|\s(\D)|\le \sum_{n=1}^{\n-r}|\s_n|\leq2\sum_{j,\,k=1}^\n\frac{b_{jk}}{\sqrt{\vk_j\vk_k}}\,,
\]
where $b_{jk}$ is the number of bridges $(v_j,v_k)$ on the fundamental graph $\G_0$, $\vk_j$ is the degree of $v_j$. Moreover,
this estimate is sharp, i.e.,  estimates \er{eq.7} become identities for some graphs.

ii)  For any $\n\ge 2$ there exists a periodic graph such that
$|\s(\D)|={8\/\n+3}$. In particular, $|\s(\D)|\ra0$ as $\n\to\iy$.

iii) Let  a fundamental graph $\G_0$ be bipartite. If $\nu>1$ is odd, then $\m=0$ is a flat band of $\D$.

\end{theorem}

\no \textbf{Remark.}
1) The estimate \er{eq.7} is effective for the case
\[
\lb{es'}
\sum\limits_{j,\,k=1}^\n\frac{b_{jk}}{\sqrt{\vk_j\vk_k}}<1.
 \]
If the inequality \er{es'} does not hold true, then \er{eq.7} yields only a simple estimate $|\s(\D)|\le 2$, which follows from the basic property $\s(\Delta)\ss[-1,1]$.
The inequality \er{es'} holds true when for each vertex $v\in V_0$  the ratio of the number of bridges starting at $v$ to
the degree of $v$ is rather small. Increasing the degree of each
vertex and fixing the bridges we can
construct a graph such that the Lebesgue measure of the spectrum of the Laplacian is
arbitrarily small.

2) The total length of spectral  bands depends essentially on the bridges on $\G_0$. If we remove the coordinate system, then the number of bridges on the corresponding fundamental graph $\G_0$ is changed in general.

\setlength{\unitlength}{1.0mm}
\begin{figure}[h]
\centering
\unitlength 1mm 
\linethickness{0.4pt}
\ifx\plotpoint\undefined\newsavebox{\plotpoint}\fi 

\begin{picture}(60,40)(0,0)

\put(30,10){\line(3,1){10.00}}
\put(40,13.3){\line(1,3){3.30}}
\put(30,10){\line(-3,1){10.00}}
\put(20,13.3){\line(-1,3){3.30}}

\put(30,36.3){\line(3,-1){10.00}}
\put(40,33.0){\line(1,-3){3.30}}
\put(30,36.3){\line(-3,-1){10.00}}
\put(20,33.0){\line(-1,-3){3.30}}

\put(40,13.3){\line(1,-1){5.00}}
\put(20,13.3){\line(-1,-1){5.00}}
\put(40,33.0){\line(1,1){5.00}}
\put(20,33.0){\line(-1,1){5.00}}

\put(30,10){\circle{1}}
\put(40,13.3){\circle{1}}
\put(16.7,23.15){\circle{1}}
\put(20,13.3){\circle{1}}
\put(30,36.3){\circle{1}}
\put(40,33.0){\circle{1}}
\put(43.3,23.15){\circle{1}}
\put(20,33.0){\circle{1}}

\put(17,13){$\scriptstyle 0$}
\put(41.3,13){$\scriptstyle 0$}
\put(17,32){$\scriptstyle 0$}
\put(41.3,32){$\scriptstyle 0$}
\put(29.3,6){$\scriptstyle 1$}
\put(29.3,38){$\scriptstyle 1$}
\put(11.5,23){$\scriptstyle -1$}
\put(44.5,23){$\scriptstyle -1$}
\end{picture}

\vspace{-0.5cm}
\caption{The cycle that is the support of the eigenfunction with the eigenvalue 0. The values of the eigenfunction in the vertices of the support are $0,\pm1$.} \label{ffc}
\end{figure}
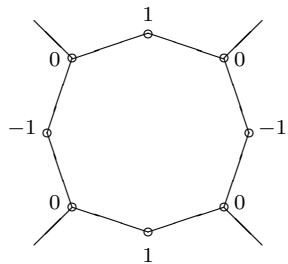

3) Recall that for each flat band there exists a
finitely supported eigenfunction.
It is well-known that for some specific case these eigenfunctions can easily be
determined using some analysis of graph cycles. For example, if the graph has a cycle of length $4n$, $n\in\N$, and the degree of each second vertices of the cycle is 2, then one can easily construct an eigenfunction supported on this cycle (see Figure \ref{ffc}) with the eigenvalue 0. Theorems \ref{T100} and \ref{T1}.iii give the simple  sufficient conditions for the existence of flat bands, not based on an analysis of cycles. This is the only test for the existence
of flat bands, that we know of.

Now we describe possible  positions of flat bands.

\begin{theorem}\lb{Tfb}
i) The first spectral band $\s_1=[\l_1^-,1]$ is open, i.e., $\l_1^-<1$.
Moreover, the  number of flat bands $r<\n$ and
the spectrum of the operator $\Delta$ has at most
$\nu-r-1$ gaps.

ii) Let $\n\geq2$. Then


1) There is a graph $\G$, such that the spectrum of the Laplacian on $\G$ has exactly 2 open spectral bands and between them, in the gap,  $\n-2$ flat bands, counting multiplicity.

2) There is a graph $\G$, such that  the spectrum of the Laplacian on $\G$ has  $[\frac\n2]$ different simple flat bands embedded in the absolutely continuous spectrum $[-1,1]$.

3) There is a graph $\G$, such that  a point $-{1\/2}$ is a flat band of the Laplacian on $\G$ (the Kagome lattice, see subsection \ref{KL}) and lies at the endpoint of the absolutely continuous spectrum $\sigma_{ac}(\Delta)=[-{1\/2},1]$.
\end{theorem}

There is an open problem: does there exist a graph with any $\n\geq2$ vertices in the fundamental graph such that the spectrum of the Laplacian on $\G$ has only 1 spectral band and $\n-1$ flat bands, counting multiplicity?

\subsection{Perturbations of the hexagonal lattice.}
Let $\bG=(V,\cE)$ be the hexagonal
lattice (Figure \ref{ff.0.3}\emph{a}).
The periods of $\bG$ are the
vectors $a_1=(3/2,\sqrt{3}/2)$, $a_2=(0,\sqrt{3}\,)$ (the
coordinates of $a_1$, $a_2$ are taken in the orthonormal basis
$e_1$, $e_2$). The vertex set and the edge set are given by
$$
\textstyle V=\Z^2\cup\big(\Z^2+\big(\frac13\,,\frac13\big)\big),
$$
$$
\textstyle \cE=\big\{\big(p,p+\big(\frac13\,,\frac13\big)\big),
\big(p,p+\big(-\frac23\,,\frac13\big)\big),\\
\big(p,p+\big(\frac13\,,-\frac23\big)\big)\quad\forall\,p\in\Z^2\big\}.
$$
Recall that the coordinates of all vertices are taken in the basis $a_1$, $a_2$. The fundamental graph $\bG_0$ consists of two vertices $v_1$,
$v_2$, \emph{multiple} edges
$\be_1=\be_2=\be_3=(v_1,v_2)$ (Figure
\ref{ff.0.3}\emph{b}) and $\bar\be_1=\bar\be_2=\bar\be_3$ with the indices
$$
\t(\be_1)=\t(\bar\be_1)=(0,0),
\qq \t(\be_2)=-\t(\bar\be_2)=(1,0),\qq \t(\be_3)=-\t(\bar\be_3)=(0,1).
$$
It is known that the spectrum of the Laplacian $\D$ on $\bG$ has the form
$\s(\D)=\s_{ac}(\D)=[-1,1]$.

\setlength{\unitlength}{1.0mm}
\begin{figure}[h]
\centering

\unitlength 1mm 
\linethickness{0.4pt}
\ifx\plotpoint\undefined\newsavebox{\plotpoint}\fi 
\begin{picture}(100,45)(0,0)

\put(5,10){(\emph{a})}

\put(14,10){\circle{1}}
\put(28,10){\circle{1}}
\put(34,10){\circle{1}}
\put(48,10){\circle{1}}

\put(18,16){\circle{1}}
\put(24,16){\circle*{1}}
\put(38,16){\circle{1}}
\put(44,16){\circle{1}}

\put(14,22){\circle{1}}
\put(28,22){\circle*{1}}
\put(34,22){\circle{1}}
\put(48,22){\circle{1}}

\put(18,28){\circle{1}}
\put(24,28){\circle{1}}
\put(38,28){\circle{1}}
\put(44,28){\circle{1}}

\put(14,34){\circle{1}}
\put(28,34){\circle{1}}
\put(34,34){\circle{1}}
\put(48,34){\circle{1}}

\put(18,40){\circle{1}}
\put(24,40){\circle{1}}
\put(38,40){\circle{1}}
\put(44,40){\circle{1}}

\put(28,10){\line(1,0){6.00}}
\put(18,16){\line(1,0){6.00}}
\put(38,16){\line(1,0){6.00}}

\put(28,22){\line(1,0){6.00}}
\put(18,28){\line(1,0){6.00}}
\put(38,28){\line(1,0){6.00}}

\put(28,34){\line(1,0){6.00}}
\put(18,40){\line(1,0){6.00}}
\put(38,40){\line(1,0){6.00}}

\put(14,10){\line(2,3){4.00}}
\put(34,10){\line(2,3){4.00}}
\put(24,16){\line(2,3){4.00}}
\put(44,16){\line(2,3){4.00}}

\put(14,22){\line(2,3){4.00}}
\put(34,22){\line(2,3){4.00}}
\put(24,28){\line(2,3){4.00}}
\put(44,28){\line(2,3){4.00}}

\put(14,34){\line(2,3){4.00}}
\put(34,34){\line(2,3){4.00}}

\put(28,10){\line(-2,3){4.00}}
\put(48,10){\line(-2,3){4.00}}
\put(38,16){\line(-2,3){4.00}}
\put(18,16){\line(-2,3){4.00}}

\put(28,22){\line(-2,3){4.00}}
\put(48,22){\line(-2,3){4.00}}
\put(38,28){\line(-2,3){4.00}}
\put(18,28){\line(-2,3){4.00}}

\put(28,34){\line(-2,3){4.00}}
\put(48,34){\line(-2,3){4.00}}

\put(30,18){$\scriptstyle a_1$}
\put(20.5,22){$\scriptstyle a_2$}

\put(24,16){\vector(0,1){12.0}}
\put(33,21.3){\vector(3,2){0.5}}

\qbezier(24,16)(29,19)(34,22)

\put(24.8,21.5){$\scriptstyle v_1$}
\put(35,21.5){$\scriptstyle v_2+a_1$}
\put(25.0,28.0){$\scriptstyle v_2+a_2$}
\put(17.8,13.5){$\scriptstyle O=v_2$}

\put(75,10){\circle*{1}}
\put(83,21){\circle*{1}}

\put(75,30){\circle{1}}
\put(95,40){\circle{1}}
\put(95,20){\circle{1}}

\put(75,10){\vector(0,1){20.0}}
\put(75,10){\vector(2,1){20.0}}

\multiput(95,20)(0,7){3}{\line(0,1){4}}
\put(75,30){\line(2,1){4.0}}
\put(82,33.5){\line(2,1){4.0}}
\put(89,37){\line(2,1){4.0}}

\qbezier(83,21)(89,20.5)(95,20)
\qbezier(83,21)(79,15.5)(75,10)
\qbezier(83,21)(79,25.5)(75,30)


\put(71,8.0){$v_2$}
\put(83,22){$v_1$}
\put(96,19.0){$v_2$}
\put(71.0,31.0){$v_2$}
\put(93.5,42.0){$v_2$}
\put(85,13){$a_1$}
\put(70,20.0){$a_2$}

\put(76,17.2){$\mathbf{e}_1$}
\put(89,21.2){$\mathbf{e}_2$}
\put(78,26.6){$\mathbf{e}_3$}

\put(64,10){(\emph{b})}
\end{picture}

\vspace{-0.5cm} \caption{ \emph{a}) Graphene $\bG$; \quad \emph{b}) the fundamental graph $\bG_0$ of the graphene.} \label{ff.0.3}
\end{figure}
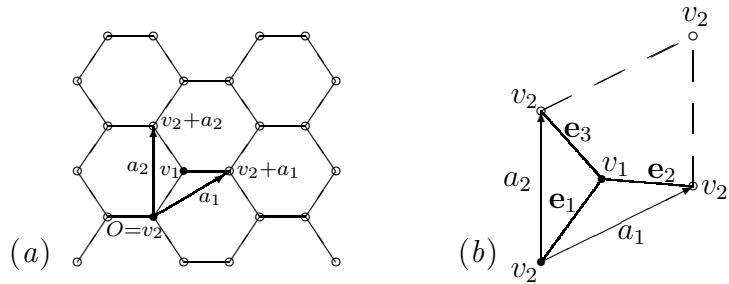

We describe the spectrum
\textbf{under the perturbation} of the graphene $\bG$ by adding an edge
to its fundamental graph.

\begin{theorem}[Perturbations of the graphene]
\lb{T10}
Let $\bG'$ be a  perturbed graph
obtained  from the hexagonal lattice $\bG$ by adding one edge (i.e., two oriented
edges $\be,\bar\be$ with the indices  $\t(\be)=-\t(\bar\be)=(\t_1,\t_2)$) to its fundamental graph $\bG_0$. Then the spectrum $\s(\D')$  of the Laplacian $\D'$ on
the perturbed graph $\bG'$ satisfies:

i) $\s(\D')=\s_{ac}(\D')$, i.e., the Laplacian on $\bG'$ has no flat bands.

ii) $\s(\D')=[-1,1]$ iff $\bG'$ is bipartite.

iii) Let $\bG'$ be non-bipartite. Then
\[
\lb{ge2} \s(\D')=[\l_2^-,0]\cup[\l_1^-,1],\qqq 0<\l_1^- \le {2\/5}\,,
\qqq if\ \t_1-\t_2\in 3\Z,
\]
\[
\lb{ge3} \s(\D')=[\l_2^-,\l_2^+]\cup[0,1], \qqq -{1\/5}\le\l_2^+<0,
\qqq if\ \t_1-\t_2\notin 3\Z,
\]
where
\[
\lb{ge1}
-1<\l_2^-\leq-\frac35\,.
\]
Moreover, the Lebesgue measure
of the spectrum satisfies
\[
\lb{ge4}
\begin{aligned}
|\s(\D')|\ge {6\/5}\,,\qqq {\rm if \qq } \t_1-\t_2\in 3\Z;\\
|\s(\D')|\ge {7\/5}\,,\qqq {\rm if \qq } \t_1-\t_2\notin 3\Z.
\end{aligned}
\]
\end{theorem}

\no \textbf{Remark.}
This theorem shows that if the perturbation
keeps the bipartition of the graph, then the spectrum of the operator does not change. If the added edge breaks the bipartition of the graph, then there appears a gap in the spectrum and the lower point of the spectrum moves to the right.

\

We present the plan of our paper.
In section 2 we prove Theorem \ref{pro2} about the decomposition of the Laplacian into a direct integral and prove some basic properties of the Floquet operator. In section 3 we estimate the Lebesgue measure of the spectrum in terms of geometric parameters of the graph in Theorem \ref{T1}. In sections 4 and 5 we consider the periodic graphs with one, two and three vertices in the fundamental graph. In particularly, we describe the spectrum of the Laplacian on the square lattice and the hexagonal lattice, perturbed by adding edges to its fundamental graphs. In Appendix we recall some well-known  properties of matrices and give some auxiliary statements, needed to prove our main results.

\section{Proof of Theorem \ref{pro2}}
\setcounter{equation}{0}

In Appendix we recall some well-known  properties of matrices, which will be used below. Also we prove that the fundamental graph of a periodic graph is finite and show that a bipartite periodic graph has a
bipartite fundamental graph. Properties of an edge index, needed to prove our main results, are given in Appendix.

\

\no {\bf Proof of Theorem \ref{pro2}.i)} Let $\mH=L^2\big(\T^2,{d\vt
\/(2\pi)^2}\,,\cH\big)=
 \int_{\T^2}^{\os}\cH\,{d\vt \/(2\pi)^2}
$ be  a constant fiber direct integral, where $\cH=\C^\nu$. It is
the Hilbert space of square integrable functions
$f:\T^2\rightarrow\mathbb{C}^\nu$, equipped with the norm
$$
\|f\|^2_{\mH}=\int_{\T^2}\|f(\vt )\|_{\C^\nu}^2\frac{d\vt
}{(2\pi)^2}\,.
$$
Let $\mS(V)$ be the set of all finitely supported functions $f\in
\ell^2(V)$. Recall that $v_1,\ldots,v_\nu$ are the vertices of the
fundamental graph $\G_0$. We identify the vertices of $\G_0$ with the vertices of the periodic graph $\G$ from the set $[0,1)^2$.
We define the operator
$U:\ell^2(V)\to\mH$ by
\[
 \lb{5001}
(Uf)_n(\vt )=\sum\limits_{p\in\mathbb{Z}^2}e^{-i\lan p,\vt\ran }f(v_n+p), \qqq
  (\vt ,n)\in \T^2\ts\N_\nu.
\]
Using standard arguments (see pp. 290--291 in \cite{RS78}), we will show that $U$ is well defined for $f\in\mS(V)$
and uniquely extendable to a unitary operator.

For  $f\in \mS(V)$ the sum \er{5001} is
clearly convergent. For such functions $f$ we compute
\begin{multline*}
\|Uf\|^2_{\mH}=\int_{\T^2}\|(Uf)(\vt )\|_{\C^\nu}^2{d\vt
\/(2\pi)^2}\\=
 \int_{\T^2}\sum_{n=1}^\nu\bigg(\sum\limits_{p\in\Z^2}e^{-i\lan p,\vt\ran }f(v_n+p)\bigg)
 \bigg(\sum\limits_{q\in\Z^2}e^{i\lan q,\vt\ran }\bar{f}(v_n+q)\bigg)
{d\vt \/(2\pi)^2}\\
=\sum_{n=1}^\nu\sum_{p\in\Z^2}\sum_{q\in\Z^2}\lt(
 f(v_n+p)\bar f(v_n+q)\int_{\T^2}e^{-i\lan p-q,\vt\ran }{d\vt \/(2\pi)^2}\rt)\\=
 \sum_{n=1}^\nu\sum_{p\in\Z^2}\big|f(v_n+p)\big|^2=\sum_{v\in V}|f(v)|^2=\|f\|_{\ell^2(V)}^2.
\end{multline*}
Here we have used the identity
$V=\big\{v_n+p: (n,p)\in \N_\nu\ts\Z^2\big\}$. Then $U$ is well defined and has a unique extension to an isometry.
To see that $U$ is onto $\mH$ we compute $U^*$.
We note that each vertex $v\in V$ can be written as
\begin{equation}\label{rep}
v=v_n+p \quad \textrm{ for some }\; (n,p)\in \N_\nu\ts\Z^2
\end{equation}
and this representation is unique. Let $g=\big(g_n(\vt):\T^2\ra\C\big)_{n\in\N_\nu}\in\mH$. Then
we define
$$
(U^*g)(v)=\int_{\T^2}e^{i\lan p,\vt\ran }g_n(\vt ){d\vt \/(2\pi)^2}\,,\qquad
 v\in V,
$$
where $p$ and $n$ are defined by (\ref{rep}). A direct computation shows that this is indeed the formula for
the adjoint of $U$. Moreover,
$$
\|U^*g\|^2_{\ell^2(V)}=\sum_{v\in
V}\big|(U^*g)(v)\big|^2=\sum_{n=1}^\nu\sum_{p\in\Z^2}\big|(U^*g)(v_n+p)\big|^2=
\sum_{n=1}^\nu\sum_{p\in\Z^2}\bigg|\int_{\T^2}e^{i\lan p,\vt\ran
}g_n(\vt ){d\vt \/(2\pi)^2}\bigg|^2
$$
$$
=\sum_{n=1}^\nu\int_{\T^2}\big|g_n(\vt )\big|^2{d\vt \/(2\pi)^2}=
 \int_{\T^2}\sum_{n=1}^\nu\big|g_n(\vt )\big|^2{d\vt \/(2\pi)^2}=\|g\|_\mH^2,
$$
where we have used the Parseval relation for the Fourier
series.

For $f\in \mS(V)$ and $j\in\N_\nu$ we obtain
\begin{multline}\label{ext1}
(U\D f)_j(\vt )=\sum_{p\in\Z^2}e^{-i\lan p,\vt\ran }(\D
f)(v_j+p) =\sum_{p\in\Z^2}e^{-i\lan p,\vt\ran }\frac1{\sqrt{\vk_j}}
\sum\limits_{(v_j+p,\,u)\in\cA}\frac1{\sqrt{\vk_u}}\,f(u)\\=\sum_{p\in\Z^2}e^{-i\lan p,\vt\ran
}\frac1{\sqrt{\vk_j}}
\sum\limits_{k=1}^\nu\sum\limits_{\be=
(v_j,v_k)\in\cA_0}\dfrac1{\sqrt{\vk_k}}\,f\big(v_k+p+\t (\be)\big)\\
=\sum\limits_{k=1}^\nu\sum\limits_{\be=
(v_j,v_k)\in\cA_0}\dfrac1{\sqrt{\vk_j\vk_k}}\,e^{i\lan\t(\be),\vt\ran }
\sum_{p\in\Z^2}e^{-i\lan p+\t(\mathbf{e}),\vt\ran }f\big(v_k+p+\t (\mathbf{e})\big)\\
=\sum\limits_{k=1}^\nu\sum\limits_{\mathbf{e}=
(v_j,v_k)\in\cA_0}\dfrac1{\sqrt{\vk_j\vk_k}}\,e^{i\lan\t
(\be),\vt\ran } (Uf)_k(\vt)=\sum\limits_{k=1}^\nu \D_{jk}(\vt)(Uf)_k(\vt),
\end{multline}
where $\D_{jk}(\vt)$ are defined by \er{l2.15}.
The identity (\ref{ext1}) yields $(U\D f)(\vt
)=\D(\vt )(Uf)(\vt )$, where $\D(\vt )=\{\D_{jk}(\vt
)\}$. Thus, we obtain
$$
 U\D U^{-1}=\frac1{(2\pi)^2}\int_{\T^2}^{\os}\D(\vt )\,d\vt.
$$
Using Proposition \ref{pro1}.i we
can write the entries $\D_{jj}$ of the matrix $\D(\vt)$ in the form (\ref{l2.15'}).
Thus, the statement i) has been proved.
\medskip

In order to prove  ii) -- vii), we discuss needed properties of the Floquet matrix $\D(0)$.

\begin{proposition}
\label{pp1}
The matrix $\D(0)$ is the Laplacian on the fundamental graph $\G_0$ and has the form
\begin{equation}\label{ll6}
\D(0)=\{\D_{jk}(0)\},\qqq \D_{jk}(0)=\frac{\vk_{jk}}{\sqrt{\vk_j\vk_k}}\,,\qqq \forall\,(j,k)\in\N_\n^2,
\end{equation}
where $\vk_{jk}\ge 1$ is the multiplicity
of the edge $(v_j,v_k)\in \cA_0$ and  $\vk_{jk}=0$ if
$(v_j,v_k)\notin \cA_0$. Moreover, they satisfy
\begin{equation}
\label{l.83} \vk_j=\sum_{k=1}^\n\vk_{jk}\ge1,\qquad \forall\,j\in \N_\nu.
\end{equation}
\end{proposition}

\no \textbf{Proof.} The identity (\ref{ll6}) for the matrix
$\D(0)$ is  obtained by a direct calculation of its entries using
the formula (\ref{l2.15}).
From the definition \er{DOL} we deduce that the Laplacian on $\G_0$ has the form \er{ll6}.
The definition of the degree ${\vk}_j$ of the vertex $v_j$
as the number of oriented edges starting at $v_j$ gives (\ref{l.83}). Since the graph $\G$ is connected, $\vk_j\geq1$.  \quad $\BBox$

\medskip

\no {\bf Proof of Theorem \ref{pro2}.ii -- vii.}
ii) Recall that we identify the vertices $v_1,\ldots,v_\nu$ of the fundamental graph $\G_0$ with the vertices of the periodic graph $\G$ from the set $[0,1)^2$
in the coordinate system with the origin $O$. Proposition \ref{pro1}.iv gives that for each $(v_j,v_k)\in\cA_0$
\begin{equation}\label{ind}
\t^{(1)}(v_j,v_k)=\t(v_j,v_k)+p_k-p_j,\qqq p_k=[v_k-b],\qqq p_j=[v_j-b],
\end{equation}
where $\t^{(1)}(v_j,v_k)$ is the index of the edge $(v_j,v_k)$ in the coordinate system with the origin $O_1$, $b=\overrightarrow{OO}_1$.
Using (\ref{ind}) we rewrite the entries of $\D^{(1)}(\vt)$ defined by (\ref{l2.15}) in the form
$$
\D_{jk}^{(1)}(\vt
)=\dfrac1{\sqrt{\vk_j\vk_k}}\sum\limits_{\be=(v_j,\,v_k)\in\cA_0}e^{i\lan\t^{(1)}
(\be),\vt\ran }\\=\dfrac{e^{i\lan p_k-p_j,\vt\ran }}{\sqrt{\vk_j\vk_k}}\sum\limits_{\be=(v_j,\,v_k)\in\cA_0}e^{i\lan\t
(\be),\,\vt\ran }=e^{i\lan p_k-p_j,\,\vt\ran }\D_{jk}(\vt
).
$$
We define the diagonal $\nu\times\nu$ matrix
$$
\cU(\vt )=\mathrm{diag}  \left(
\begin{array}{ccc}
   e^{-i\lan p_1,\,\vt\ran}, &  \ldots,& e^{-i\lan p_\nu,\,\vt\ran}
\end{array}\right),\qquad \forall\, \vt \in\T^2.
$$
A direct calculation yields
$$
\cU(\vt )\,\D(\vt )\,\cU^{-1}(\vt
)=\D^{(1)}(\vt ),\qqq \forall\, \vt \in\T^2.
$$
Thus, for each $\vt \in\T^2$ the matrices $\D^{(1)}(\vt )$ and $\D(\vt )$ are unitarily equivalent.

iii) Since $\D(0)$ is the Laplacian on $\G_0$, $\|\D(0)\|=1$.
From this fact and the formula (\ref{l2.15}) it follows that
the entries of each matrix $\D(\vt )$, $\vt \in\T^2$, satisfy
\[
\lb{vv1}
|\D_{jk}(\vt )|\leq\D_{jk}(0)\leq1, \qquad \forall\,(j,k)\in\N_\nu^2.
\]
Then, Proposition \ref{MP}.i  implies that the spectral radius $\r\big(\D(\vt )\big)$ satisfies
$$
\r\big(\D(\vt )\big)\leq\r\big(\D(0)\big)=\|\D(0)\|=1,
$$
which yields $\s\big(\D(\vt
)\big)\ss[-1,1]$ for each $\vt
\in\T^2$.

iv) This statement follows immediately from (\ref{l2.15}) and the definition of the bridge.

v) We prove by the contradiction. Let the point 1 be an eigenvalue of the Laplacian
on a graph $\G$. Then there exists an eigenfunction $0\neq f\in\ell^2(V)$ with eigenvalue 1 and with a finite support $\cB\ss V$ (see Theorem 3.2 in \cite{HN09}). Let $M=\max\limits_{v\in\cB}\frac{f(v)}{\sqrt{\vk_v}}=\frac{f(\tilde v)}{\sqrt{\vk_{\tilde v}}}$ for some $\tilde v\in\cB$.
Thus, we have
$$
f(\tilde v)=\big(\D f\big)(\tilde v)=\frac1{\sqrt{\vk_{\tilde v}}} \sum\limits_{(\tilde v,\,u)\in\cA_0}\frac1{\sqrt{\vk_u}}\,f(u)\leq\frac1{\sqrt{\vk_{\tilde v}}} \sum\limits_{(\tilde v,\,u)\in\cA_0}\frac1{\sqrt{\vk_{\tilde v}}}\,f(\tilde v)=M\sqrt{\vk_{\tilde v}}=f(\tilde v).
$$
We conclude that the inequality has to be an equality, and therefore
$$
\frac{f(u)}{\sqrt{\vk_u}}=\frac{f(\tilde v)}{\sqrt{\vk_{\tilde v}}}\,, \qqq \forall\,u\sim
\tilde v.
$$
Repeating this
argument until we reach a vertex from $V\setminus\cB$, we conclude that $f=0$.
We have a contradiction.

vi) We prove by the contradiction. Let all entries $\D_{jk}(\cdot)$, $1\leq j\leq k\leq\n$, be constant. Due to the self-adjointness of the matrix $\D(\vt)$ all its entries are constants.
Then the point 1 is an eigenvalue of $\D(\vt)=\D(0)$ for each $\vt\in\T^2$, i.e., the point 1 is an eigenvalue of $\D$ with infinite multiplicity.
This contradicts item v).

vii) Let $\G_0$ be a bipartite fundamental graph with the parts
$V_1=\{v_1,\ldots,v_{k}\}$ and $V_2=\{v_{k+1},\ldots,v_\nu\}$. Since vertices from the same part of
$\G_0$ are not adjacent, each matrix $\D(\vt )$, $\vt \in\T^2$, has
the form
\begin{equation}\label{m1}
\D(\vt )=\left(
\begin{array}{cc}
  \O_{kk} & A(\vt )\\[10pt]
  A^\ast(\vt ) & \O_{\nu-k,\nu-k}
\end{array}\right),
\end{equation}
for some $k\ts(\n-k)$ matrix $A(\vt)$. Here $\O_{jk}$ is the zero $j\ts k$ matrix. We define the matrix
$$
\cU=\left(
\begin{array}{cc}
  I_k & \O_{k,\nu-k}\\[10pt]
  \O_{\nu-k,k} & -I_{\nu-k}
\end{array}\right)=\cU^{-1},
$$
where $I_k$ is the identity $k\ts k$ matrix. By a direct calculation,
one can verify that
$$
\cU\D(\vt )\cU^{-1}=-\D(\vt ),
$$
which yields that $\s(\D(\vt ))$ is symmetric with respect to 0.

Conversely, since $\D(0)$ is the Laplacian on $\G_0$ and its spectrum is symmetric with respect to 0,
$\G_0$ is bipartite due to the main property 4) of the Laplacian.
\quad
$\BBox$ \medskip

In the following theorem we show unitary equivalence of Laplacians on graphs with multiple indices.

\begin{theorem}\label{tt2}
Let $\G_0=(V_0,\cE_0)$ and $\G'_0=(V'_0,\cE'_0)$ be fundamental
graphs of periodic graphs $\G$ and $\G'$, respectively, satisfying
the following conditions:

1) $\G_0$ and $\G'_0$ are isomorphic, i.e., there exists a bijection $\vp:V_0\ra V'_0$ that preserves the adjacency of vertices;

2) there exists $n\in\N$ such that for any $\mathbf{e}=(u,v)\in\cA_0$
\begin{equation}\label{rat}
\t ({\bf e})=n\,\t ({\bf e}'), \quad \textrm{where} \quad \mathbf{e}'=\big(\varphi(u),\varphi(v)\big)\in\cA'_0.
\end{equation}
Then the Laplace operators on $\G$ and $\G'$ are unitarily equivalent.
\end{theorem}
\no \textbf{Proof.} Denote the Laplace operators on
$\G$ and $\G'$ by $\D$ and $\D'$, respectively, and the Floquet matrices for
$\G$ and $\G'$ by
$\D(\vt)$ and $\D'(\vt)$, $\vt
\in\T^2$,  respectively. Since $\G_0$ and $\G'_0$ are isomorphic and the indices of their edges satisfy the identity (\ref{rat}),
$$
\D(\vt)=\D'(n\vt), \qqq \forall\,\vt\in\T^2.
$$
Then using (\ref{raz}) we obtain
$$
U\D U^{-1}=\frac1{(2\pi)^2}\int^\oplus_{\T^2}\D(\vt )d\vt =
\frac1{(2\pi)^2}\int^\oplus_{\T^2}\D'(n\vt )d\vt .
$$
If we make the change of variables $\vt'=n\vt$, $\vt'\in\T'^2=\R^2/(2\pi n\Z)^2=[-\pi n;\pi n]^2$, then we rewrite the last
identity in the form
$$
U\D U^{-1}= \frac1{(2\pi)^2}\int^\oplus_{\T^2}\D'(n\vt )d\vt = \frac{1}{(2\pi
n)^2}\int^\oplus_{\T'^2}\D'(\vt ')d\vt'=
\frac1{(2\pi)^2}\int^\oplus_{\T^2}\D'(\vt')d\vt '=U'\D'U'^{-1}
$$
for some unitary operator $U'$.
Thus, $\D$ and $\D'$ are unitarily equivalent. \qquad $\BBox$

\section{Spectral analysis of Laplacians}
\setcounter{equation}{0}

Below we need the following representation of the Floquet matrix $\D(\vt)$, $\vt\in\T^2$:
\[
\label{eq.1}
\D(\vt)=\D_0+\wt\D(\vt),\qqq \D_0={1\/(2\pi)^2}\int_{\T^2}\D(\vt
)d\vt.
\]
From \er{eq.1}, \er{l2.15} it follows that the entries of the matrix $\wt\D(\vt)=\{\wt\D_{jk}(\vt)\}$ have the form
\[
\label{tl2.15}
\wt\D_{jk}(\vt )=
{1\/\sqrt{\vk_j\vk_k}}\sum\limits_{{\bf e}=(v_j,v_k)\in{\cA}_0\atop\t
(\be)\neq0 }e^{\,i\lan\t
({\bf e}),\,\vt\ran }.
\]

\no {\bf Proof of Theorem \ref{T100}.}
i)-ii) Due to \er{eq.1} we have $\D(\vt)=\D_0+\wt\D(\vt)$, where  the matrix $\wt\D(\vt)$, $\vt \in \T^2$, has the form
$$
\wt\D(\vt)=\diag\big(\wt\D_{11},\ldots,\wt\D_{\n\n}\big)(\vt).
$$
The identity \er{tl2.15} for $\wt\D_{jj}$ has the form
$$
\wt\D_{jj}(\vt )=
{1\/\vk_j}\sum\limits_{{\bf e}=(v_j,v_j)\in{\cA}_0\atop\t
(\be)\neq0 }\cos\lan\t
({\bf e}),\,\vt\ran,
$$
which yields
$$
\wt\D_{jj}(\vt_0)\leq\wt\D_{jj}(\vt)\leq\wt\D_{jj}(0),\qqq \forall\,j\in\N_\n.
$$
Then $\wt\D(\vt_0)\le\wt\D(\vt)\le\wt\D(0)$ and we have
\[
\D(\vt_0)=\D_0+\wt\D(\vt_0)\le\D(\vt )=\D_0+\wt\D(\vt)\le \D_0+\wt\D(0)= \D(0).
\]
Then Proposition \ref{MP}.ii gives
\[
\l_n(\vt_0)\le\l_n(\vt )\le \l_n(0),\qquad \forall\,
(\vt,n)\in\T^2\ts\N_\n.
\]
Thus,
$\l_n^+=\max\limits_{\vt
\in\T^2}\l_n(\vt )=\l_n(0)$,  $\l_n^-=\min\limits_{\vt
\in\T^2}\l_n(\vt )=\l_n(\vt_0)$. The last statement of the item i) follows from the definition of flat bands.

iii) Since $\G$ is bipartite, the spectrum of the Laplacian on $\G$ is symmetric with respect to zero. From item i) it follows that $\l_1(0)\geq\ldots\geq\l_\n(0)$ are the upper endpoints of the spectral bands. Then $-\l_1(0)\leq\ldots\leq-\l_\n(0)$ are the lower endpoints of the spectral bands. Thus, the endpoints of the spectral bands $\l_n^\pm$, $n\in\N_\n$, are the eigenvalues of the matrix $\pm\D(0)$.
\qq  \BBox

\begin{lemma}\label{t8.1'}
Let all edges $(v_j,v_k)\in\cA_0$, $1\leq j,k\leq\nu-1$, of the fundamental graph $\G_0$ have zero indices. Then

i) The Floquet matrix $\D(\vt)$ has the form
\[
\lb{ww.10}
\Delta(\vt )=\left(
\begin{array}{cc}
  A & y(\vt ) \\
  y^\ast(\vt ) & a(\vt )
\end{array}\right),
\]
where for each $\vt\in\T^2$ the entry $y(\vt )\in\C^{\nu-1}$ is a vector
and $a(\vt )$ is a real number, $A$ is a self-adjoint $(\nu-1)\ts(\nu-1)$ matrix not depending on $\vt$.

ii) If $A$ has
an eigenvalue $\m$ with multiplicity $\ge 2$, then $\m$ is  a flat
band of the Laplacian $\D$ on the periodic graph $\G$.

\end{lemma}

\no {\bf Proof.} i) This follows from Theorem \ref{pro2}.iv and the self-adjointness of $\D(\vt)$.

ii) Let the matrix $A$ have a multiple eigenvalue $\mu$. Due to Proposition \ref{MP}.iii, there exists an
eigenvalue $\l(\vt )$ of $\D(\vt )$
satisfying
$
\mu\leq\l(\vt )\leq\mu $ for all $\vt \in \T^2$,
which yields
$\l(\cdot)=\mu=\mathrm{const}$, i.e., $\mu$ is a flat band of the Laplacian $\D$ on $\G$. \quad $\BBox$

\

\setlength{\unitlength}{1.0mm}
\begin{figure}[h]
\centering
\unitlength 1mm 
\linethickness{0.4pt}
\ifx\plotpoint\undefined\newsavebox{\plotpoint}\fi 
\begin{picture}(120,50)(0,0)

\put(10,10){\line(1,0){40.00}}
\put(10,30){\line(1,0){40.00}}
\put(10,50){\line(1,0){40.00}}
\put(10,10){\line(0,1){40.00}}
\put(30,10){\line(0,1){40.00}}
\put(50,10){\line(0,1){40.00}}

\put(10,10){\circle{1}}
\put(30,10){\circle{1}}
\put(50,10){\circle{1}}

\put(10,30){\circle{1}}
\put(30,30){\circle{1}}
\put(50,30){\circle{1}}

\put(10,50){\circle{1}}
\put(30,50){\circle{1}}
\put(50,50){\circle{1}}

\put(10,10){\vector(1,0){20.00}}
\put(10,10){\vector(0,1){20.00}}

\put(7.5,7.5){$\scriptstyle v_\n$}
\put(1,29){$\scriptstyle v_\n+a_2$}
\put(29,7){$\scriptstyle v_\n+a_1$}
\put(20,8){$\scriptstyle a_1$}
\put(6,20){$\scriptstyle a_2$}
\put(14,26.5){$\scriptstyle v_1$}
\put(23.8,17){$\scriptstyle v_{\nu-2}$}

\put(24.0,12){$\scriptstyle v_{\nu-1}$}
\put(20,26.5){$\scriptstyle v_2$}

\put(10,10){\line(2,3){10.00}}
\put(10,10){\line(3,2){15.00}}
\put(10,10){\line(1,3){5.00}}
\put(10,10){\line(3,1){15.00}}
\put(20.5,21){$\ddots$}
\put(15,25){\circle{1}}
\put(20,25){\circle{1}}
\put(25,15){\circle{1}}
\put(25,20){\circle{1}}

\put(30,10){\line(2,3){10.00}}
\put(30,10){\line(3,2){15.00}}
\put(30,10){\line(1,3){5.00}}
\put(30,10){\line(3,1){15.00}}
\put(40.5,21){$\ddots$}
\put(35,25){\circle{1}}
\put(40,25){\circle{1}}
\put(45,15){\circle{1}}
\put(45,20){\circle{1}}

\put(10,30){\line(2,3){10.00}}
\put(10,30){\line(3,2){15.00}}
\put(10,30){\line(1,3){5.00}}
\put(10,30){\line(3,1){15.00}}
\put(20.5,41){$\ddots$}
\put(15,45){\circle{1}}
\put(20,45){\circle{1}}
\put(25,35){\circle{1}}
\put(25,40){\circle{1}}

\put(30,30){\line(2,3){10.00}}
\put(30,30){\line(3,2){15.00}}
\put(30,30){\line(1,3){5.00}}
\put(30,30){\line(3,1){15.00}}
\put(40.5,41){$\ddots$}
\put(35,45){\circle{1}}
\put(40,45){\circle{1}}
\put(45,35){\circle{1}}
\put(45,40){\circle{1}}
\put(-4,10){(\emph{a})}
\put(64,25){(\emph{b})}

\multiput(80,45)(4,0){5}{\line(1,0){2}}
\multiput(100,25)(0,4){5}{\line(0,1){2}}
\put(80,25){\line(1,0){20.00}}
\put(80,25){\line(0,1){20.00}}
\put(80,25){\circle{1}}
\put(100,25){\circle{1}}

\put(80,45){\circle{1}}
\put(100,45){\circle{1}}

\put(80,25){\vector(1,0){20.00}}
\put(80,25){\vector(0,1){20.00}}

\put(77.0,23.0){$\scriptstyle v_\n$}
\put(76.5,45){$\scriptstyle v_\n$}
\put(101,23){$\scriptstyle v_\n$}
\put(101,45){$\scriptstyle v_\n$}
\put(90,23){$\scriptstyle a_1$}
\put(76,35){$\scriptstyle a_2$}

\put(84,41.5){$\scriptstyle v_1$}
\put(93.8,32.5){$\scriptstyle v_{\nu-2}$}

\put(94.0,27){$\scriptstyle v_{\nu-1}$}
\put(90,41.5){$\scriptstyle v_2$}

\put(80,25){\line(2,3){10.00}}
\put(80,25){\line(3,2){15.00}}
\put(80,25){\line(1,3){5.00}}
\put(80,25){\line(3,1){15.00}}
\put(90.5,36){$\ddots$}
\put(85,40){\circle{1}}
\put(90,40){\circle{1}}
\put(95,30){\circle{1}}
\put(95,35){\circle{1}}
\put(64,10){(\emph{c})}

\put(80,10){\line(1,0){30.00}}
\put(80,9){\line(0,1){2.00}}
\put(110,9){\line(0,1){2.00}}
\put(90,9){\line(0,1){2.00}}
\put(100,9){\line(0,1){2.00}}
\put(95,10){\circle*{1}}

\put(80,9.8){\line(1,0){10.00}}
\put(80,10.2){\line(1,0){10.00}}

\put(100,9.8){\line(1,0){10.00}}
\put(100,10.2){\line(1,0){10.00}}

\put(103,12){$\s_1$}
\put(83,12){$\s_2$}
\put(94,12){$\m$}
\put(77,6){$\scriptstyle-1$}
\put(94.3,6){$\scriptstyle0$}
\put(109.3,6){$\scriptstyle1$}
\end{picture}
\vspace{-0.5cm}
\caption{\emph{a}) The periodic graph $\G$;\quad \emph{b}) the fundamental graph $\G_0$; only two loops in the vertex $v_\n$ are bridges;\quad \emph{c}) the spectrum of the Laplacian.} \label{ff.10}
\end{figure}
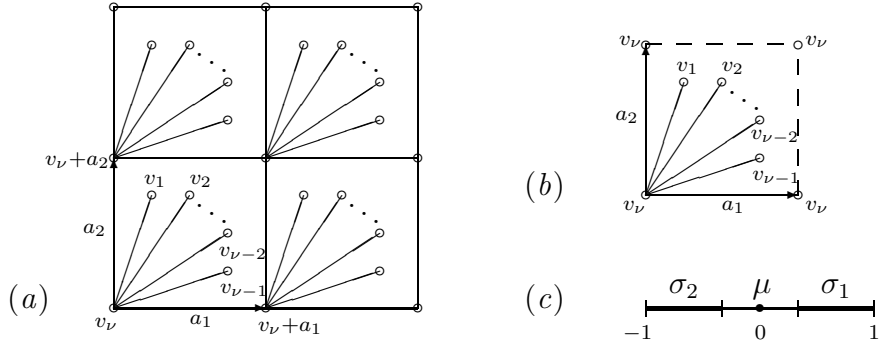

\begin{proposition}\lb{TG1}
Let $\G_0$ be obtained from the fundamental graph $\bS_0$ of the square
lattice $\bS$ by adding $\n-1$ ($\n\geq2$) vertices $v_1,\ldots,v_{\n-1}$ and $\n-1$ unoriented edges $(v_1,v_\n),\ldots,(v_{\n-1},v_\n)$  with zero indices (see Figure~\ref{ff.10}), $v_\n$ is a single vertex of $\bS_0$.
Then the spectrum of the Laplacian on $\G$  has the form
\[
\lb{acs1}
\s(\D)=\s_{ac}(\D)\cup\s_{fb}(\D),\qqq \s_{fb}(\D)=\{0\},
\]
where the point $0$ is a flat band of multiplicity $\n-2$ and
the absolutely continuous part $\s_{ac}(\D)$ has only two spectral bands $\s_1$ and $\s_2$ given by
\[
\lb{acs2}
\s_{ac}(\D)=\s_2\cup\s_1, \qqq \s_1=-\s_2=[{\textstyle\frac{\n-1}{\n+3}}\,,1].
\]

\end{proposition}
\no {\bf Proof.} The fundamental graph $\G_0$ consists of $\n\geq2$ vertices $v_1,v_2,\ldots,v_\n$; $\n-1$ unoriented edges ${(v_1,v_\n),\ldots,(v_{\n-1},v_\n)}$ with zero indices and 2 unoriented loops in the vertex $v_\n$ with the indices $(\pm1,0)$, $(0,\pm1)$. All bridges of $\G_0$ are loops and the graph $\G$ is bipartite. Then, by Theorem \ref{T100}.iii, the spectrum of the Laplacian is completely defined by the eigenvalues of $\D(0)$. According to (\ref{ll6}) we have
$$
\Delta(0)=\left(
\begin{array}{cccc}
  0 & 0  & \ldots & \frac1{\sqrt{\nu+3}} \\[2pt]
  0 & 0 & \ldots & \frac1{\sqrt{\nu+3}} \\[2pt]
  \ldots & \ldots & \ldots & \ldots  \\[2pt]
  \frac1{\sqrt{\nu+3}} & \frac1{\sqrt{\nu+3}} & \ldots & \frac{4}{\nu+3}  \\
\end{array}\right).
$$
Using the formula \er{det} we obtain
$$
\det\big(\D(0)-\l I_\n\big)=\frac{(-1)^\n\,\l^{\n-2}}{\n+3}
\big((\n+3)\l^2-4\l+1-\n\big).
$$
Then the eigenvalues of $\D(0)$ have the form
$$
\l_1(0)=1,\qqq \l_2(0)=\ldots=\l_{\n-1}(0)=0,\qqq \l_\n(0)=-\frac{\n-1}{\n+3}\,.
$$
Thus, due to Theorem \ref{T100}.iii, the spectrum of the Laplacian on $\G$  has the form \er{acs1}, \er{acs2}.
\BBox

\

\setlength{\unitlength}{0.8mm}
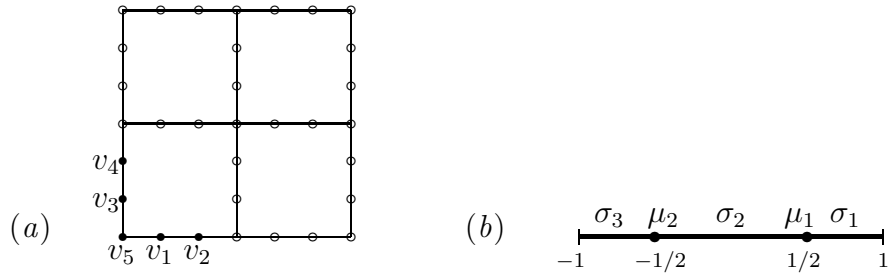
\begin{figure}[h]
\centering

\unitlength 1mm 
\linethickness{0.4pt}
\ifx\plotpoint\undefined\newsavebox{\plotpoint}\fi 
\begin{picture}(90,40)(0,0)
\put(-5,10){(\emph{a})}
\put(10,10){\line(1,0){30.00}}
\put(10,25){\line(1,0){30.00}}
\put(10,40){\line(1,0){30.00}}
\put(10,10){\line(0,1){30.00}}
\put(25,10){\line(0,1){30.00}}
\put(40,10){\line(0,1){30.00}}

\put(10,10){\circle*{1}}
\put(15,10){\circle*{1}}
\put(20,10){\circle*{1}}
\put(25,10){\circle{1}}
\put(30,10){\circle{1}}
\put(35,10){\circle{1}}
\put(40,10){\circle{1}}

\put(10,25){\circle{1}}
\put(15,25){\circle{1}}
\put(20,25){\circle{1}}
\put(25,25){\circle{1}}
\put(30,25){\circle{1}}
\put(35,25){\circle{1}}
\put(40,25){\circle{1}}

\put(10,40){\circle{1}}
\put(15,40){\circle{1}}
\put(20,40){\circle{1}}
\put(25,40){\circle{1}}
\put(30,40){\circle{1}}
\put(35,40){\circle{1}}
\put(40,40){\circle{1}}

\put(10,15){\circle*{1}}
\put(10,20){\circle*{1}}
\put(10,30){\circle{1}}
\put(10,35){\circle{1}}

\put(25,15){\circle{1}}
\put(25,20){\circle{1}}
\put(25,30){\circle{1}}
\put(25,35){\circle{1}}

\put(40,15){\circle{1}}
\put(40,20){\circle{1}}
\put(40,30){\circle{1}}
\put(40,35){\circle{1}}
\put(8.0,7.0){$v_5$}
\put(13.0,7.0){$v_1$}
\put(18.0,7.0){$v_2$}
\put(6.0,14.0){$v_3$}
\put(6.0,19.0){$v_4$}

\put(55,10){(\emph{b})}
\put(70,10){\line(1,0){40.00}}
\put(70,9.8){\line(1,0){40.00}}
\put(70,10.2){\line(1,0){40.00}}

\put(70,9){\line(0,1){2.00}}
\put(110,9){\line(0,1){2.00}}

\put(80,10){\circle*{1.5}}
\put(100,10){\circle*{1.5}}

\put(103,12){$\s_1$}
\put(88,12){$\s_2$}

\put(72,12){$\s_3$}
\put(79,12){$\m_2$}
\put(97,12){$\m_1$}
\put(67,6){$\scriptstyle-1$}
\put(97.3,6){$\scriptstyle1/2$}
\put(77.3,6){$\scriptstyle-1/2$}
\put(109.3,6){$\scriptstyle1$}
\end{picture}
\vspace{-0.5cm}
\caption{\emph{a}) Graph $\G$ obtained by adding two vertices on all edges of the square lattice $\bS$;\quad \emph{b}) the spectrum of the Laplacian.} \label{f.10}
\end{figure}

\begin{proposition}\lb{TG2}
Let $\G$ be the graph obtained from the square
lattice $\bS$ by adding $N$ vertices on each edge of $\bS$
(for $N=2$ see Figure~\ref{f.10}). Then the fundamental graph of $\G$ has $\n=2N+1$ vertices and the spectrum of the Laplacian on $\G$ has the form
\[
\lb{acs3}
\s(\D)=\s_{ac}(\D)\cup\s_{fb}(\D),
\]
where the absolutely continuous part $\s_{ac}(\D)=[-1,1]$ consists of $N+1$ non-degenerate spectral bands $\s_1,\ldots,\s_{N+1}$,
\[
\lb{acs4}
\s_{fb}(\D)=\Big\{\cos{\pi n\/ N+1}\;: \;n=1,\ldots,N\Big\},
\]
where each flat band is simple. There are no other flat bands.
\end{proposition}
\no {\bf Proof.}
The case $N=1$ will be considered in Corollary \ref{TC1}.

Let $N\geq2$. The fundamental graph $\G_0$ has $\n=2N+1$ vertices. The matrix $\D(\vt)$ is given by (\ref{ww.10}), where the $2N\ts2N$ matrix $A$ and the vector $y(\vt)$ have the form
\[
\lb{sl0}
\begin{aligned}
A=\left(
\begin{array}{cc}
  A_N & \mathbb{O}_{NN} \\
  \mathbb{O}_{NN} & A_N
\end{array}\right),\qqq
A_N={1\/2}\left(
\begin{array}{cccc}
  0 & 1 & 0 &\ldots \\
  1 & 0 & 1 &\ldots\\
  0 & 1 & 0 &\ldots \\
   \ldots & \ldots & \ldots &\ldots \\
\end{array}\right),\\
y(\vt)=\big(y_1(\vt_1),y_2(\vt_2)\big)^T,\qqq y_s(\vt_s)=\frac1{2\sqrt{2}}\big(1,0,\ldots,e^{i\vt_s}\big)^T,\qqq s=1,2,
\end{aligned}
\]
$a(\vt)=0$.
According to \er{eig} all eigenvalues of $A_N$ have the form
\begin{equation}\label{ww.17}
\mu_n=\cos\frac{\pi n}{N+1}\,,\qquad n=1,\ldots, N,
\end{equation}
and they are different. Then the matrix $A$ has $N$ different eigenvalues $\m_n$ of the multiplicity 2. Thus, due to Lemma \ref{t8.1'}.ii, the Laplacian on $\G$ has at least $N=[\frac\n2]$ flat bands $\m_n, n\in \N_N$.

We describe $\s_{ac}(\D)$. Identity \er{det} yields
\[
\lb{sl1}
\det\big(\D(\vt)-\l I_\n\big)=\det\big(A-\l I_{2N}\big)\big(-\l-y^\ast(\vt)
(A-\l I_{2N})^{-1}y(\vt)\big).
\]
By a direct calculation we get
\[
\lb{sl2}
\det\big(A-\l I_{2N}\big)=\cD^2_N(\l),\qqq \cD_N(\l)=\det\big(A_N-\l I_N\big),
\]
\[
\lb{sl3}
(A-\l I_{2N})^{-1}=\left(
\begin{array}{cc}
  B & \mathbb{O}_{NN} \\
  \mathbb{O}_{NN} & B
\end{array}\right), \qqq B=(A_N-\l I_N)^{-1}\,,
\]
\[
\lb{sl5}
y^\ast(\vt)
(A-\l I_{2N})^{-1}y(\vt)=y_1^\ast(\vt_1)
By_1(\vt_1)+y_2^\ast(\vt_2)
By_2(\vt_2),
\]
\[\lb{sl16}
y_s^\ast(\vt_s)
By_s(\vt_s)=\frac1{4\cD_N(\l)}\Big(\cD_{N-1}(\l)+
\Big(-\frac12\Big)^{N-1}\cos\vt_s\,\Big),\qqq s=1,2.
\]
Substituting \er{sl2}, \er{sl5} and \er{sl16}  into \er{sl1} we obtain
\[
\lb{sl4}
\det\big(\D(\vt)-\l I_\n\big)=\cD_N(\l)\Big(-\l\cD_N(\l)-
\frac12\,\cD_{N-1}(\l)-
\Big(-\frac12\Big)^{N+1}A(\vt)
\Big),
\]
where
$$
A(\vt)=\cos\vt_1+\cos\vt_2.
$$
From the identity \er{rel} it follows that
$$
\cD_N(\l)=-\l\,\cD_{N-1}(\l)-\frac14\,\cD_{N-2}(\l).
$$
Using this formula we rewrite \er{sl4} in the form
\[
\lb{sl15}
\det\big(\D(\vt)-\l I_\n\big)=\cD_N(\l)\Big(\cD_{N+1}(\l)-
\frac14\,\cD_{N-1}(\l)-
\Big(-\frac12\Big)^{N+1}A(\vt)
\Big).
\]
The identity \er{det1} gives
\[
\lb{sl6}
\cD_N(\l)=\frac{(-1)^N}{2^N}\cdot\frac{\sin(N+1)\vp}{\sin\vp}\,,
\]
where $\vp$ is given by
$$
-\lambda=2\cos\vp.
$$
The formulas \er{sl15} and \er{sl6}  give
\begin{multline}
\lb{sl7}
\det\big(\D(\vt)-\l I_\n\big)=\cD_N(\l)\bigg(\frac12\bigg)^{N+1}\left(
\frac{\sin(N+2)\vp}{\sin\vp}-\frac{\sin N\vp}{\sin\vp}-A(\vt)
\right)\\=\cD_N(\l)\bigg(\frac12\bigg)^{N+1}\big(2\cos(N+1)\vp-A(\vt)\big).
\end{multline}
Then the eigenvalues of the matrix $\D(\vt)$ are defined by
$$
\cD_N(\l)=0,\qqq 2\cos(N+1)\vp=A(\vt).
$$
The first identity gives all flat bands of $\D$ defined by \er{ww.17}. From the second identity and the fact that the range of the function $A(\vt)$  is $[-2,2]$  it follows that
for any $\l\in[-1,1]$ there exists $\vt\in\T^2$ such that
\[
\lb{sl9}
2\cos(N+1)\vp=A(\vt),
\]
i.e., each $\l\in[-1,1]$ is an eigenvalue
of $\D(\vt)$ for some $\vt\in\T^2$. Thus, $\s_{ac}(\D)=[-1,1]$.\qq \BBox

\

In order to prove Theorem \ref{T1}, we need the following lemma.

\begin{lemma}\lb{TT1}
Let $V=\{V_{jk}\}$ be a self-adjoint $\n\ts \n$  matrix such that $\sum\limits_{j,k=1}^\n|V_{jk}|<\iy$. Then the following estimate holds
true:
\[
\begin{aligned}
\lb{V}
&-B\le V\le B,\\
{\rm where}\qqq   & B=\diag \{B_1,\ldots, B_\n\},\qqq
B_j=\sum_{k=1}^\n |V_{jk}|.
\end{aligned}
\]
\end{lemma}
\no {\bf Proof.} Let $a=(j,k)$, $j,k=1,\ldots,\n$. We have
the identity
\[
\lb{b}
\begin{aligned}
V=\sum_{a=(j,k)} I(a)V_{a}=\sum_{a=(j,k)}Q(a),\qqq \qqq
Q(a)={1\/2}\big(I(a)V_{a}+I^\ast(a)\bar V_{a}\big),
\end{aligned}
\]
where the matrix $I(a)=\{I_{mn}(a)\}$, $a=(j,k)$ is given by
$$
I_{mn}(a)=\ca 1, & {\rm if} \ (m,n)=a\\
             0,  & {\rm if} \ (m,n)\ne a \ac.
$$
For any $a=(j,k)$ the following estimate holds true:
\[
\lb{Q}
\begin{aligned}
Q(a)={1\/2}\big(I(a)V_{a}+I^\ast(a)\bar V_{a}\big)\le |Q(a)|,\qqq\\
{\rm where }\qqq |Q(a)|=\sqrt{Q(a)\,Q^\ast(a)}={|V_{a}|\/2}\big(I(j,j)+I(k,k)\big).
\end{aligned}
\]
Summing \er{Q}, we obtain
$$
V=\sum_{a=(j,k)}Q(a)\le
\sum_{a=(j,k)}|Q(a)|=\sum_{a=(j,k)}{1\/2}\,|V_{a}|\big(I(j,j)+I(k,k)\big)=
\sum_{j=1}^\n I(j,j)\sum_{k=1}^\n|V_{jk}|=B,
$$
which yields $V\le B$. Moreover, this
yields $-B\le V$, since $-V\le B$. \qq \BBox

\

\

\no {\bf Proof of Theorem \ref{T1}.}
i) For each $\vt\in\T^2$ the estimate \er{V} yields
\[
\lb{eq.2}
-B(\vt) \le \wt \D(\vt )\le B(\vt),\qqq B(\vt )=\diag \{B_1(\vt ),\ldots, B_\n(\vt )\},
\]
\[
\lb{eq.2'}
B_j(\vt )=\sum_{k=1}^\n |\wt\D_{jk}(\vt)|.
\]
For all $(j,k,\vt)\in\N_\nu^2\ts\T^2$ it follows from \er{tl2.15} that
\[
\lb{wvv1}
|\wt\D_{jk}(\vt )|\leq\wt\D_{jk}(0))={b_{jk}\/\sqrt{\vk_j\vk_k}}.
\]
The estimate \er{wvv1} yields that the entries $B_j(\vt)$ defined by \er{eq.2'} satisfy
$$
B_j(\vt )\le B_j(0),\qqq \forall\,(j,\vt)\in\N_\nu\ts\T^2,
$$
and then
$$
B(\vt)\leq B(0),\qqq \forall\,\vt\in\T^2.
$$
Using this estimate we rewrite \er{eq.2} in the form
\[\label{eq'.2}
-B(0) \le \wt \D(\vt )\le B(0),\qq B(0)=\diag \{B_1(0),\ldots, B_\n(0)\},\qq
B_j(0)=\sum_{k=1}^\n \wt\D_{jk}(0).
\]
We use some arguments from \cite{Ku10}. Combining \er{eq.1} and \er{eq'.2}, we
obtain
\[
\D_0-B(0)\le \D(\vt )\le \D_0+B(0).
\]
Thus, Proposition \ref{MP}.ii implies
\[
\l_n(\D_0-B(0))\leq\l_n^-\le \l_n(\vt)\leq\l_n^+\le \l_n(\D_0+B(0)),
\qqq \forall \ \vt \in \T^2
\]
and then
\[
\lb{009}
 \big|\s(\D)\big|\le\sum_{n=1}^{\nu}(\l_n^+-\l_n^-)\leq
 \sum_{n=1}^\nu\big(\l_n(\D_0+B(0))-\l_n(\D_0-B(0))\big)=
 2\Tr B(0).
\]
Identity in
 \er{wvv1} and relations  \er{eq'.2} and \er{009} give
$$
\big|\s(\D)\big|\le2\Tr B(0)=2\sum_{j=1}^\n B_j(0)=2\sum_{j,k=1}^\n \wt\D_{jk}(0)=
2\sum_{j,k=1}^\n{b_{jk}\/\sqrt{\vk_j\vk_k}}\,,
$$
which yields \er{eq.7}.
In the proof of ii) we show that the estimates \er{eq.7} are sharp.

ii) Consider the graph $\G$
shown in Figure~\ref{ff.10}. Proposition \ref{TG1} gives
\[
\lb{sp}
|\s(\D)|=2\Big(1-\frac{\n-1}{\n+3}\,\Big)=\frac8{\n+3}\,.
\]
On the other hand, we estimate $|\s(\D)|$ using (\ref{eq.7}). For the graph $\G$ we have
$$
\vk_1=\ldots=\vk_{\n-1}=1,\qq \vk_\n=\n+3.
$$
The fundamental graph $\G_0$ has only 4 oriented bridges, which are the loops in the vertex $v_\n$. Thus,
$$
b_{\n\n}=4; \qqq b_{jk}=0, \qq \forall\, (j,k)\in\N^2_\n\setminus(\n,\n).
$$
Then the estimate (\ref{eq.7})
for the graph $\G$ has the form
\[
\lb{sp2}
|\s(\D)|\le\frac{8}{\nu+3}\,.
\]
Thus, \er{sp} and \er{sp2} show that for the graph $\G$ the estimate
(\ref{eq.7}) becomes an
identity. The last statement of the item is a direct consequence of \er{sp}.

iii) For each $\vt\in\T^2$ the matrix $\D(\vt)$ has $\nu$ eigenvalues,
where $\nu$ is odd. Theorem \ref{pro2}.vii gives that the spectrum $\s(\D(\vt))$ is
symmetric with respect to 0. Then $0\in\s\big(\D(\vt )\big)$ for
any $\vt\in\T^2$. Therefore, $\m=0$ is a flat band of $\D$. \qq \BBox

\

\no {\bf Proof of Theorem \ref{Tfb}.}
i) This follows from the facts that the point 1 is never a flat band (see Theorem \ref{pro2}.v) and $1\in\sigma(\Delta)$. The number of open spectral bands of the operator is $\n-r$. Some of them may overlap. Then the number of gaps between them is at most $\n-r-1$.

ii) This statement is a direct consequence of Propositions \ref{TG1}, \ref{TG2}, \ref{kag}. \qq \BBox

\section{Perturbations of square lattice}
\setcounter{equation}{0}

\subsection{Fundamental graphs with one vertex.} Let a fundamental
graph $\Gamma_0=(V_0,\cE_0)$ consist of one vertex $v$ and
any number of edges. We note that in this case the
index $\t(\be)$ of each edge $\be\in\cA$ coincides with the vector of its coordinates in the basis $a_1,a_2$ (the periods of~$\G$) and  all
edges of $\G_0$ are loops. We consider the spectrum of the Laplacian
$\D$ and describe all isospectral graphs (with one vertex in the fundamental graph) on which the spectrum
$\s(\D)=\s_{ac}(\D)=[-1,1]$.

Since $\n=1$, we deduce that $\Delta(\vt )$ is a scalar function given by
\[
\label{ll1}
\Delta(\vt )=\dfrac1{\vk_v}\sum_{\be\in\cA_0}\cos \lan\t
(\be),\vt\ran , \qquad \vt \in \T^2,
\]
where $\t(\be)=\big(\t_1(\be),\t_2(\be)\big)$ is the vector of the coordinates of $\bf e$ in the basis $a_1,a_2$. The spectrum of the operator $\D$ on the graph $\Gamma$ has the form
\begin{equation}\label{lap2}
\s(\D)=\s_{ac}(\D)=[\l^-,1], \qquad \l^-=\min\limits_{\vt
\in\T^2}\D(\vt )<1.
\end{equation}

\begin{theorem}\label{t1}
Let the fundamental graph $\Gamma_0=(V_0,\mathcal{E}_0)$ consist of one
vertex $v$ of the degree $\vk_v$. Then

i)  $\s(\D)=[-1,1]$ $\Leftrightarrow$ the graph $\Gamma$ is bipartite.

ii)  $\s(\D)=[-1,1]$ if  one of the following conditions holds true:

a) $\t_1(\be)$ is odd for all $\be\in\cA_0$;

b) $\t_2(\be)$ is odd for all $\be\in\cA_0$;

c) $\t_1(\be)+\t_2(\be)$ is odd for all
$\be\in\cA_0$.
\end{theorem}

\no \textbf{Proof.} i)  The spectrum of the Laplacian consists of one spectral band. The point $-1\in\s(\D)$ iff the graph is bipartite (see the main property 4) of the Laplacian). It gives the required statement.

Item ii) can be proved using Theorem \ref{T100} and the formula \er{l-}. But we give another proof (by contradiction). A graph is bipartite iff there are no cycles of odd length in it (see p.105 in \cite{Or62}). Let one of the conditions a) -- c) hold true and let the graph $\G$ be non-bipartite. Then there is a cycle with edges $\be_1,\ldots,\be_J\in\cA$ of odd length $J$ in it. It gives the identities
$$
\t_s(\be_1)+\ldots+\t_s(\be_J)=0, \qquad s=1,2.
$$
But it contradicts all conditions a) -- c), because the sum of an odd number $J$ of odd terms is not 0. Thus, the graph is bipartite and item i) gives $\s(\D)=[-1,1]$.
\quad $\BBox$

\medskip

We consider the spectrum of the Laplacian on
the square lattice $\mathbf{S}=(V,\mathcal{E})$, where the vertex set and
the edge set are given by
$$
V=\Z^2,\qquad  \cE=\big\{(p,p+e_1),(p,p+e_2) \quad
\forall\,p\in\Z^2\big\},
$$
the orthonormal basis $e_1,e_2$ coincides with the periods $a_1,a_2$ of $\mathbf{S}$,
see Figure \ref{ff.0.1}\emph{a}. The fundamental graph $\mathbf{S}_0$ of the
the square lattice $\mathbf{S}$ consists of one vertex $v$ and two unoriented
edges-loops $\mathbf{e}_1=\mathbf{e}_2=(v,v)$, see Figure
\ref{ff.0.1}\emph{b}. It is known that the spectrum of the Laplacian $\Delta$ on $\mathbf{S}$ has the form $
\s(\D)=\s_{ac}(\D)=[-1,1]$.
\setlength{\unitlength}{1.0mm}
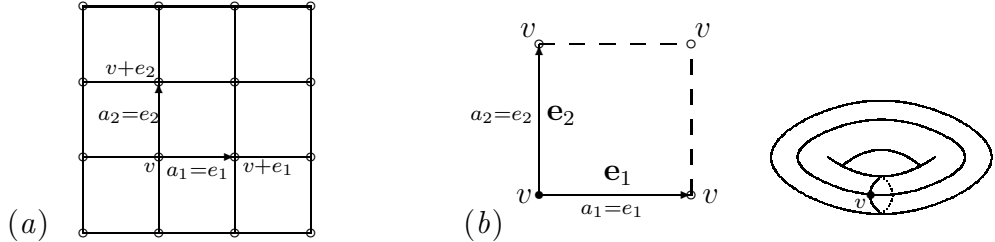
\begin{figure}[h]
\centering

\unitlength 1mm 
\linethickness{0.4pt}
\ifx\plotpoint\undefined\newsavebox{\plotpoint}\fi 
\begin{picture}(140,45)(0,0)

\put(10,10){\line(1,0){30.00}}
\put(10,20){\line(1,0){30.00}}
\put(10,30){\line(1,0){30.00}}
\put(10,40){\line(1,0){30.00}}
\put(10,10){\line(0,1){30.00}}
\put(20,10){\line(0,1){30.00}}
\put(30,10){\line(0,1){30.00}}
\put(40,10){\line(0,1){30.00}}

\put(10,10){\circle{1}}
\put(20,10){\circle{1}}
\put(30,10){\circle{1}}
\put(40,10){\circle{1}}

\put(10,20){\circle{1}}
\put(20,20){\circle{1}}
\put(30,20){\circle{1}}
\put(40,20){\circle{1}}

\put(10,30){\circle{1}}
\put(20,30){\circle{1}}
\put(30,30){\circle{1}}
\put(40,30){\circle{1}}

\put(10,40){\circle{1}}
\put(20,40){\circle{1}}
\put(30,40){\circle{1}}
\put(40,40){\circle{1}}

\put(20,20){\vector(1,0){10.00}}
\put(20,20){\vector(0,1){10.00}}

\put(18.0,18.0){$\scriptstyle v$}
\put(13,31){$\scriptstyle v+e_2$}
\put(31,18){$\scriptstyle v+e_1$}
\put(21,17.5){$\scriptstyle a_1=e_1$}
\put(12,25){$\scriptstyle a_2=e_2$}

\put(0,10){(\emph{a})}

\put(70,15){\vector(1,0){20.00}}
\put(70,15){\vector(0,1){20.00}}
\put(70,15){\circle*{1}}
\put(70,35){\circle{1}}
\put(90,35){\circle{1}}
\put(90,15){\circle{1}}

\multiput(70,35)(4,0){5}{\line(1,0){2}}
\multiput(90,15)(0,4){5}{\line(0,1){2}}

\put(67.0,14.0){$v$}
\put(91.5,14.0){$v$}
\put(67.5,36.0){$v$}
\put(90.5,36.0){$v$}
\put(75.5,12.5){$\scriptstyle a_1=e_1$}
\put(61.0,25.0){$\scriptstyle a_2=e_2$}
\put(78.5,17.0){$\mathbf{e}_1$}
\put(71.0,25.0){$\mathbf{e}_2$}


\bezier{120}(105,15)(115,10)(125,15)
\bezier{90}(105,15)(96,20)(105,25)
\bezier{90}(125,15)(134,20)(125,25)
\bezier{120}(105,25)(115,30)(125,25)

\bezier{120}(108,20)(115,15)(122,20)
\bezier{120}(110,19)(115,23)(120,19)
\put(113.6,15){\circle*{1}}
\bezier{15}(115,12.5)(118,15)(115,17.5)
\bezier{70}(115,12.5)(112,15)(115,17.5)

\bezier{120}(106,17.5)(115,12.5)(124,17.5)
\bezier{90}(106,17.5)(102,20)(106,22.5)
\bezier{90}(124,17.5)(128,20)(124,22.5)
\bezier{120}(106,22.5)(115,27.5)(124,22.5)
\put(111.5,13.3){$\scriptstyle v$}

\put(60,10){(\emph{b})}
\end{picture}

\vspace{-0.5cm} \caption{\emph{a}) The square lattice
$\mathbf{S}$;\quad \emph{b})  the fundamental graph $\mathbf{S}_0$.}
\label{ff.0.1}
\end{figure}

We describe the spectrum
\textbf{under the perturbation} of the graph $\bS$ by adding some edges
to its fundamental graph.

\begin{theorem}
{\bf (Perturbations of the square lattice).}
\label{t1'}
Let $\bS'=(\Z^2,\cE')$ be a  perturbed graph
obtained  from the square lattice $\bS$ by adding $N$ unoriented
edges to its fundamental graph $\bS_0$ and let
$\cA'_0$ be the set of all oriented edges of the fundamental graph of
$\bS'$. Then the spectrum $\s(\D')$  of the Laplacian $\D'$ on
the perturbed graph $\bS'$ satisfies:

i) $\s(\D')=[-1,1] \quad \Lra  \quad \t_1(\be)+\t_2(\be) \text{ is odd for all  } \be\in\cA'_0$, where $(\t_1(\be),\t_2(\be))$ is the vector of the coordinates of the edge $\be$.

ii) Let we add one oriented edge, i.e., two oriented
edges $\be,\bar\be$ and let $\t=(\t_1,\t_2)$ be the vector of the coordinates of the edge $\be$, where $\t_1+\t_2$ is even. Then $\s(\D')=[\l^-(\t),1]$,
where $\l^-(\t)$ satisfies
\[
\lb{ffff}
-1<\l^-(\t)\leq-\cos\dfrac{\pi}{|\t_1|+|\t_2|+1}\leq-\dfrac12\,.
\]
Moreover, the following asymptotics holds true:
\begin{equation}\label{as}
\l^-(\t)=-1+\frac{\pi^2}{6\,|\t|^2}+{O(1)\/|\t|^4}\quad
\textrm{as} \quad |\t|\rightarrow\infty.
\end{equation}
In particular, if $\t_1=\t_2$, then
\[
\lb{f1} \l^-(\t)=\frac13\min\limits_{\varphi\in[0,\pi]}
(2\cos\varphi+\cos2\t_1\varphi\,).
\]
\end{theorem}

\no\textbf{Proof.}
i) Let the edge set $\cA'_0$ consists of edges $\be_1,\be_2,\bar\be_1,\bar\be_2$ of the fundamental graph of the square lattice and $N$ additional unoriented edges, i.e., $2N$ oriented edges each of which also has an odd sum of the coordinates. Then by Theorem \ref{t1}.ii (the condition c), $\s(\D')=[-1,1]$.

In order to prove the converse we will use the proof by contradiction. Let
$\t_1(\be)+\t_2(\be)$ be even for some edge $\be\in\cA'_0$. Then there exists the cycle
$$
\underbrace{\be_1,\ldots,\be_1,}_{\t_{1}(\be)\textrm{ times}}\underbrace{\be_2,\ldots,\be_2,}_{\t_{2}(\be)\textrm{ times}}\bar\be
$$
of odd length in the graph $\bS'$. Thus, the graph is non-bipartite and Theorem \ref{t1}.i yields $\s(\D')\neq[-1,1]$.

ii) Without loss of generality we may assume that $0\leq \t_1\le \t_2$.
From \er{ll1} and the fact that $\vk_v=6$ we deduce that
\[
\lb{Del}
\D'(\vt)=\dfrac13
\big(\cos\vt_1+\cos\vt_2+\cos(\t_1\vt_1+\t_2\vt_2)\big).
\]
We will show \er{ffff}. Using \er{Del} we have
\[
\lb{lam}
\l^-(\t)=\min\limits_{\vt\in\T^2}\D'(\vt)
\leq\min\limits_{\vt\in\T^2:\atop\vt_1=\vt
_2=\vp}\D'(\vt)=\frac13\min\limits_{\vp\in[-\pi,\pi]}
\big(2\cos\vp+\cos(\t_1+\t_2)\vp\,\big)\leq\cos\vp_0,
\]
where $\vp_0$ is a solution of the equation
\[
\lb{co} \cos\vp_0=\cos(\t_1+\t_2)\vp_0,\qqq \vp_0\in \R.
\]
The solutions of this equation have the form
$$
\vp_0=\frac{2\pi n}{\t_1+\t_2\pm 1}\,,\qqq
n\in\Z.
$$
We take the solution $\vp_0$ given by
$$
\vp_0=\frac{\pi(\t_1+\t_2)}{\t_1+\t_2+1}=\pi-\frac{\pi}{\t_1+\t_2+1}\,,
$$
i.e, the nearest to $\pi$.
Then this identity and \er{lam} give
$$
\l^-(\t)\leq\cos\vp_0=
-\cos\frac{\pi}{\t_1+\t_2+1}\leq-\frac12\,.
$$
Thus, the inequality (\ref{ffff}) has been proved.

If $\t_1=\t_2=0$, then \er{Del} gives
\[
\lb{zer}
\l^-(0)=\min\limits_{\vt\in\T^2}\D'(\vt)=-\,\frac13\,.
\]
Let now $\t_1=\t_2\neq0$. Differentiating $\D'(\vt)$ given by \er{Del} at $\t_1=\t_2$ we obtain the necessary conditions for a minimum of $\D'$:
\[
\lb{nes}
\left\{
\begin{array}{l}
\sin \t_1(\vt_1+\vt_2)=\,\dfrac{\sin\vt_2}{\t_1}\,,\\[6pt]
\sin\vt_1=\sin\vt_2\,.
\end{array}\right.
\]
Using the second condition in \er{nes} we obtain two cases.

Firstly, if
$\vt_1=\vt_2$, then the function $\D'(\vt)$ has the form
\begin{equation}
\label{w1}
\D'(\vt )=\frac13\big(2\cos\vt_1+\cos 2\t_1\vt_1\,\big).
\end{equation}
 Secondly, if $\vt_1+\vt_2=\pi$, then $\D'(\vt
)=\frac{(-1)^{\t_1}}3$\,. But inequality (\ref{ffff}) gives that
$\min\limits_{\vt \in\T^2}\D'(\vt )\leq-\frac12$\,. Therefore the
minimum point of the function $\Delta'(\vt )$ is on the line $\vt
_1=\vt _2$ and thus,
$$
\l^-(\t)=\min\limits_{\vt\in\T^2,\ \vt_1=\vt
_2}\D'(\vt)=\frac13\min\limits_{\vt_1\in[0,\pi]}(2\cos\vt_1+\cos2\t_1\vt_1\,).
$$
This and \er{zer} yield \er{f1}.

We determine the asymptotics \er{as}. We have $\t_2\to\infty$ as
$|\t|\to\infty$, since $0\leq \t_1\leq \t_2$. Note that $\D'(-\vt)=\D'(\vt)$. Hence
\[
\lb{om0}
\l^-(\t)=\min\limits_{\vt\in\T^2}\D'(\vt)=\min\limits_{\vt\in\Theta}\D'(\vt),\qqq \Theta=\{-\pi\le\vt_1\le\pi,\ \vt_1\le\vt_2\le\pi\}.
\]
We introduce the local coordinates $\ve$ (see Figure \ref{Om}) by
$$
\vt=(\pi,\pi)-\ve, \qqq \ve=(\ve_1,\ve_2)\in\tilde{\Om}=
\{ 0\le\ve_1\le2\pi,\ 0\le\ve_2\le\ve_1\}.
$$
Then, using that $\t_1+\t_2$ is even, from \er{Del} we obtain
\[
\lb{a2}
\D'(\vt
)=\D'(\vt,\t)=\frac13\big(-\cos\ve_1-\cos\ve_2+\cos\b(\ve)\big),
\]
where
$$
\b(\ve)=\t_1\ve_1+\t_2\ve_2.
$$

Firstly, we show that the function $\D'(\vt,\t)$ achieves its global minimum on the torus $\T^2$ in the domain (see Figure \ref{Om})
\[
\lb{om5}
\Om=\{ \ve=(\ve_1,\ve_2)\in [0,\pi]^2: 0\le \b(\ve)\le\pi\le\pi\}.
\]

\setlength{\unitlength}{1.0mm}
\begin{figure}[h]
\centering
\unitlength 1mm 
\linethickness{0.4pt}
\ifx\plotpoint\undefined\newsavebox{\plotpoint}\fi 
\begin{picture}(50,55)(0,0)

\linethickness{0.1pt}
\put(12,10){\line(0,1){8.5}}
\put(14.2,10){\line(0,1){7.5}}
\put(16.0,10){\line(0,1){6.5}}
\put(18.0,10){\line(0,1){5.5}}
\put(20.2,10){\line(0,1){4.5}}
\put(22.0,10){\line(0,1){3.5}}
\put(24.0,10){\line(0,1){2.5}}
\put(26.0,10){\line(0,1){1.5}}
\put(28.0,10){\line(0,1){0.5}}

\linethickness{0.4pt}
\put(10,10){\vector(1,0){47.00}}
\put(10,10){\vector(0,1){47.00}}

\put(10,50){\line(1,-1){40.00}}
\put(10,30){\line(1,0){20.00}}
\put(30,10){\line(0,1){20.00}}

\bezier{80}(5,30)(25,20)(45,10)
\put(5,22){\line(2,-1){24.0}}
\bezier{40}(5,18)(13,14)(21,10)

\put(20,22.5){\circle{1}}
\put(13.1,14){\circle{1}}
\bezier{100}(10,10)(15,16.25)(20,22.5)

\put(8.0,8.0){$\scriptstyle 0$}
\put(-6.0,23.0){$\scriptstyle \b(\ve)=\pi$}
\put(-6.0,15.5){$\scriptstyle \b(\ve)=\b^0$}
\put(-6.0,30.0){$\scriptstyle \b(\ve)=\tilde\b$}

\put(18.0,24.0){$\scriptstyle (\tilde\ve_1,\,\tilde\ve_2)$}
\put(12.0,11.0){$\scriptstyle A$}

\put(30.0,8.0){$\scriptstyle \pi$}
\put(7.5,29.5){$\scriptstyle \pi$}
\put(6.0,49.0){$\scriptstyle 2\pi$}
\put(5.5,56.0){$\ve_2$}
\put(48.0,7.5){$\scriptstyle 2\pi$}
\put(56.0,6.8){$\ve_1$}

\put(15.0,35){$\tilde{\Om}$}
\put(18.0,12){$\scriptstyle \Om$}
\end{picture}
\vspace{-0.5cm}
\caption{The local coordinates $(\ve_1,\ve_2)$ belong to $\tilde{\Om}$; the shaded domain $\Om$ contains a global minimum point of the function $\D'(\vt,\t)$. The point $A=(\ve_1^0,\ve_2^0)$.} \label{Om}
\end{figure}
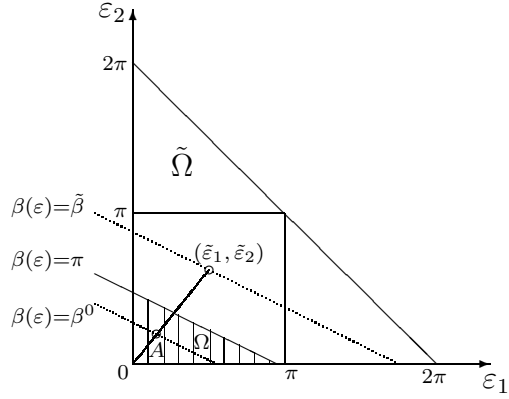

The function $\D'(\vt,\t)$ at the minimum point $(\ve_1,\ve_2)$ satisfies
$$
\left\{
\begin{array}{l}
\dfrac{\partial\D'}{\partial\ve_1}=
\dfrac13\big(\sin\ve_1-\t_1\sin\b(\ve)\big)=0,\\[12pt]
\dfrac{\partial\D'}{\partial\ve_2}=
\dfrac13\big(\sin\ve_2-\t_2\sin\b(\ve)\big)=0,
\end{array}\right.
$$
that is equivalent to
\[
\lb{a1}
\dfrac{\sin\ve_1}{\t_1}=\dfrac{\sin\ve_2}{\t_2}=\sin\b(\ve).
\]

From the fact that $0\leq \t_1\leq \t_2$ and the first identity in (\ref{a1}) it follows that $\sin\ve_1$ and $\sin\ve_2$ at the minimum point $(\ve_1,\ve_2)$ have the same sign. Thus, the minimum point $(\ve_1,\ve_2)\in[0,\pi]^2$.

We show that in the domain $\Om$ there exists a global minimum point of the function $\D'(\vt,\t)$. Let $\D'(\vt,\t)$ achieve its global minimum at a point $\tilde\ve=(\tilde\ve_1,\tilde\ve_2)\in[0,\pi]^2\sm\Om$ (see Figure \ref{Om}).
Show that there exists a point $\ve^0\in\Om$ such that
\[
\lb{mi}
\D'(\ve^0,\t)\leq\D'(\tilde\ve,\t).
\]

We denote
\[
\lb{om1}
\tilde\b=\b(\tilde\ve),\qqq \b^0=\arccos(\cos\tilde\b)\in [0,\pi],
\]
where the branch of $\arccos$ is fixed by the condition $\arccos(0)={\pi\/2}$.
Note that $\tilde\b>\pi$, since $\tilde\ve\in[0,\pi]^2\sm\Om$. Thus,
\[
\lb{om3}
q={\b^0\/\tilde\b}\in[0,1).
\]
We define the point $\ve^0$ (the point $A$ in Figure \ref{Om}) by
\[
\lb{om2}
\ve^0=\big(\ve_1^0,\ve_2^0\big)=q\,\tilde\ve.
\]
 Using \er{om1} -- \er{om2}, we have
$$
0\le\ve_s^0=q\,\tilde\ve_s<\tilde\ve_s\le\pi,\qq s=1,2,
$$
\[
\lb{om10}
\b(\ve^0)=q\,\b(\tilde\ve)=\b^0,
\]
and hence $\ve^0\in\Om$. Furthermore, the identity \er{a2} gives
\begin{multline*}
\D'\big(\ve^0,\t\big)=\frac13\big(-\cos\ve_1^0-\cos\ve_2^0+
\cos\b(\ve^0)\big)=\frac13\big(-\cos(q\,\tilde\ve_1)
-\cos(q\,\tilde\ve_2)+
\cos\b^0\big)\\
\le\frac13\big(-\cos\tilde\ve_1
-\cos\tilde\ve_2+
\cos\b^0\big)=\frac13\big(-\cos\tilde\ve_1
-\cos\tilde\ve_2+
\cos\tilde\b\,\big)=\D'\big(\tilde\ve,\t\big).
\end{multline*}
Here we have used the second identity in \er{om1}, \er{om10} and the following simple inequality
$$
\cos x\leq \cos (qx), \qqq\forall\,(x,q)\in[0,\pi]\ts [0,1],
$$
since $\cos $ is monotonic on the segment $[0,\pi]$. Thus, \er{mi} holds true,
i.e., the function $\D'(\vt,\t)$ also achieves its global minimum at the point $\ve^0\in\Om$. Then we rewrite \er{om0} as
\[
\lb{om4}
\l^-(\t)=\min\limits_{\ve\in\Om}\D'(\vt,\t).
\]

Secondly, we show that at the global minimum point $\ve=(\ve_1,\ve_2)\in\Om$ of the function $\D'(\vt,\t)$
\[
\lb{om6}
(\ve_1, \ve_2)={O(1)\/|\t|}\qqq  \qqq \textrm{as}\qq
|\t|\to\infty.
\]
The minimum point $\ve\in\Omega$, hence
$$
\b(\ve)\leq\pi,\qqq 0\le\ve_s,\qqq s=1,2.
$$
Since $0\leq \t_1\leq \t_2$, the last inequalities give
$$
\t_2\ve_2 \le \t_1\ve_1+\t_2\ve_2=\b(\ve)\leq\pi,\qqq 0\le\ve_2\le\frac{\pi}{\t_2}\,,
$$
which yields the asymptotics for the second component in (\ref{om6}). From this asymptotics and the first identity in (\ref{a1}) it follows that
$$
\sin\ve_1=\dfrac{\t_1}{\t_2}\,\sin\ve_2={O(1)\/|\t|},
$$
and hence $\ve_1\to0$ or $\ve_1\to\pi$ as $|\t|\to\infty$. But in the second case
$$
\l^-(\t)=\frac13\big(-\cos\ve_1-\cos\ve_2+\cos(\t_1\ve_1+\t_2\ve_2)\big)\to
\frac{(-1)^{\t_1}}3\cos \t_2\ve_2\qqq \textrm{as}\qq
|\t|\to\infty.
$$
Since this contradicts the estimate \er{ffff}, we conclude that $\ve_1\to0$ as $|\t|\to\infty$ and the asymptotics for the first component in (\ref{om6}) holds true. Thus, (\ref{om6}) has been proved.

Thirdly, we obtain the asymptotics for $\b(\ve)$ at the global minimum point $\ve$ as $|\t|\to\infty$.
The second identity in (\ref{a1}) and the asymptotics for the second component in (\ref{om6}) give
\[
\lb{a4} \sin\b(\ve)={O(1)\/|\t|^2}\,,
\]
which yields
$$
\b(\ve)=\t_1\ve_1+\t_2\ve_2=\pi n(\t)+{O(1)\/|\t|^2}, \qqq {\rm for  \ some } \
n(\t)\in \Z.
$$
Since the minimum point $\ve\in\Omega$, $n(\t)=0$ or $n(\t)=1$. If $n(\t)=0$, then
$$
\l^-(\t)=\frac13\big(-\cos\ve_1-\cos\ve_2+\cos(\t_1\ve_1+\t_2\ve_2)\big)\to
-\,\frac13\qqq \textrm{as}\qq
|\t|\to\infty.
$$
This again contradicts the estimate \er{ffff}. Thus, the function $\D'(\vt,\t)$ achieves its global minimum on the curve
\[
\lb{om7}
\b(\ve)=\t_1\ve_1+\t_2\ve_2=\pi+{O(1)\/|\t|^2}\,.
\]
This curve comes arbitrarily close to the line $\t_1\ve_1+\t_2\ve_2=\pi$ when $|\t|$ is rather large.

Finally, we obtain the asymptotics for the value of $\D'(\vt,\t)$ at the minimum point $\ve=(\ve_1,\ve_2)$ as $|\t|\to\infty$. The identity \er{om7} yields
\[
\lb{ww.3}
\cos(\t_1\ve_1+\t_2\ve_2)=\cos\bigg(\pi +{O(1)\/|\t|^2}\bigg)=-1+{O(1)\/|\t|^4}\,.
\]
Using the asymptotics \er{om6} and \er{ww.3} we rewrite the identity (\ref{a2}) as
\[
\lb{ww.5} \D'(\vt,\t)=
\frac13\bigg(-1+\frac{\ve_1^2}2-1+\frac{\ve_2^2}2-1+{O(1)\/|\t|^4}\bigg)=
-1+{\ve_1^2+\ve_2^2\/6}+{O(1)\/|\t|^4}\,.
\]

This function achieves the global minimum at the point $\ve=(\ve_1,\ve_2)$ of the curve \er{om7}, such that the square of the distance $\ve_1^2+\ve_2^2$ from the point $(0,0)$ to $\ve$ is minimal. Therefore, $\ve_1^2+\ve_2^2$ is equal to
the distance from the point $(0,0)$ to the curve \er{om7}.
The distance $d$ from a point $(x_0,y_0)\in\R^2$ to a line $Ax+By+C=0$, $A,B,C\in\R$, is given by
\[
\lb{dis}
d=\frac{|Ax_0+By_0+C|}{\sqrt{A^2+B^2}}\,.
\]
Using this formula we obtain
\[
\lb{ww.4}
\ve_1^2+\ve_2^2=\bigg(\frac{\pi}{|\t|}+{O(1)\/|\t|^4}\bigg)^2=
\frac{\pi^2}{|\t|^2}+{O(1)\/|\t|^5}\,.
\]
Substituting (\ref{ww.4}) into (\ref{ww.5}), we have
$$
\l^-(\t)=\min\limits_{\ve\in\Om}\D'(\vt,\t)=-1+\frac{\pi^2}{6\,|\t|^2}+{O(1)\/|\t|^4}\,.
$$
Thus, the asymptotics \er{as} has been proved.\qq
$\BBox$

\setlength{\unitlength}{1.0mm}
\begin{figure}[h]
\centering
\unitlength 1mm 
\linethickness{0.4pt}
\ifx\plotpoint\undefined\newsavebox{\plotpoint}\fi 
\begin{picture}(100,45)(0,0)

\put(10,10){\line(1,0){30.00}}
\put(10,20){\line(1,0){30.00}}
\put(10,30){\line(1,0){30.00}}
\put(10,40){\line(1,0){30.00}}
\put(10,10){\line(0,1){30.00}}
\put(20,10){\line(0,1){30.00}}
\put(30,10){\line(0,1){30.00}}
\put(40,10){\line(0,1){30.00}}

\put(10,10){\circle{1}}
\put(20,10){\circle{1}}
\put(30,10){\circle{1}}
\put(40,10){\circle{1}}

\put(10,20){\circle{1}}
\put(20,20){\circle*{1}}
\put(30,20){\circle{1}}
\put(40,20){\circle{1}}

\put(10,30){\circle{1}}
\put(20,30){\circle{1}}
\put(30,30){\circle{1}}
\put(40,30){\circle{1}}

\put(10,40){\circle{1}}
\put(20,40){\circle{1}}
\put(30,40){\circle{1}}
\put(40,40){\circle{1}}

\put(10,10){\line(2,1){30.00}}
\put(20,10){\line(2,1){20.00}}
\put(30,10){\line(2,1){10.00}}
\put(10,20){\line(2,1){30.00}}
\put(10,15){\line(2,1){30.00}}
\put(10,30){\line(2,1){10.00}}
\put(10,25){\line(2,1){30.00}}
\put(10,30){\line(2,1){20.00}}
\put(10,35){\line(2,1){10.00}}

\linethickness{1.2pt}
\put(20,20){\vector(1,0){10.00}}
\put(20,20){\vector(0,1){10.00}}
\put(20,20.2){\line(2,1){20.00}}
\put(20,20){\line(2,1){20.00}}
\put(20,19.8){\line(2,1){20.00}}
\linethickness{0.4pt}
\put(18.0,17.5){$\scriptstyle v$}
\put(22.5,18.0){$\scriptstyle e_1$}
\put(17.0,27.0){$\scriptstyle e_2$}
\put(0,10){(\emph{a})}
\put(70,10){\line(1,0){30.00}}
\put(70,20){\line(1,0){30.00}}
\put(70,30){\line(1,0){30.00}}
\put(70,40){\line(1,0){30.00}}
\put(70,10){\line(0,1){30.00}}
\put(80,10){\line(0,1){30.00}}
\put(90,10){\line(0,1){30.00}}
\put(100,10){\line(0,1){30.00}}

\put(70,10){\line(1,1){30.00}}
\put(70,20){\line(1,1){20.00}}
\put(80,10){\line(1,1){20.00}}
\put(70,30){\line(1,1){10.00}}
\put(90,10){\line(1,1){10.00}}

\put(70,10){\circle{1}}
\put(80,10){\circle{1}}
\put(90,10){\circle{1}}
\put(100,10){\circle{1}}

\put(70,20){\circle{1.0}}
\put(80,20){\circle*{1.0}}
\put(90,20){\circle{1.0}}
\put(100,20){\circle{1}}

\put(70,30){\circle{1}}
\put(80,30){\circle{1.0}}
\put(90,30){\circle{1.0}}
\put(100,30){\circle{1}}

\put(70,40){\circle{1}}
\put(80,40){\circle{1}}
\put(90,40){\circle{1}}
\put(100,40){\circle{1}}

\linethickness{1.2pt}
\put(80,20){\vector(1,0){10.00}}
\put(80,20){\vector(0,1){10.00}}
\put(80,20.2){\line(1,1){10.00}}
\put(80,20){\line(1,1){10.00}}
\put(80,19.8){\line(1,1){10.00}}
\linethickness{0.4pt}
\put(80.5,18.0){$\scriptstyle v$}
\put(84.5,18.0){$\scriptstyle e_1$}
\put(77.0,25.0){$\scriptstyle e_2$}
\put(60,10){(\emph{b})}
\end{picture}

\vspace{-0.5cm} \caption{Graph $\mathbf{S}'$, the edges of its
fundamental graph $\mathbf{S}'_0$ are marked by bold lines.}
\label{f.0.2}
\end{figure}
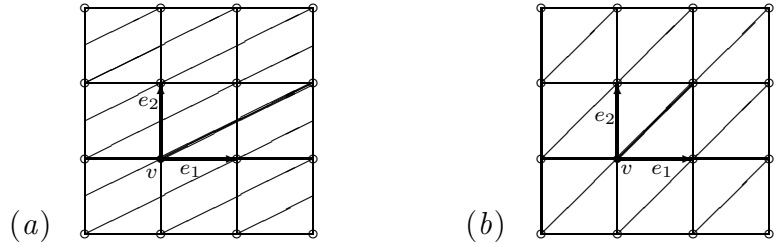

\no \textbf{Remark.}
1) Item i) describes all iso-spectral  perturbations.

2) In ii) we  estimate the end point $\l^-(\t)$ of the spectrum of the operator on the perturbed  square lattice. In the case when the added edges has equal coordinates we determine $\l^-(\t)$.

3) The asymptotics (\ref{as}) shows that if we add a "very long"\, edge to the fundamental graph of the square lattice, then the spectrum of the Laplacian is almost unchanged.

4) For the graph $\mathbf{S}'$  (on Figure~\ref{f.0.2}\emph{a})
$\sigma(\Delta')=\s_{ac}(\D')=[-1,1]$, since $\mathbf{S}'$ is obtained from the
square lattice $\mathbf{S}$ by adding the loop with the coordinates
$(2,1)$ (and its inverse loop) to the fundamental graph $\mathbf{S}_0$, $2+1=3$ is odd.

The
graph $\mathbf{S}'$, shown on Figure~\ref{f.0.2}\emph{b}, is
obtained from $\mathbf{S}$ by adding the loop with the coordinates $(1,1)$
to $\mathbf{S}_0$. Then Theorem~\ref{t1'}.ii yields
$$
\l^-(\t)=\frac13\min\limits_{\vp\in[0,\pi]}(2\cos\vp+\cos2\vp\,)=
-{1\/2}\qqq {\rm and} \qqq \sigma(\Delta')=\s_{ac}(\D')=[-1/2,1].
$$

\medskip

\section{Perturbations of Hexagonal lattice}
\setcounter{equation}{0}

\medskip

\subsection{Fundamental graphs with two vertices.}
We consider a bipartite periodic graph $\Gamma$. By Lemma \ref{l1}.ii, there exists a bipartite fundamental graph
$\Gamma_0=(V_0,\mathcal{E}_0)$. Assume that $\Gamma_0$ consists of two vertices $v_1$, $v_2$ and any number of edges. In this case the vertices $v_1$ and $v_2$ has the same degree $\vk$ and there are no loops in $\Gamma_0$.

For each $\vt\in \T^2$ the $2\times2$ matrix $\D(\vt)$ defined by (\ref{l2.15}) is
given by
\begin{equation}\label{l.82''}
\D(\cdot)=\ma
 0 & \D_{12}\\
  \bar\D_{12} & 0\\
\am, \qquad \D_{12}(\vt )=\dfrac1{\vk}
\sum\limits_{\mathbf{e}=(v_1,\,v_2)\in \cA_0}e^{i\lan\t
(\mathbf{e}),\vt\ran }.
\end{equation}
The diagonal entries of $\D(\vt)$ are zeroes, since there are no loops on $\G_0$.
For each $\vt \in\T^2$ the eigenvalues $\l_{1}(\vt ),\l_{2}(\vt )$ of
the matrix $\Delta(\vt )$ have the form
$$
\l_{1}(\vt )=-\l_{2}(\vt )=|\D_{12}(\vt)|.
$$
Recall that the spectrum  $\s(\D)$ of the Laplace operator on a bipartite graph with two vertices in the fundamental graph is symmetric with respect to 0 and consists of two spectral bands. Moreover, the point $1\in\s(\D)$. Thus we obtain
\[
\label{t2}
\begin{aligned}
&i)\qqq  \qqq \s(\D)=\s_{ac}(\D)=[-1,-\l_0]\cup[\l_0,1],\qq {\rm where} \qq
\l_0=\min\limits_{\vt
\in\T^2}\big|\D_{12}(\vt)\big|,\\
&ii)\qqq \qqq  \s(\D)=[-1,1]\ \Lra \ \D_{12}(\vt_0 )=0,\qq {\rm for\ some}
\ \vt_0\in\T^2.
\end{aligned}
\]

\medskip

The following statement gives the method for constructing of periodic graphs with $\l_0>0$.
\begin{proposition}\label{sug6}
Let the bipartite fundamental graph $\G_0$ consist of two vertices
$v_1,v_2$ and $N^2$ multiple oriented edges $(v_1,v_2)$ with indices
running over all values in the set
$$
\mathcal{D}=\{(\t_{1j},\t_{2j})\in\mathbb{Z}^2: j=1,\ldots,N\}
$$
and their inverse edges.
Then

i) The function $\D_{12}$ defined by (\ref{l.82''}) has  the form
$$
\D_{12}(\vt )=\frac1{N^2}\,P_1(e^{i\vt _1})P_2(e^{i\vt _2}),
$$
where
$$
P_s(z)=\sum\limits_{j=1}^{N}z^{\t_{sj}}, \qquad |z|=1,\qquad s=1,2.
$$

ii) Let $P_s(z)\neq0$ for each $(z,s)\in \{z\in \C: |z|=1\}\ts
\{1,2\}$. Then
 the spectrum $\s(\D)$ is symmetric with respect  to 0
and consists of exactly two spectral bands separated by a gap
$(-\l_0,\l_0)$, where $\l_0=\min\limits_{\vt
\in\T^2}\big|\D_{12}(\vt )\big|>0$.
\end{proposition}
\no \textbf{Proof.} i) The vertices $v_1$ and $v_2$ has the same degree $\vk=N^2$. For the function $\D_{12}$ defined by (\ref{l.82''}) we have
\begin{multline*}
\D_{12}(\vt)=\frac1\vk\sum\limits_{\mathbf{e}=(v_1,\,v_2)\in
\cA_0}e^{i\lan\t (\mathbf{e}),\vt\ran
}=\frac1{N^2}\sum\limits_{(\t_{1},\t_{2})\in\mathcal{D}}e^{i(\t_1\vt _1+\t_2\vt
_2)}\\=\frac1{N^2}\sum\limits_{j=1}^Ne^{i\t_{1j}\vt
_1}\sum\limits_{k=1}^Ne^{i\t_{2k}\vt _2}=\frac1{N^2}\,P_1(e^{i\vt
_1})P_2(e^{i\vt _2}).
\end{multline*}

ii) Since the functions $P_1,P_2$ have no zeroes, $\D_{12}(\vt )=\frac1{N^2}\, P_1(e^{i\vt_1})P_2(e^{i\vt_2})\neq0$ for any
$\vt\in\T^2$. Due to   \er{t2}, the
spectrum $\s(\D)$ is symmetric with
respect to 0 and consists of exactly two spectral bands separated by a
gap $(-\l_0,\l_0)$, where $\l_0=\min\limits_{\vt
\in\T^2}\big|\D_{12}(\vt )\big|>0$. \qq \BBox

\setlength{\unitlength}{1.2mm}
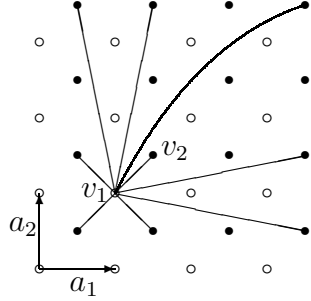
\begin{figure}[h]
\centering
\unitlength 1mm 
\linethickness{0.4pt}
\ifx\plotpoint\undefined\newsavebox{\plotpoint}\fi 
\begin{picture}(46,40)(0,0)

\put(0,0){\vector(0,1){10.0}}
\put(0,0){\vector(1,0){10.0}}

\put(0,0){\circle{1}}
\put(10,0){\circle{1}}
\put(20,0){\circle{1}}
\put(30,0){\circle{1}}

\put(0,10){\circle{1}}
\put(10,10){\circle{1}}
\put(20,10){\circle{1}}
\put(30,10){\circle{1}}

\put(0,20){\circle{1}}
\put(10,20){\circle{1}}
\put(20,20){\circle{1}}
\put(30,20){\circle{1}}

\put(0,30){\circle{1}}
\put(10,30){\circle{1}}
\put(20,30){\circle{1}}
\put(30,30){\circle{1}}

\put(5,5){\circle*{1}}
\put(15,5){\circle*{1}}
\put(25,5){\circle*{1}}
\put(35,5){\circle*{1}}

\put(5,15){\circle*{1}}
\put(15,15){\circle*{1}}
\put(25,15){\circle*{1}}
\put(35,15){\circle*{1}}

\put(5,25){\circle*{1}}
\put(15,25){\circle*{1}}
\put(25,25){\circle*{1}}
\put(35,25){\circle*{1}}

\put(5,35){\circle*{1}}
\put(15,35){\circle*{1}}
\put(25,35){\circle*{1}}
\put(35,35){\circle*{1}}

\put(10,10){\line(-1,-1){5.0}}
\put(10,10){\line(1,-1){5.0}}
\put(10,10){\line(-1,1){5.0}}
\put(10,10){\line(1,1){5.0}}

\put(10,10){\line(5,-1){25.0}}
\put(10,10){\line(5,1){25.0}}
\put(10,10){\line(-1,5){5.0}}
\put(10,10){\line(1,5){5.0}}

\qbezier(10.0,10.0)(20,30)(35.0,35.0)

\put(5.5,9.5){$v_1$}
\put(16.0,15.0){$v_2$}

\put(4.0,-3.0){$a_1$}
\put(-4.0,5.0){$a_2$}
\end{picture}
\caption{A graph $\Gamma$; only edges of the fundamental
graph $\Gamma_0$ are shown.} \label{fig3}
\end{figure}

\no \textbf{Example of a graph such that $\s(\D)$ consists of two spectral bands separated by a
gap.} The functions
$$
P_1(z)=P_2(z)=z^{-1}+1+z^2
$$
have no zeroes on the unit circle. Let the bipartite fundamental graph $\G_0$ consist of two vertices $v_1$, $v_2$ and $3^2=9$ oriented edges $(v_1,v_2)$ with the indices
$$
(0,0),\quad (-1,0),\quad (2,0), \quad (0,-1),\quad (-1,-1),\quad
(2,-1), \quad (0,2),\quad (-1,2),\quad (2,2)
$$
and their inverse edges (see Figure \ref{fig3}). Then by Proposition \ref{sug6}, the function $\D_{12}$ defined by (\ref{l.82''}) has the form
$$
\D_{12}(\vt )=\frac1{N^2}\, P_1(e^{i\vt _1})P_2(e^{i\vt _2})=
\frac1{N^2}\,\big(1+e^{-i\vt_1}+e^{2i\vt_1}\big)\big(1+e^{-i\vt_2}+e^{2i\vt_2}\big),
$$
and the spectrum $\sigma(\Delta)$ on the graph $\Gamma$ is symmetric with respect to 0 and consists of exactly two spectral bands
separated by a gap
$(-\l_0,\l_0)$, where $\l_0>0$. Direct calculations yield
$$
\l_0=\min\limits_{\vt
\in\T^2}\big|\D_{12}(\vt )\big|=\big|\D_{12}(\vt_1^\ast,\vt_2^\ast)\big|\approx0{,}04;
$$
where $\vt_s^\ast\in[0,\pi]$ is defined by $\cos\vt_s^\ast=\frac{-1+\sqrt{7}}6$\,, $s=1,2$.

\label{Gr}
\subsection{Hexagonal lattice.} We discuss the Laplacian on the
hexagonal lattice.

\begin{proposition}
\lb{TD}
The Laplacian on the hexagonal lattice $\textbf{G}$ satisfies
\[
\lb{Di}
\D(\vt)=H_D(t)+O(|t|^2) \qqq
\textrm{as}\qq |t|\to 0,
\]
\[
\lb{Di1}
t=(t_1,t_2)\in \R^2,\qqq
t_1=-\,\frac16\,(\vt_1+\vt_2),\qqq
t_2=\frac{\sqrt{3}}6\,\Big(-\vt_1+\vt_2+\frac{4\pi}3\Big),
\]
where $H_D(t)$ is the $2D$ Dirac operator given by
$$
H_D(t)=\s_1t_1+\s_2t_2,\qqq \s_1=\ma 0 & 1 \\
         1 & 0 \am, \qqq \s_2=\ma 0 & -i \\
         i & 0 \am.
$$
\end{proposition}
\no {\bf Proof}. The Floquet matrix $\Delta(\vt)$ for the hexagonal
lattice has the form
$$
\Delta(\vt)=\left(
\begin{array}{cc}
 0 & \D_{12}(\vt )\\[10pt]
 \bar{\D}_{12}(\vt ) & 0 \\
\end{array}\right),\qq
\D_{12}(\vt )=\frac13\,(1+e^{i\vt _1}+e^{i\vt _2}),\qq \vt=(\vt_1,\vt_2).
$$
It is easy to show that
$$
\D_{12}(\vt)=0 \qq
\Leftrightarrow \qq \vt=\pm\vt^0, \qqq \vt^0=(\vt^0_1,\vt^0_2)= \Big(\frac{2\pi}3\,,-\frac{2\pi}3\Big)\in\T^2.
$$
The Taylor expansion for the entry $\D_{12}(\vt)$ about the point $\vt^0$ is given by
\[
\lb{DDD}
\D_{12}(\vt)=\frac13\,(1+e^{i\vt _1}+e^{i\vt _2})=\frac13\big(1+e^{i\vt^0_1}(1+\vt_1-\vt^0_1)+
e^{i\vt^0_2}(1+\vt_2-\vt^0_2)\big)+O(|\vt-\vt^0|^2).
\]
Using the identity
$
e^{\,\pm i\frac{2\pi}3}={1\/2}(-1\pm i\sqrt{3}),
$
we rewrite \er{DDD} in the form
\[
\lb{ex}
\D_{12}(\vt)=\frac13+\frac16\,\Big[(-1+ i\sqrt{3}\;)(1+\vt_1-\vt^0_1)-(1+ i\sqrt{3}\;)(1+\vt_2-\vt^0_2)\Big]+O(|\vt-\vt^0|^2)=
t_1-it_2+O(|t|^2).
\]
Thus,
$$
\D(\vt)=\ma 0 & t_1-it_2 \\
         t_1+it_2 & 0 \am+O(|t|^2),
$$
which yields \er{Di}. Finally, we note that the Taylor expansion for $\D_{12}(\vt)$ about the point $-\vt^0$ is given by the same asymptotics \er{ex}, but $t_2$ is defined by $t_2=-\vt_1+\vt_2-{4\pi\/3}$\,.
\qq  $\BBox$

The result similar to the asymptotics \er{Di} for the 2D hexagonal lattice was described earlier by Wallace \cite{W47} and Slonczewski-Weiss [SW58]. Namely, they predicted that such a 2D monolayer material should present the branching points of the electron spectrum ($K$-points), where a degeneracy of the valence and conductivity band states takes place, the electronic state dispersion law asymptotically has a form of the double-napped cone, while the quasi-wave vector approaches the $K$-point. The equation set for the electron states in the vicinity of the $K$-point is mathematically similar to the Dirac equation for a zero-mass particle. These features stem from the specific symmetry of the hexagon lattice, which has two Bravais sublattices, and all atoms of it are situated in identical positions. These properties of this 2D material do not depend on approximate procedures used in calculation of the electronic spectrum.
The mentioned above similarity to the Dirac equation was used by many authors in the study of the graphene electronic properties: bound and resonance states, electron scattering, conductivity and other transport coefficients, see \cite{FKP09}, \cite{FK10} and references therein.

\subsection{Perturbed Hexagonal lattice.}
In order to prove Theorem \ref{T10} we consider the function $\cF$ given by
\[
\lb{FF1}
\cF(\vt)=\big|1+e^{i\vt_1}+e^{i\vt_2}\big|^2,\qqq \vt=(\vt_1,\vt_2)\in\T^2.
\]
We have the following
\[
\lb{FF2}
\min\limits_{\vt\in\T^2}\cF(\vt)=\cF(\pm \vt^0)=0,\qqq \textrm{where} \qq \vt^0=(2\pi/3,-2\pi/3).
\]
We rewrite $\cF$ in the form
$$
\cF(\vt)=1+8\cos\frac{\vt_1-\vt_2}2\cos\frac{\vt_1}2\cos\frac{\vt_2}2\,,
$$
which yields
\[
\lb{FF3}
\max\limits_{\vt\in\T^2}\cF(\vt)=\cF(0)=9.
\]

\no\textbf{Proof of Theorem \ref{T10}.}
Item i) is a direct consequence of ii) and iii).

We show ii). Recall that the fundamental graph $\bG_0$ of the hexagonal lattice $\bG$ consists of two vertices $v_1$,
$v_2$, three multiple oriented edges
$\mathbf{e}_1=\mathbf{e}_2=\mathbf{e}_3=(v_1,v_2)$ (Figure
\ref{ff.0.3}\emph{b}) with the indices  $\t (\mathbf{e}_1)=(0,0)$,
$\t (\mathbf{e}_2)=(1,0)$, $\t (\mathbf{e}_3)=(0,1)$ and their inverse edges.

The graph $\bG'$ remains bipartite iff the adding edge $\be$ connects the vertices $v_1$ and $v_2$. Let we add the edge $\be=(v_1,v_2)$ with an index
$\t=(\t_1,\t_2)\in\Z^2$ (and its inverse edge) to $\bG_0$.
For each $\vt=(\vt_1,\vt_2)\in \T^2$ the matrix $\D(\vt)$ defined by (\ref{l2.15}) has the
form
$$
\Delta(\vt )=\left(
\begin{array}{cc}
 0 & \D_{12}(\vt )\\[10pt]
 \bar{\D}_{12}(\vt ) & 0 \\
\end{array}\right),\qqq
\D_{12}(\vt )=\frac14\,(1+e^{i\vt _1}+e^{i\vt _2}+e^{i\lan\t,\vt\ran}).
$$
We define the point
$$
\vt_0=\left\{
\begin{array}{cl}
  (0,\pi), & {\rm if \qq } \t_2 {\rm \; is \; odd}\\[4pt]
  (\pi,0), & {\rm if \qq } \t_1 {\rm \; is \; odd}\\[4pt]
  (\pi,\pi), & {\rm otherwise }
  \end{array}\right..
$$
One can verify by a direct calculation that $\D_{12}(\vt_0)=0$. From ii) in \er{t2}, it follows that
$\s(\D')=[-1,1]$. Thus, we have proved that if $\bG'$ is bipartite then $\s(\D')=[-1,1]$. The converse follows from item iii) of this theorem.

iii) The graph $\bG'$ is non-bipartite
iff we add a loop with an index
$\t=(\t_1,\t_2)\in\Z^2$ (and its inverse loop) to $\bG_0$.
Without loss of generality we may assume that we add a loop $\mathbf{e}=(v_2,v_2)$.
For the graph $\bG'$ the matrix $\D(\vt)$ has the
form
$$
\Delta(\vt )=\left(
\begin{array}{cc}
 0 & \D_{12}(\vt )\\[10pt]
  \bar\D_{12}(\vt ) & \D_{22}(\vt )\\
\end{array}\right),
$$
where
$$
\D_{12}(\vt )=\dfrac1{\sqrt{15}}\,(1+e^{i\vt _1}+e^{i\vt _2}),\qquad \D_{22}(\vt )=\dfrac2{5}\cos\lan \t,\vt\ran,\qquad \vt \in\T^2.
$$
The eigenvalues of the matrix $\Delta(\vt )$ are given by
\[
\lb{ge.1}
\l_{s}(\vt )=\frac15\,\cos\lan \t,\vt\ran+(-1)^{s+1}\sqrt{\frac1{25}\,\cos^2\lan \t,\vt\ran+\frac1{15}\,\cF(\vt)}\,,\qqq s=1,2,
\]
where $\cF(\vt)$ is defined by \er{FF1}. Thus, the spectrum of the Laplacian on the non-bipartite perturbed graphene has the form
$$
\s(\D)=[\l_2^-,\l_2^+]\cup[\l_1^-,\l_1^+],\qqq
\l_s^-=\min\limits_{\vt\in\T^2}\l_s(\vt),\qqq \l_s^+=\max\limits_{\vt\in\T^2}\l_s(\vt),\qqq s=1,2.
$$
Using \er{ge.1} and  \er{FF2}, \er{FF3}, we have
$$
\l_1^+=\max\limits_{\vt\in\T^2}\l_1(\vt)=\l_1(0)=1,\qqq
\l_2^-=\min\limits_{\vt\in\T^2}\l_2(\vt)
\leq\l_2(0)=-{3\/5}\,.
$$
The last inequality and non-bipartition of the graph $\bG'$ yield \er{ge1}.

We will show \er{ge2}, \er{ge3}. From (\ref{ge.1}) and \er{FF2}, \er{FF3},
it follows that
\[
\lb{ge.4}
\l_1(\vt)\ge0,\qqq \l_2(\vt)\le0, \qqq \forall\,\vt\in\T^2,
\]
and
\[
\lb{ge.3}
\l_s\big(\vt^0\big)
=\frac15\,\cos\Big({2\pi\/3}\,(\t_1-\t_2)\Big)
+\frac{(-1)^{s+1}}{5}\,\Big|\cos\Big({2\pi\/3}\,(\t_1-\t_2)\Big)\Big|,\qq s=1,2,
\]
where $\vt^0=({2\pi\/3}\,,-{2\pi\/3})$.
If $\t_1-\t_2\in3\Z$,  then  \er{ge.4}, \er{ge.3} give
$$
\l_1^-=\min\limits_{\vt\in\T^2}\l_1(\vt)
\leq\l_1\big(\vt^0\big)=\frac25\,,\qqq 0\ge\l_2^+=\max\limits_{\vt\in\T^2}\l_2(\vt)\ge\l_2\big(\vt^0\big)=0.
$$
If $\t_1-\t_2\notin3\Z$, then
$$
0\le\l_1^-=\min\limits_{\vt\in\T^2}\l_1(\vt)\le\l_1\big(\vt^0\big)=0,\qqq \l_2^+=\max\limits_{\vt\in\T^2}\l_2(\vt)
\ge\l_2\big(\vt^0\big)=-\frac15\,.
$$
This yields \er{ge2}, \er{ge3}.

The estimates \er{ge4} follow directly from \er{ge2} -- \er{ge1}.
\qq $\BBox$

\no \textbf{Example of a non-bipartite perturbed graphene.} Let the graph
$\mathbf{G}'$ be obtained from the hexagonal lattice $\mathbf{G}$ by
adding the loop ${\bf e}=(v_2,v_2)$ with the index $\t({\bf e})=(\t_1,\t_2)=(1,0)$ (and its inverse loop) to the fundamental
graph $\mathbf{G}_0$ (Figure~\ref{f.0.5}).

\setlength{\unitlength}{1.1mm}
\begin{figure}[h]
\centering
\unitlength 1mm 
\linethickness{0.4pt}
\ifx\plotpoint\undefined\newsavebox{\plotpoint}\fi 
\begin{picture}(60,42)(0,0)

\put(14,10){\circle{1}}
\put(28,10){\circle{1}}
\put(34,10){\circle{1}}
\put(48,10){\circle{1}}

\put(18,16){\circle{1}}
\put(24,16){\circle*{1}}
\put(38,16){\circle{1}}
\put(44,16){\circle{1}}

\put(14,22){\circle{1}}
\put(28,22){\circle*{1}}
\put(34,22){\circle{1}}
\put(48,22){\circle{1}}

\put(18,28){\circle{1}}
\put(24,28){\circle{1}}
\put(38,28){\circle{1}}
\put(44,28){\circle{1}}

\put(14,34){\circle{1}}
\put(28,34){\circle{1}}
\put(34,34){\circle{1}}
\put(48,34){\circle{1}}

\put(18,40){\circle{1}}
\put(24,40){\circle{1}}
\put(38,40){\circle{1}}
\put(44,40){\circle{1}}
\qbezier(14,10)(29,19)(44,28)
\qbezier(34,10)(39,13)(44,16)
\qbezier(14,22)(29,31)(44,40)
\qbezier(14,34)(19,37)(24,40)
\put(28,10){\line(1,0){6.00}}
\put(18,16){\line(1,0){6.00}}
\put(38,16){\line(1,0){6.00}}

\put(28,22){\line(1,0){6.00}}
\put(18,28){\line(1,0){6.00}}
\put(38,28){\line(1,0){6.00}}

\put(28,34){\line(1,0){6.00}}
\put(18,40){\line(1,0){6.00}}
\put(38,40){\line(1,0){6.00}}

\put(14,10){\line(2,3){4.00}}
\put(34,10){\line(2,3){4.00}}
\put(24,16){\line(2,3){4.00}}
\put(44,16){\line(2,3){4.00}}

\put(14,22){\line(2,3){4.00}}
\put(34,22){\line(2,3){4.00}}
\put(24,28){\line(2,3){4.00}}
\put(44,28){\line(2,3){4.00}}

\put(14,34){\line(2,3){4.00}}
\put(34,34){\line(2,3){4.00}}

\put(28,10){\line(-2,3){4.00}}
\put(48,10){\line(-2,3){4.00}}
\put(38,16){\line(-2,3){4.00}}
\put(18,16){\line(-2,3){4.00}}

\put(28,22){\line(-2,3){4.00}}
\put(48,22){\line(-2,3){4.00}}
\put(38,28){\line(-2,3){4.00}}
\put(18,28){\line(-2,3){4.00}}

\put(28,34){\line(-2,3){4.00}}
\put(48,34){\line(-2,3){4.00}}

\put(30,18){$\scriptstyle a_1$}
\put(20.5,22){$\scriptstyle a_2$}

\put(24,16){\vector(0,1){12.0}}
\put(33,21.3){\vector(3,2){0.5}}

\put(24.8,21.5){$\scriptstyle v_1$}
\put(35,21.0){$\scriptstyle v_2+a_1$}
\put(25.0,27.0){$\scriptstyle v_2+a_2$}
\put(25.0,15.0){$\scriptstyle O=v_2$}
\end{picture}
\vspace{-1cm} \caption{The non-bipartite perturbed graphene $\mathbf{G}'$.} \label{f.0.5}
\end{figure}
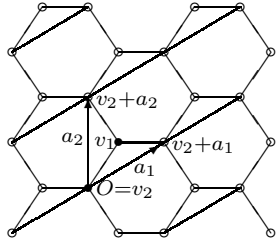

Since $\t_1-\t_2\notin3\Z$, by Theorem \ref{T10}, we have
$$
\s(\D')=\s_{ac}(\D')=[\l_2^-,\l_2^+]\cup[0,1].
$$
We obtain numerically that
$$
\l_2^-\approx-0{,}627;\qqq \l_2^+\approx-0{,}106.
$$

\subsection{Bipartite fundamental graphs with three vertices.}\lb{ss1}
Consider a bipartite graph $\G$. It has a bipartite fundamental graph $\Gamma_0=(V_0,\mathcal{E}_0)$ with parts $V_1$ and $V_2$. Let $\Gamma_0$ consist of three vertices $v_1$, $v_2$, $v_3$ with degrees $\vk_1$, $\vk_2$, $\vk_3$, respectively,
and any number of edges. Without loss of generality we may assume that $v_1\in V_1$ and $v_2,v_3\in V_2$.

For each $\vt\in \T^2$ the matrix $\Delta(\vt)$ defined by (\ref{l2.15}) is
the following $3\times3$ matrix
\begin{equation}\label{z2}
\Delta(\vt )=\left(
\begin{array}{ccc}
   0 & \D_{12}(\vt ) & \D_{13}(\vt ) \\
  \bar\D_{12}(\vt ) &  0 & 0 \\
  \bar \D_{13}(\vt ) & 0 & 0
\end{array}\right),
\end{equation}
where
\begin{equation}
\label{fff}
\D_{1s}(\vt )=
\dfrac1{\sqrt{\vk_1\vk_s}}\sum\limits_{\mathbf{e}=(v_1,\,v_s)\in\cA_0}e^{i\lan\t
(\mathbf{e}),\vt\ran },\qquad s=2,3.
\end{equation}

In the following proposition we determine the spectrum of the operator $\Delta$ on the graph $\Gamma$ and formulate the necessary and sufficient conditions when $\sigma(\Delta)=[-1,1]$.

\begin{proposition}
\label{t5}
Let the bipartite fundamental graph $\G_0$ consist of three vertices and some number of edges. Then
\[
\lb{all}
\begin{aligned}
&\s(\D)=\s_{ac}(\D)\cup\s_{fb}, \qqq \s_{ac}(\D)=[-1,-\l_0]\cup[\l_0,1], \qqq \s_{fb}=\{0\},\\
&\l_0=\min\limits_{\vt
\in\T^2}\sqrt{|\D_{13}(\vt )|^2+|\D_{12}(\vt )|^2};\\
&\sigma(\Delta)=[-1,1]\ \Lra \ \D_{12}(\vt_0)=\D_{13}(\vt_0)=0 \; \textrm{ for some }\; \vt_0\in\T^2,
\end{aligned}
\]
where the functions $\D_{12},\D_{13}$  are defined by (\ref{fff}).
\end{proposition}

\no\textbf{Proof.} For each $\vt \in\T^2$ the eigenvalues $\l_{1}(\vt ),\l_{2}(\vt ),\l_3(\vt )$ of
the matrix $\Delta(\vt )$ given by (\ref{z2}) have the form
$$
\l_{1}(\vt )=-\l_{2}(\vt
)=\sqrt{|\D_{13}(\vt )|^2+|\D_{12}(\vt )|^2}\,, \qqq \l_3(\cdot)=\m=0,
$$
which yields \er{all}.
\quad $\BBox$

\begin{corollary}\lb{TC1}
Let $\G$ be the graph obtained from the square lattice $\bS$ by adding one vertex on each edge of $\bS$. Then the spectrum of the Laplacian on $\G$ has the form
$$
\s(\D)=\s_{ac}(\D)\cup\s_{fb}(\D), \qqq
\s_{ac}(\D)=[-1,1], \qqq \s_{fb}(\D)=\{0\}.
$$
\end{corollary}
\no\textbf{Proof.} The graph $\G$
is a bipartite periodic graph with
three vertices on the fundamental graph. Then by Proposition \ref{t5} the spectrum has the form
$$
\s(\D)=\s_{ac}(\D)\cup\s_{fb},\qqq \s_{ac}(\D)=[-1,-\l_0]\cup[\l_0,1], \qqq \s_{fb}=\{0\},
$$
where
$\l_0=\min\limits_{\vt
\in\T^2}\sqrt{|\D_{13}(\vt )|^2+|\D_{12}(\vt )|^2}$. The functions
$\D_{12},\D_{13}$ defined by (\ref{fff}) have the form
$$
\D_{1s}(\vt)=\frac1{2\sqrt{2}}\,(1+e^{-i\vt_s}).
$$
Due to identities
$
\D_{12}(\vt_0)=\D_{13}(\vt_0)=0,\  \vt_0=(\pi,\pi),
$
Proposition \ref{t5} gives $\l_0=0$ and $\s_{ac}(\D)=[-1,1]$.\qq \BBox
\subsection{The Kagome lattice\lb{KL}.} \lb{kl} As an example of non-bipartite graph with three vertices in the fundamental graph consider
the graph $\G$ shown in Figure~\ref{f.0.6}\emph{a}. This graph is called the Kagome lattice. It is a lattice structure found in many natural minerals' molecular arrangements.

The fundamental graph of the Kagome lattice consists of three vertices $v_1$, $v_2$, $v_3$ each of which has the degree 4, six oriented edges
$$
\mathbf{e}_1=\mathbf{e}_2=(v_1,v_2),\qquad \mathbf{e}_3=\mathbf{e}_4=(v_1,v_3),\qquad \mathbf{e}_5=\mathbf{e}_6=(v_2,v_3)
$$
and their inverse edges.
The indices of the fundamental graph edges in the coordinate system with the origin $O$ (Figure~\ref{f.0.6}\emph{a}) are given by
\begin{multline*}
\t (\mathbf{e}_1)=(0,0),\quad \t (\mathbf{e}_2)=(-1,0),\quad \t (\mathbf{e}_3)=(0,0),\\[4pt]
\t (\mathbf{e}_4)=(0,-1),\quad \t (\mathbf{e}_5)=(0,0),\quad \t (\mathbf{e}_6)=(1,-1).
\end{multline*}

\setlength{\unitlength}{1.0mm}
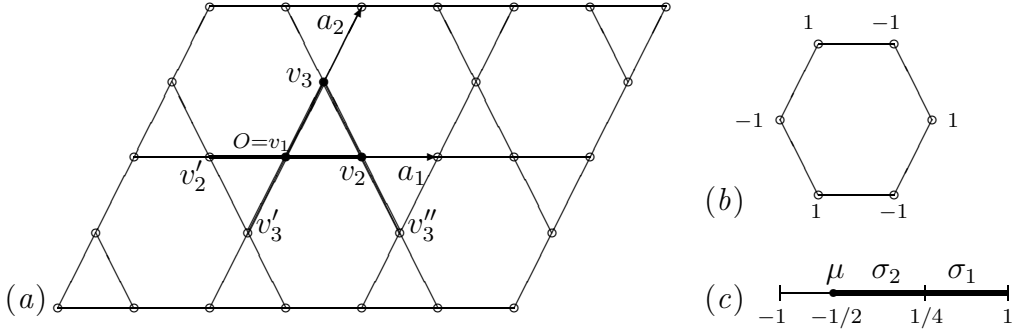
\begin{figure}[h]
\centering
\unitlength 1mm 
\linethickness{0.4pt}
\ifx\plotpoint\undefined\newsavebox{\plotpoint}\fi 
\begin{picture}(140,45)(0,0)

\put(0,0){\line(1,0){60.0}}
\put(10,20){\line(1,0){60.0}}
\put(20,40){\line(1,0){60.0}}
\put(0,0){\line(1,2){20.0}}
\put(20,0){\line(1,2){20.0}}
\put(40,0){\line(1,2){20.0}}
\put(60,0){\line(1,2){20.0}}

\put(10,0){\line(-1,2){5.0}}
\put(30,0){\line(-1,2){15.0}}
\put(50,0){\line(-1,2){20.0}}
\put(65,10){\line(-1,2){15.0}}
\put(75,30){\line(-1,2){5.0}}

\put(30,20){\vector(1,0){20.0}}
\put(30,20){\vector(1,2){10.0}}

\put(0,0){\circle{1}}
\put(10,0){\circle{1}}
\put(20,0){\circle{1}}
\put(30,0){\circle{1}}
\put(40,0){\circle{1}}
\put(50,0){\circle{1}}
\put(60,0){\circle{1}}

\put(5,10){\circle{1}}
\put(25,10){\circle{1}}
\put(45,10){\circle{1}}
\put(65,10){\circle{1}}

\put(10,20){\circle{1}}
\put(20,20){\circle{1}}
\put(30,20){\circle{1}}
\put(40,20){\circle{1}}
\put(50,20){\circle{1}}
\put(60,20){\circle{1}}
\put(70,20){\circle{1}}

\put(15,30){\circle{1}}
\put(35,30){\circle{1}}
\put(55,30){\circle{1}}
\put(75,30){\circle{1}}

\put(20,40){\circle{1}}
\put(30,40){\circle{1}}
\put(40,40){\circle{1}}
\put(50,40){\circle{1}}
\put(60,40){\circle{1}}
\put(70,40){\circle{1}}
\put(80,40){\circle{1}}

\put(30,20){\circle*{1}}
\put(40,20){\circle*{1}}
\put(35,30){\circle*{1}}

\linethickness{1.2pt}
\put(20,20){\line(1,0){20.0}}
\put(25,10.3){\line(1,2){10.0}}
\put(25,9.7){\line(1,2){10.0}}
\put(45,10.3){\line(-1,2){10.0}}
\put(45,9.7){\line(-1,2){10.0}}
\linethickness{0.4pt}

\put(23.0,21.2){${\scriptstyle O=v_1}$}
\put(37.0,17.0){$v_2$}
\put(16.0,16.5){$v'_2$}
\put(30.0,30.0){$v_3$}
\put(26.0,10.0){$v'_3$}
\put(46.0,10.0){$v''_3$}

\put(44.7,17.0){$a_1$}
\put(34.0,37.0){$a_2$}
\put(-7.0,0.0){(\emph{a})}
\put(85.0,13.0){(\emph{b})}

\put(100,15){\line(1,0){10.0}}
\put(100,35){\line(1,0){10.0}}
\put(100,15){\line(-1,2){5.0}}
\put(115,25){\line(-1,2){5.0}}
\put(115,25){\line(-1,-2){5.0}}
\put(100,35){\line(-1,-2){5.0}}

\put(100,15){\circle{1}}
\put(110,15){\circle{1}}
\put(100,35){\circle{1}}
\put(110,35){\circle{1}}
\put(115,25){\circle{1}}
\put(95,25){\circle{1}}

\put(99.0,12.0){$\scriptstyle1$}
\put(108.0,12.0){$\scriptstyle-1$}

\put(89.0,24.0){$\scriptstyle-1$}
\put(117.0,24.0){$\scriptstyle1$}

\put(98.0,37.0){$\scriptstyle1$}
\put(107.0,37.0){$\scriptstyle-1$}
\put(85.0,0){(\emph{c})}
\put(95,2){\line(1,0){30.00}}
\put(102,1.8){\line(1,0){23.00}}
\put(102,2.2){\line(1,0){23.00}}

\put(95,1){\line(0,1){2.00}}
\put(125,1){\line(0,1){2.00}}
\put(114,1){\line(0,1){2.00}}

\put(102,2){\circle*{1}}

\put(117,4){$\s_1$}
\put(107,4){$\s_2$}
\put(101,4){$\m$}
\put(92,-2){$\scriptstyle-1$}
\put(112.0,-2){$\scriptstyle1/4$}
\put(99.0,-2){$\scriptstyle-1/2$}
\put(124.2,-2){$\scriptstyle1$}
\end{picture}
\caption{(\emph{a}) The Kagome lattice, the edges of the fundamental graph are marked by bold lines;
(\emph{b}) the support of the eigenfunction;\quad \emph{c}) the spectrum of the Laplacian.} \label{f.0.6}
\end{figure}

\begin{proposition}\label{kag}
The spectrum of the Laplace operator on the Kagome lattice has the form
\begin{equation}\label{l.71}
\s(\D)=\s_{ac}(\D)\cup\s_{fb}(\D),\qqq\sigma_{ac}(\Delta)=[-1/2;1],\qqq \s_{fb}(\D)=\{-1/2\}.
\end{equation}
\end{proposition}

\no\textbf{Proof.} For each $\vt\in \T^2$ the matrix $\Delta(\vt)$ has the form
$$
\Delta(\vt )=\frac14\left(
\begin{array}{ccc}
   0 & 1+e^{-i\vt _1} & 1+e^{-i\vt _2} \\
  1+e^{i\vt _1} & 0 & 1+e^{i(\vt _1-\vt _2)}\\
  1+e^{-i\vt _2} & 1+e^{i(\vt _1-\vt _2)} & 0
\end{array}\right).
$$
The eigenvalues of this matrix are given by
$$
\l_1(\vt )=\frac14\big(1+\sqrt{\cF(\vt )}\,\big), \qqq
\l_2(\vt )=\frac14\big(1-\sqrt{\cF(\vt )}\,\big), \qqq
\l_3(\vt )=\mu=-\frac12\,,\qqq \forall\,\vt\in\T^2,
$$
where $\cF(\vt)$ is defined by \er{FF1}. From \er{FF2}, \er{FF3} it follows that
the ranges of the functions $\l_1(\vt)$
and $\l_2(\vt)$ are $[{1\/4};1]$ and $[-{1\/2};{1\/4}]$, respectively. Thus,
the spectrum of the operator has the form (\ref{l.71}) and
$\mu=-{1\/2}$ is a flat band. \quad $\BBox$

\no \textbf{Remark.}
1) The eigenfunctions corresponding to the eigenvalue
 $\mu_1=-1/2$ have the finite support shown in Figure \ref{f.0.6}\emph{b}.
 The values of the eigenfunction in the vertices of the support
 are pointed out.\\
2) As sets the spectra of the Laplacians on the Kagome lattice and on
the triangular lattice (Figure \ref{f.0.2}\emph{b}) are the same.

\section{Appendix}
\setcounter{equation}{0}

\subsection{Properties of  matrices.}
We recall some well-known  properties of matrices (see e.g., \cite{HJ85} and \cite{HC96}), which will be used below.
Let $\r(A)$ be the spectral radius of~$A$.

\begin{proposition}\lb{MP}
i) Let $A=\{A_{jk}\}$ and $B=\{B_{jk}\}$ be $\nu\ts\nu$ matrices. If ${|A_{jk}|\leq B_{jk}}$ for all $j,k\in\N_\n$, then $\r(A)\leq\r(B)$ (see Theorem 8.1.18 in \cite{HJ85}).


ii) Let $A,B$ be $\nu\ts\nu$ self-adjoint  matrices and let $B\geq0$. Then the eigenvalues $\l_n(A)\leq\l_n(A+B)$ for all $n\in\N_\n$ (see Corollary 4.3.3 in \cite{HJ85}).

iii) Let $B$ be the self-adjoint $(\nu+1)\ts(\nu+1)$ matrix given by
$$
B=\left(
\begin{array}{cc}
  A & y \\
  y^\ast & a
\end{array}\right)
$$
for some self-adjoint $\nu\ts \nu$ matrix $A$, some  real number $a$ and some vector $y\in\C^{\nu}$.

Let the eigenvalues of $A$ and $B$ be denoted by $\{\mu_j\}$ and
$\{\l_j\}$, respectively, and assume that they have been arranged in
decreasing order
$$
\mu_\nu\leq\ldots\leq\mu_1,\qquad
\l_{\nu+1}\leq\l_{\nu}\leq\ldots\leq\l_1.
$$
Then
$$
\l_{\nu+1}\leq\mu_\nu\leq\l_{\nu}\leq\ldots\leq\l_2\leq\mu_1\leq\l_1,
$$
(see Theorem 4.3.8 in \cite{HJ85}).

%

iv) Let $M$ be a $\n\ts\n$ matrix having the form
$$
M=\ma
  A & B \\
  C & D
\am
$$
for some square matrices $A,D$ and some matrices $B,C$. Then
\[
\lb{det}
\det M=\det A\cdot\det\big(D-CA^{-1}B\big)
\]
(see pp.21--22 in \cite{HJ85}).

v) Let $M_\n$ be a $\n\ts\n$ finite Jacobi  matrix given by
$$
M_\n=\left(
\begin{array}{cccc}
  b & 1 & 0 &\ldots \\
  1 & b & 1 &\ldots\\
  0 & 1 & b &\ldots \\
   \ldots & \ldots & \ldots &\ldots \\
\end{array}\right),
$$
where $b=2\cos\vp\in(-2,2)$. Then
\[
\lb{det1}
\det M_\n=(-1)^\n\,\frac{\sin(\n+1)\vp}{\sin\vp}\,,
\]
which satisfies
\[
\lb{rel}
\det M_{n+1}=b\det M_{n}-\det M_{n-1},\qqq \forall \ n\in\N,
\]
with the initial conditions $\det M_0=1$, $\det M_1=b$ (see pp.1511--1512 in \cite{HC96}).

Moreover, the eigenvalues of the matrix $M_\n$ have the form
\[
\lb{eig}
\l_n=b+2\cos{\pi n\/\n+1}\,, \qqq \forall \ n\in\N_\n.
\]

\end{proposition}

\

\subsection{The basic properties of periodic graphs and fundamental
graphs.}
In this subsection we discuss main properties of the periodic graphs.

\begin{proposition}\label{pro0}
The fundamental graph of a $\Z^2$-periodic graph is finite.
\end{proposition}

\no \textbf{Proof.}
Recall that we identify the vertices of the fundamental graph $\G_0$ with the vertices of the periodic graph from the bounded domain $[0,1)^2$. So the number of vertices of $\G_0$ is finite (see item 1 of the definition of a $\Z^2$-periodic graph).
The number of edges of $\G_0$ is also finite, because $\G_0$ contains a finite number of vertices, each of which has a finite degree (see item 2 of the definition). \quad $\BBox$
\vspace {3mm}

Below we need the following properties of an edge index.
\begin{proposition}\label{pro1}
Let $(u,v)\in\cA$ and let $p,q\in\Z^2$. Then the following identities hold true.

i) $\t(u,v)=-\t(v,u)$.

ii) Let $(v+p,v+q)\in\cA$. Then the index $\t(v+p,v+q)=q-p$.

iii) The edge $(u+p,v+p)\in\cA$ and its index $\t(u+p,v+p)=\t(u,v)$.

iv) Let $\t^{(1)}(\tilde\be)$ be the index of an edge $\tilde\be=(u,v)\in\cA_0$ in the coordinate system with the origin $O_1$. Then
\begin{equation}\label{ind'}
\t^{(1)}(\tilde\be)=\t(\tilde\be)+[v-b]-[u-b],\qqq\textrm{ where }\;
b=\overrightarrow{OO}_1.
\end{equation}
\end{proposition}

\no \textbf{Proof.}
Using the definition (\ref{in}) of an edge index, we have
$$
\begin{array}{ll}
{\rm i)} \quad &\t(u,v)=[v]-[u]=-([u]-[v])=-\t(v,u), \\
{\rm ii)}\quad&\t(v+p,v+q)=[v]+q-[v]-p=q-p. \\
\end{array}
$$

iii) The periodicity of the graph yields that the edge $(u+p,v+p)\in\cA$  and
$$
\t(u+p,v+p)=[v]+p-[u]-p=[v]-[u]=\t(u,v).
$$

\setlength{\unitlength}{1.0mm}
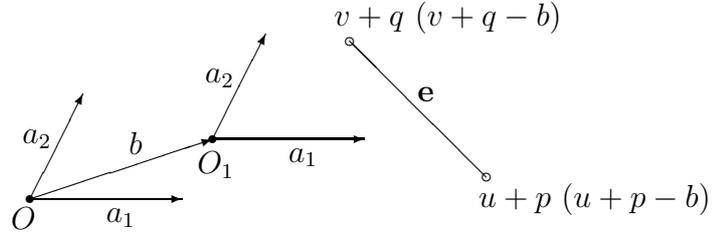
\begin{figure}[h]
\centering
\unitlength 1mm 
\linethickness{0.4pt}
\ifx\plotpoint\undefined\newsavebox{\plotpoint}\fi 
\begin{picture}(60,30)(0,0)

\put(10,5){\circle*{1}}
\put(7.5,1.0){$O$}
\put(20.0,2.0){$a_1$}
\put(9.0,12.5){$a_2$}
\put(10,5){\vector(1,0){20.00}}
\put(10,5){\vector(1,2){7.00}}

\put(10,5){\vector(3,1){24.00}}

\put(34,13){\circle*{1}}
\put(32.0,8.5){$O_1$}
\put(44.0,10.0){$a_1$}
\put(33.0,20.5){$a_2$}
\put(34,13){\vector(1,0){20.00}}
\put(34,13){\vector(1,2){7.00}}

\put(23.0,11.0){$b$}

\put(70,8){\line(-1,1){18.00}}
\put(70,8){\circle{1}}

\put(52,26){\circle{1}}
\put(50,28.0){$v+q$ ($v+q-b$)}
\put(61.0,18.0){$\be$}

\put(69,4.0){$u+p$ ($u+p-b$)}
\end{picture}

\caption{The edge $\be=(u+p,v+q)$ in the coordinate system with the origin~$O$; $\be=(u+p-b,v+q-b)$ in the coordinate system with the origin $O_1$; $a_1,a_2$ are the periods of the graph.} \label{gg.00}
\end{figure}

iv) Let $\tilde\be=(u,v)$ be an oriented edge of the fundamental graph with an index $\t(\tilde\be)$. Then by the definition of the fundamental graph and the formulas \er{in}, \er{inf} there is an edge $\be=(u+p,v+q)\in\cA$, where  $p,q\in\Z^2$ are some integer vectors such that
\[
\lb{in4}
\t(\be)=q-p=\t(\tilde\be)
\]
(Figure~\ref{gg.00}). Recall that we identify the vertices of the fundamental graph with the vertices of the periodic graph from the set $[0,1)^2$
in the coordinate system with the origin $O$.
In the coordinate system with the origin $O_1$ the edge $\be$ has the form
$$
\be=(u+p-b,v+q-b\big).
$$
Using \er{in}, \er{inf} and \er{in4} we obtain
$$
\t^{(1)}(\tilde\be)=\t^{(1)}(\be)=q+[v-b]-p-[u-b]=\t(\tilde{\bf e})+[v-b]-[u-b].
$$
Thus, the identity \er{ind'} has been proved.
\qq
$\BBox$

\setlength{\unitlength}{1.0mm}
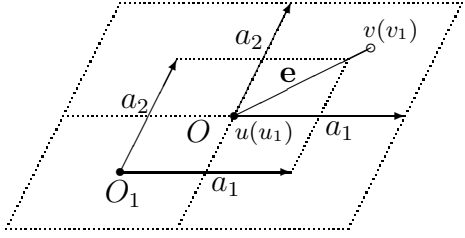
\begin{figure}[h]
\centering
\unitlength 1mm 
\linethickness{0.4pt}
\ifx\plotpoint\undefined\newsavebox{\plotpoint}\fi 
\begin{picture}(60,30)(0,0)

\bezier{70}(0,0)(22.5,0)(45,0)
\bezier{70}(7.5,15)(30,15)(52.5,15)
\bezier{70}(15,30)(37.5,30)(60,30)

\bezier{50}(0,0)(7.5,15)(15,30)
\bezier{50}(22.5,0)(30,15)(37.5,30)
\bezier{50}(45,0)(52.5,15)(60,30)

\put(30,15){\circle*{1}}
\put(23.7,11.5){$O$}
\put(42.0,12.5){$a_1$}
\put(30.1,24.5){$a_2$}
\put(30,15){\vector(1,0){22.5}}
\put(30,15){\vector(1,2){7.5}}

\put(30,15){\line(2,1){18.00}}

\put(30.1,12){\scriptsize $u$($u_1$)}
\put(48,24){\circle{1}}
\put(47,25.5){\scriptsize$v$($v_1$)}
\put(36,19.5){$\be$}

\put(15,7.5){\circle*{1}}
\put(13.0,3.5){$O_1$}
\put(27.0,5.0){$a_1$}
\put(15.1,16.0){$a_2$}
\put(15,7.5){\vector(1,0){22.5}}
\put(15,7.5){\vector(1,2){7.5}}

\bezier{40}(22.5,22.5)(33.75,22.5)(45,22.5)
\bezier{30}(37.5,7.5)(41.25,15)(45,22.5)

\end{picture}

\caption{The edge $\be=(u,v)$ in the coordinate system with the origin $O$; $\be=(u_1,v_1)$ in the coordinate system with the origin $O_1$.} \label{ff.00}
\end{figure}

\no \textbf{Remark.}
1) From item iv) it follows that an edge index, generally
speaking, depends on the choice of the coordinate origin.
Indeed, let the oriented
edge $\be=(u,v)$ in the coordinate system with the origin $O$
and the same edge $\be=(u_1,v_1)$
in the system with the origin $O_1$ (Figure~\ref{ff.00}). Then the
index of the edge $\be=(u,v)$ in the first system is defined by
$$
\t(\be)=[v]-[u]=(0,0)-(0,0)=(0,0)
$$
and the index of the same edge $\be=(u_1,v_1)$ in the second system has the form
$$
\t^{(1)}(\be)=[v_1]-[u_1]=(1,1)-(0,0)=(1,1).
$$
Item ii) shows that the index of the edge
$(v+p,v+q)\in\cA$ does not depend on the choice of the coordinate origin $O$. It also means that the indices of all loops
on the fundamental graph do not depend on the choice of the point
$O$.

2) Under the group $\Z^2$ action the set $\cA$ of oriented edges of the
graph  $\G$ is divided into equivalence classes. Each equivalence
class is an oriented edge $\be\in\cA_0$ of the fundamental graph $\G_0$. From item iii) it follows that all edges from one equivalence class $\mathbf{e}$
have the same index that is also the index of the fundamental graph
edge $\mathbf{e}$.

\

Now we show that a bipartite periodic graph has a
bipartite fundamental graph. Recall that a graph $\G=(V,\cE)$ is
\emph{bipartite} if its vertex  set  $V$ can be divided into two
disjoint sets $V_1$ and $V_2$ (called \emph{parts} of the graph) such that every
edge connects vertices from different sets.

\begin{lemma}\label{l1}
i) If $v$ belongs to some part  $V_1$ of the graph, then $v+2p\in V_1$ for each $p\in\Z^2$.

ii) For a bipartite periodic graph there exists a bipartite
fundamental graph.
\end{lemma}
\no \textbf{Proof.} i) The distance between two vertices in a graph
is the number of edges in the shortest path connecting these
vertices. It is known that the distance between two
vertices from one part of a bipartite graph is even and the distance
between two vertices from the different parts is odd (see p.105 in \cite{Or62}).

 For the vertex $v+p$ we have two cases.

 Firstly, let  $v+p\in V_1$. Then the distance $d(v,v+p)$ between $v$ and $v+p$ is even. Due to the periodicity of the graph the shortest path connecting $v+p$ and $v+2p$ is obtained by the
translation of the shortest path between $v$ and $v+p$ through the
vector $p$. Thus, the distance $d(v+p,v+2p)$ is also even and hence $v+2p\in V_1$.

Secondly, if $v+p\in V_2$,
then $d(v,v+p)$ is odd. The distance $d(v+p,v+2p)$ is also odd and hence $v+2p\in V_1$. Thus, in both case $v+2p\in V_1$.

ii) Let $\G_0=(V_0,\cE_0)$ be a fundamental graph of a bipartite
periodic graph $\G$. If it is non-bipartite, then we consider another
fundamental graph $\G'_0=(V'_0,\cE'_0)=\G/(2\Z)^2$.
We will show that $\G'_0$ is bipartite. We identify the vertices of $\G'_0$ with the vertices of the periodic graph $\G$ from the set $[0,2)^2$. Let us divide $V'_0$ into two
disjoint sets $V'_1$ and $V'_2$ in the following way:
\[
\lb{part}
v\in V'_s \ \Lra \ v\in V_s\qqq {\rm for \ each } \ s=1,2.
\]
We show that each edge of the fundamental graph $\G'_0$ connects vertices from the different sets $V'_1$ and $V'_2$. Assume the contrary, that
there exists an edge
$\be=(u,v)\in\cE'_0$ connecting the vertices $u,v$ from one set $V'_1$ or $V'_2$. Without loss of generality suppose that $u,v\in V'_1$.
Then by the definition of the fundamental graph $\G'_0$ there is an edge $(u+p,v+q)\in\cE$ for some $p,q\in2\Z$. Since $u,v\in V'_1$, by \er{part} and item i) we have
$u+p,v+q\in V_1$. Thus, the edge $(u+p,v+q)$ on $\G$ connects the vertices from one part of $\G$. This contradicts the bipartition of $\G$. Thus, there is no edge connecting vertices
from the same part of $\G'_0$ and hence $\G'_0$ is a bipartite
fundamental graph of $\G$. \qq $\BBox$

\medskip

\no\textbf{Acknowledgments.}
\footnotesize
Various parts of this paper were written during Evgeny Korotyaev's stay in the Mathematical Institute of Tsukuba University, Japan  and Mittag-Leffler Institute, Sweden. He is grateful to the institutes for the hospitality.
His study was supported by The Ministry of education and science of Russian Federation, project   07.09.2012  No  8501 No «2012-1.5-12-000-1003-016»
and the RFFI grant "Spectral and asymptotic methods
for studying of the differential operators" No 11-01-00458.

This work was also supported by the federal target program "Scientific and scientific-pedagogical personnel of innovative Russia" in 2009 -- 2013 of the Ministry of Education and Science of the Russian Federation, state contract 14.740.11.0581.


\begin{thebibliography}{9999}
\setlength{\itemsep}{-\parskip}
\footnotesize

\bibitem[A12]{A12} Ando, K. Inverse scattering theory for discrete Schr\"odinger operators on the hexagonal lattice, arXiv:1110.3922.


\bibitem[Ba98]{Ba98} Bollobas, B. Modern graph theory. Graduate Texts in Mathematics, 184. Springer-Verlag, New York, 1998. xiv+394 pp.

\bibitem[Bo08]{Bo08}    Bondy, J. A.; Murty, U. S. R. Graph theory. Graduate Texts in Mathematics, 244. Springer, New York, 2008. xii+651 pp.

\bibitem[BS99]{BS99}
A. Boutet de Monvel and J. Sahbani, On the spectral properties of discrete Schr\"odinger operators : (The multi-dimensional case), Review in Math. Phys., 11 (1999), 1061--1078.

\bibitem[BH12]{BH12} Brouwer, A.; Haemers, W. Spectra of graphs. Universitext. Springer, New York, 2012.

\bibitem[C97]{C97} Cattaneo, C. The spectrum of the continuous Laplacian
 on a graph, Monatsh. Math. 124 (1997), 215--235.

\bibitem[Ch97]{Ch97} Chung, F. Spectral graph theory, American Mathematical Society,
 Providence, Rhode. Island, 1997.

\bibitem[CDS95]{CDS95} Cvetkovic, D.; Doob, M.; Sachs, H. Spectra of graphs. Theory and applications. Third edition. Johann Ambrosius Barth, Heidelberg, 1995. ii+447 pp.

\bibitem[CDGT88]{CDGT88}
  Cvetkovic, D.; Doob, M.;  Gutman, I,; Torgaљev, A. Recent results in the theory of graph spectra. Annals of Discrete Mathematics, 36. North-Holland Publishing Co., Amsterdam, 1988.

\bibitem[FK10]{FK10} Firsova, N.; Ktitorov, S.
Electrons scattering in the monolayer graphene with short-range impurities. Physics Letters A. 374 (2010), 1270--1273.

\bibitem[FKP09]{FKP09} Firsova, N.; Ktitorov, S.; Pogorelov, P.
Bound electron states in the monolayer gapped graphene with short-range impurities. Physics Letters A. 373 (2009), 525--528.



\bibitem  [GKT93] {GKT93} Gieseker, D.; Kn\"orrer, H.; Trubowitz, E.
The geometry of algebraic Fermi curves. Perspectives in Mathematics,
 14. Academic Press, Inc., Boston, MA, 1993.

\bibitem[Ha85]{Ha85} Harris, P.  Carbon nano-tubes and related structure,
Cambridge, Cambridge University Press, 2002.

\bibitem[HC96]{HC96} Hu, G.; O'Connell, R. Analytical inversion of symmetric tridiagonal matrices. J. Phys. A: Math. Gen. 29 (1996), 1511--1513.


\bibitem[HS04]{HS04} Higuchi, Y. Shirai, T. Some spectral and geometric properties for infinite graphs, AMS Contemp. Math. 347 (2004), 29--56.

\bibitem[HN09]{HN09} Higuchi, Y.; Nomura, Y.
Spectral structure of the Laplacian on a covering graph.
European J. Combin. 30 (2009), no. 2, 570--585.

\bibitem[HJ85]{HJ85} Horn, R; Johnson, C. Matrix analysis.
Cambridge University Press, 1985.

\bibitem[IK12]{IK12} Isozaki, H. ;  Korotyaev, E. Inverse problems,
trace formulae for discrete Schr\"odinger operators,  Annales Henri
Poincare, 13(2012), No 4 ,  751--788.


\bibitem[KK10]{KK10}  Korotyaev, E.;  Kutsenko, A. Zigzag nanoribbons in
external electric fields, Asympt. Anal. 66 (2010), no. 3--4,
187--206.


\bibitem[KK10a]{KK10a} Korotyaev, E.; Kutsenko, A. Zigzag nanoribbons
 in external   electric and magnetic fields. Int. J. Comput. Sci. Math.  3 (2010),
no. 1--2, 168--191.

\bibitem[KK10b]{KK10b} Korotyaev, E.;  Kutsenko, A. Zigzag and
armchair nanotubes in external fields, "Differential Equations:
Advances in Mathematics Research, Volume 10 (2010) Nova Science
Publishers, Inc. 273--302.


\bibitem[KS]{KS} Korotyaev, E.; Saburova, N. Laplacians on periodic quantum graphs, in preparation.

\bibitem[KS1]{KS1}  Korotyaev, E.; Saburova, N. Schr\"odinger operators on
 $\Z^d$ periodic discrete graphs, in preparation.


\bibitem[Ku10]{Ku10} Kutsenko A. Sharp spectral estimates for periodic
matrix-valued Jacobi operators, preprint 2010.
http://arxiv.org/abs/1007.5412.


\bibitem[MRA07]{MRA07} Mantoiu, M.; Richard, S.; Tiedra de Aldecoa, R.
Spectral analysis  for adjacency operators on graphs.
  Ann. Henri Poincare' 8 (2007), no. 7, 1401--1423.

\bibitem[Me94]{Me94} Merris, R.  Laplacian matrices of graphs: a survey,
Linear algebra and its applications, 197-198(1994),143--176.



\bibitem[M92]{M92} B. Mohar, Laplace eigenvalues of
graphs: a survey, Discrete mathematics 109 (1992), 171--183.

\bibitem[MW89]{MW89} Mohar, B.; Woess, W. A survey on spectra of infinite graphs,
Bull. London Math. Soc., 21 (1989), 209--234.

\bibitem[NG04]{NG04} Novoselov K.S.; Geim A.K. et al,
Electric field effect in atomically thin carbon films, Science 22 October,
306(2004), no. 5696, 666--669.

\bibitem[Or62]{Or62} Ore,O. Theory of graphs. AMS Colloquium Publications 38. AMS 1962. 270p.

\bibitem[P12]{P12} Post, O. Spectral analysis on graph-like spaces. Lecture Notes in Mathematics, 2039. Springer, Heidelberg, 2012.

%

\bibitem[RR07]{RR07} Rabinovich, V. S.; Roch, S.
Essential spectra of difference operators on $Z^n$-periodic graphs.
J. Phys. A 40 (2007), no. 33, 10109--10128.

\bibitem[RS78]{RS78} Reed, M.; Simon, B. Methods of modern mathematical
physics, Vol.IV, Analysis of operators, Academic Press, New York, 1978.

\bibitem[RoS09]{RoS09} Rosenblum, G.; Solomjak, M. On the spectral estimates for the Schr{\"o}dinger operator on ${\bf Z}^d$, $d \geq 3$, Problems in Mathematical Analysis, No. 41, J. Math. Sci. N. Y. 159 (2009), No. 2, 241--263.

\bibitem[SDD98]{SDD98} Saito R, Dresselhaus, G.; Dresselhaus, M.
Physical properties of carbon nanotubes, London, Imperial College Press,  1998.

\bibitem[SW58]{SW58} Slonczewski, J.C.; Weiss, P.R. Band structure of
graphite, Phys. Rev., 109(1958), no 2, 272--279.



\bibitem[W47]{W47} Wallace, P. The band theory of
graphite, Phys. Rev., 71(1947), no 9, 622--634.



\end{thebibliography}
\end{document}